\documentclass[12pt,a4paper,reqno]{amsart}

\usepackage[utf8]{inputenc}
\usepackage[T1]{fontenc}
\usepackage[english,french]{babel}
\usepackage[margin=2.5cm]{geometry}

\usepackage{amsmath} 
\usepackage{amssymb}
\usepackage{amsthm} 
\usepackage{mathtools}
\usepackage{mathrsfs} 
\usepackage{tabularx}
\usepackage{yfonts}
\usepackage{enumerate} 
\usepackage{tikz}
\usepackage{stmaryrd} 
\allowdisplaybreaks
\usepackage{subcaption}
\usepackage{textcomp} 
\usepackage{gensymb}

\usepackage{empheq} 

\usepackage[pdfborder={0 0 0}]{hyperref}				
\usepackage{amsfonts} 
\usepackage{graphicx} 
\DeclareGraphicsExtensions{.eps,.pdf,.png,.jpg,.JPG,.jpeg,.ps,.svg}
\usepackage{color}
\usepackage{xcolor}
\usepackage{float}
\usepackage{blkarray}
\usepackage{mathrsfs} 
\usepackage{mathdots}
\usepackage{minitoc}
\usepackage{appendix}
\DeclareTextSymbol{\degre}{T1}{6}

\newcommand{\R}{\mathbb{R}}

\newcommand{\dsp}{\displaystyle}

\newcommand{\image}[3]  
{\begin{figure}[H]                                               %
	\begin{center}                                                     %
	\includegraphics[width=#3\columnwidth]{#1}                         %
	\caption{#2} \label{#1}                                        %
	\end{center}                                                       %
	\end{figure}                                                       %
}

\newcommand{\BB}[1]{\begin{equation}  \label{#1} \begin{aligned}}
\newcommand{\EE}{ \end{aligned} \end{equation}}

\newcommand{\BBl}[1]
{
\begin{equation} \label{#1}
\left|
 \begin{aligned}
 }
 
\newcommand{\EEl}{
 \end{aligned} 
 \right. 
 \end{equation}
 }
 
 \newcommand{\ald}{\[ \begin{aligned}}
   \newcommand{\alf}{ \end{aligned} \]}

\theoremstyle{plain}

\newtheorem{Proposition}{Proposition}

\newcommand{\edit}[1]{\textcolor{black}{#1}}
\newcommand{\vv}{{\mathtt{v}}}
\newcommand{\dvv}{{\partial_{n_y}\mathtt{v}}}
\newcommand{\dvvx}{{\partial_{n_{x^\ast}}\mathtt{v}}}
\begin{document}
\selectlanguage{english}

\title{Modified representations for the close evaluation problem}
\author{C. Carvalho\\ Applied Math Department, Univeristy of California Merced, USA}
\maketitle

\begin{abstract}
 When using boundary integral equation methods, we represent solutions of a linear partial differential equation as layer potentials. It is well-known that the approximation of layer potentials using quadrature rules suffer from poor resolution when evaluated closed to (but not on) the boundary. To address this challenge, we provide modified representations of the problem's solution. Similar to Gauss's law used to modify Laplace's double-layer potential, we use modified representations of Laplace's single-layer potential and Helmholtz layer potentials that avoid the close evaluation problem. {Some techniques have been developed in the context of the representation formula or using interpolation techniques\edit{. W}e provide alternative modified representations of the layer potentials directly (or when only one density is at stake).} Several numerical examples illustrate the efficiency of the technique in two and three dimensions.
\end{abstract}
\section{Introduction}
One can represent the solution of partial differential boundary-value problems using boundary integral equation methods, which involves integral operators defined on the domain's boundary called layer potentials. 
Using layer potentials, the solution can be evaluated anywhere in the domain without restriction to a particular mesh. For that reason boundary integral equations
have found broad applications, including in fluid mechanics, 
electromagnetics, and plasmonics~\cite{smith2009boundary,barnett2015spectrally,marple2016fast,keaveny2011applying,akselrod2014probing,novotny2011antennas,mayer2008label,sannomiya2008situ}.

The close evaluation problem refers to the nonuniform error produced
by high-order quadrature rules used to discretize layer potentials. This phenomenon arises when computing the solution close to the boundary (i.e. at close evaluation points). It is well understood that this growth in error is due to the fact that the integrands of the layer potentials become increasingly peaked as the evaluation point approaches the boundary (nearly singular behavior), leading in limit cases to an $O(1)$ error {\cite{barnett2014evaluation}}.

There exists a plethora of manners to address the close evaluation problem: using extraction methods based on Taylor series expansions~\cite{schwab1999extraction}, regularizing the nearly singular behavior of the integrand and adding corrections~\cite{beale2001method,beale2016simple}, compensating quadrature rules via interpolation~\cite{helsing2008evaluation}, using Quadrature By Expansion related techniques (QBX)~\cite{barnett2014evaluation,klockner2013quadrature,epstein2013convergence,af2016fast,rachh2017fast,af2017error,wala20183DQBX}, {using adaptive methods~\cite{GORV21},} using singularity subtraction techniques and interpolation~\cite{Perez18,PeTuFa19,PeFaTu19}, or using asymptotic approximations~\cite{carvalho2018asymptotic,ckk2020asymp,KKCC19}, to name a few. Most techniques rely on either providing corrections to the \textit{kernel} (related to the fundamental solution of the PDE at stake), or to the \textit{density} (solution of the boundary integral equation).

In the latter category, it is well-known that Laplace's double-layer potential can be straightforwardly modified via a \textit{density subtraction technique} \edit{based on Gauss' law (e.g.~\cite{Hwang13})}. This modification alleviates the close evaluation problem, and provides a better approximation for any given numerical method. However this identity technique is specific to Laplace's double-potential. Other identities have been derived for other problems, such as for the elastostatic problem \cite{Liu99}.

In this paper we provide modified representations of layer potentials, and \edit{we give} guidance to address the close evaluation problem in two and three dimensions. In particular, \edit{we modify} Laplace's {single-layer potential} (representing the solution of the exterior Neumann Laplace problem) \edit{and} Helmholtz layer potentials (in the context of a sound-soft scattering problem)\edit{. With some given quadrature rule, the resulted modified representations allow us to }obtain better approximations compared to standard representations. The proposed modifications are based on subtracting {specific solutions} (or \textit{auxiliary functions}) of the {PDE at stake}. The use of auxiliary functions have been developed in the context of Boundary Regularized Integral Equation Formulation (BRIEF)~\cite{KSC12,SKKC14,SKKC15} to regularize the representation formula on the boundary, or in the context of density interpolation techniques~\cite{Perez18,PeTuFa19,FaPeBo21} to regularize layer potentials {(generalization of density subtractions)}. Those techniques commonly consider multiple auxiliary functions, and may require to solve additional problems to find such functions. The proposed work concentrates on regularizing nearly singular integrals using explicitly one analytic auxiliary function, and when representing the solution with layer potentials involving only one density (no representation formula). We provide several examples of {auxiliary} functions (and compare them), and provide guidelines to find them. The proposed modified representations {are simple and easy to implement, it }allows {one} to straightforwardly gain accuracy in evaluating the solution{, especially when computational resources are limited. }This work provides valuable insights into Laplace and Helmholtz layer potentials. Additionally this can also be applied to modify boundary integral equations to avoid weakly singular integrals.

The paper is organized as follows: {Section \ref{S:motiv} presents some context and motivation for the proposed modified representations, Section \ref{sec:modified_rep} establishes the modified representations and general guidelines to find appropriate auxiliary functions}. Sections \ref{S:application} and \ref{S:BIE} illustrate the efficiency of the modified representations for Laplace and Helmholtz in two and three dimensions, off and on boundary. Finally, Section \ref{S:conclu} presents our concluding remarks, Appendices \ref{sec:appdx2D} and \ref{sec:appdx3D} provide a brief summary of the Nyström methods used in two and three dimensions, and Appendix~\ref{sec:appdx_proofs} details some proofs for Section \ref{sec:modified_rep}. 
\section{Motivation for modified representations}\label{S:motiv}
Consider a \textit{domain} $D \subset \R^d$, $d = 2,3$, that is a bounded simply connected open set with smooth boundary (of class $\mathcal{C}^2$), and a linear elliptic partial differential equation of the form $\mathcal{L}u = 0$. It is common to represent the solution $\vv$ of that PDE using the so-called representation formula (e.g. \cite[Theorem 6.5]{Kress89},~\cite[Theorem 3.1]{CoKr13}). In particular for $\vv$ satisfying $\mathcal{L}\vv = 0$ in $D$, we have the following identities:
\begin{equation}\label{eq:repres_form}
 \begin{aligned}
& \dsp \int_{\partial D}  \dsp  \partial_{n_y}G(x,y)\vv (y) d \sigma_y - \dsp \int_{\partial D}  G(x,y)\dvv (y) d \sigma_y  =  \begin{cases}
 -\vv (x)   &\, x \in D ,\\
- \dsp \frac{1}{2} \vv (x)&\, x \in \partial D,\\
 \,\,\,\,\,0 &\,x \in E: = \mathbb{R}^d \setminus \bar{D},
 \end{cases}
 \end{aligned}
\end{equation} 
where $G$ denotes the fundamental solution of considered PDE, $n_y$ is the unit outward normal of $D$ at $y$, and $d \sigma_y$ is the integration surface element. For instance, \eqref{eq:repres_form} holds true for $\mathcal{L} := \Delta$ and $\mathcal{L} := \Delta + k^2$, the Laplace and the Helmholtz equation, respectively. The goal of this paper is to use \eqref{eq:repres_form} with well-chosen $\vv$ to modify the representation of the solution of boundary value problems associated to $\mathcal{L}$. Let us illustrate the strategy with for example the Exterior Neumann Laplace problem:
\begin{equation}\label{eq:exterior_bvp}
\left|
\begin{aligned}
&\mbox{Find } u \in \mathcal{C}^2(E) \cap \mathcal{C}^1(\bar{E}:=\mathbb{R}^d \setminus D) \mbox{ such that:}\\
&\Delta u = 0 \quad \text{in }E,\quad \partial_{n} u= g \quad \text{on }\partial D,\quad  \lim \limits_{|x| \to \infty} u(x) =  o(1),
\end{aligned}
\right.
\end{equation}
with some smooth data $g$ (with null average). The solution of Problem \eqref{eq:exterior_bvp} can be represented using the Green's formula~\cite{guenther1996partial,CoKr13}:
\begin{equation}\label{eq:repres_laplaceN}
\begin{aligned}
u(x) & =  \displaystyle \int_{\partial D} \partial_{n_y} G(x,y) u (y) \, d\sigma_y -  \displaystyle \int_{\partial D} G(x,y) \partial_{n_y} u (y) \, d\sigma_y , \quad x \in E, \\
& =    \displaystyle \int_{\partial D} \partial_{n_y} G(x,y) u (y) \, d\sigma_y - \displaystyle \int_{\partial D} G(x,y) g(y) \, d\sigma_y, \quad x \in E,
\end{aligned}
\end{equation}
where 
\begin{equation}\label{eq:kernel_Laplace}
G(x,y) = \begin{cases} &  \dsp -\frac{1}{2\pi} \log |x-y| \quad \mbox{for }d = 2,\\
& \dsp \frac{1}{4 \pi} \frac{1}{|x-y|} \quad \mbox{for }d = 3,
\end{cases}
\end{equation}
and the trace on the boundary satisfies the boundary integral equation of the second kind: 
\begin{equation}\label{eq:BIE_laplace_Ne}
\frac{1}{2} u(x^\ast) - \displaystyle \int_{\partial D} \partial_{n_y} G(x^\ast,y) u(y)\, d\sigma_y = \displaystyle \int_{\partial D} G(x^\ast,y) g (y) \, d\sigma_y, \quad x^\ast \in \partial D.
\end{equation}
The fundamental solution $G$ is singular when $y = x^\ast$. \edit{For $x \in \mathbb{R}^d \setminus \partial D$, assume we can write $ x = x^\ast \pm \ell n_{x^\ast}$ with $n_{x^\ast}$ the unit outward normal at $x^\ast$, and $\ell>0$ the distance from the boundary. Then $G$ is \textit{nearly singular} at $y = x^\ast$ when $|x-y| = \ell \ll 1$ (i.e. when $x$ is close to the boundary).} A layer potential is said to be a weakly singular integral (resp. a nearly singular integral) when its {kernel} ($G$ or $\partial_{n}G$ in the cases above) is singular at $y = x^\ast$ (resp. nearly singular \edit{at $y = x^\ast$}). There exist high-order quadrature rules to approximate weakly singular integrals with very high accuracy (e.g. \cite{atkinson1997numerical,bruno2001fast,ganesh2004high,bremer2010nonlinear}). However, high accuracy is lost for nearly singular integrals: this is the so-called close evaluation problem. Assuming we have solved \eqref{eq:BIE_laplace_Ne}, we can modify \eqref{eq:repres_laplaceN} using \eqref{eq:repres_form} to address the close evaluation problem. Taking the difference we obtain
\begin{equation}\label{eq:modified_repres_laplaceN}
\begin{aligned}
u(x) & =  \displaystyle \int_{\partial D} \partial_{n_y} G(x,y) [u (y) -\vv(y)] \, d\sigma_y -  \displaystyle \int_{\partial D} G(x,y) [g(y)- \dvv (y)] \, d\sigma_y, \quad x \in E.
\end{aligned}
\end{equation}
If one finds $\vv$ such that $\vv(x^\ast)  = u(x^\ast)$ and $\dvvx(x^\ast) = g(x^\ast)$, where $x^\ast \in \partial D$ denotes the closest boundary point of the evaluation point $x$ ($ x = x^\ast + \ell n_{x^\ast}$), then \eqref{eq:modified_repres_laplaceN} doesn't suffer from the close evaluation problem. \\
Similarly, one can represent the solution of Problem \eqref{eq:exterior_bvp} using \textit{a single-density representation} \edit{given by} the single-layer potential:
\begin{equation}\label{eq:SLP_laplaceN}
u(x) = \displaystyle \int_{\partial D} G(x,y) \rho(y)\, d\sigma_y, \quad x \in D,
\end{equation}
with $\rho$ a continuous density solution of the boundary integral equation of the second-kind:
\begin{equation}\label{eq:BIE_laplace_NS}
-\frac{1}{2} \rho(x^\ast) + \displaystyle \int_{\partial D} \partial_{n_x^\ast} G(x^\ast,y) \rho(y)\, d\sigma_y = g(x^\ast), \quad x^\ast \in \partial D.
\end{equation}
Assuming we have solved \eqref{eq:BIE_laplace_NS} for $\rho$, subtracting \eqref{eq:repres_form} from \eqref{eq:SLP_laplaceN} we obtain
\begin{equation}\label{eq:modified_SLP_laplaceN}
\begin{aligned}
u(x) & = \displaystyle \int_{\partial D} G(x,y) [\rho(y)- \dvv (y)] \, d\sigma_y +  \displaystyle \int_{\partial D} \partial_{n_y} G(x,y) \vv(y) \, d\sigma_y, \quad x \in E.
\end{aligned}
\end{equation}
If one finds $\vv$ such that $\vv(x^\ast)  = 0$ and $\dvvx(x^\ast) = \rho(x^\ast)$, then \eqref{eq:modified_SLP_laplaceN} doesn't suffer from the close evaluation problem. \\
Representations \eqref{eq:modified_repres_laplaceN} and \eqref{eq:modified_SLP_laplaceN} are attractive representations, and several works have provided guidelines on how to build appropriate solutions $\vv$. For \eqref{eq:modified_repres_laplaceN} one can use Taylor-like functions $\vv(x) = u(x^\ast)\tilde{g}(x) + \partial_{n_{x^\ast}}u(x^\ast)\tilde{f}(x)$, with $\tilde{g}$ and $\tilde{f}$ solutions of some Laplace boundary value problems~\cite{KSC12,SKKC14,SKKC15}. This technique has been first developed in the context of Boundary Regularized Integral Equation Formulation (BRIEF) (namely to solve \eqref{eq:BIE_laplace_Ne} using the same subtraction technique on boundary) and applied to evaluate the solution near the boundary. For \eqref{eq:modified_SLP_laplaceN} one can use density interpolation methods~\cite{Perez18,PeTuFa19,FaPeBo21}: $\vv = \vv(x^\ast,y) = \sum_{j=0}^J c_j(y) H_j(x^\ast-y)$ where $(H_j)_j$ satisfy the PDE (in the above case $(H_j)_j$ are harmonic functions). In both methods \edit{the} chosen auxiliary functions $\vv$ necessarily depend on the trace $u$ (and/or normal trace $\partial_{n}u$), or the density $\rho$ at the closest evaluation point\edit{. Furthermore they} require to satisfy at least two conditions (two boundary value problems or two boundary conditions). \\
In this paper we provide another construction of modified representations for single-density representations of Laplace and Helmholtz boundary value problems. The construction relies on auxiliary functions $\vv$ that are independent of the density (solution of the boundary integral equation), and requires fewer constraints in the context of \eqref{eq:SLP_laplaceN}. \edit{As a consequence, our approach provides more freedom in choosing $\vv$.} The proposed modified representations are also simple to implement and do not add significant computational costs. In what follows we provide modified representations for Laplace and Helmholtz in 2D and 3D, and provide several examples to illustrate the efficiency of the method.
\section{Modified representations}\label{sec:modified_rep}
We present modified representations for single-density representations of Laplace and Helmholtz boundary value problems. In particular, we consider the interior Dirichlet Laplace problem (where one can represent the solution using the double-layer potential), the exterior Neuman Laplace problem \eqref{eq:exterior_bvp} (using the single-layer potential \eqref{eq:SLP_laplaceN}), and the sound-soft scattering problem.
 \subsection{\edit{Modified representation for the Laplace double-layer potential}}\label{ssec:DLP_ds}
The interior Dirichlet problem for Laplace consists in finding $u \in \mathcal{C}^2(D) \cap \mathcal{C}^1(\overline{D})$ such that
\begin{equation}\label{eq:interior_pb}
\left|
\begin{aligned}
&\Delta u = 0 \quad \text{in }D,\quad u=f \quad \text{on }\partial D,
\end{aligned}
\right.
\end{equation}
with some smooth data $f$. The solution of Problem \eqref{eq:interior_pb} can be represented as a double-layer potential~\cite{guenther1996partial,CoKr13}:
\begin{equation}\label{eq:DLP_laplace}
u(x) = \displaystyle \int_{\partial D} \partial_{n_y} G(x,y) \mu(y)\, d\sigma_y, \quad x \in D,
\end{equation} with $G$ defined in \eqref{eq:kernel_Laplace}, and 
$\mu$ a continuous density solution of the boundary integral equation: 
\begin{equation}\label{eq:BIE_laplace}
- \frac{1}{2} \mu(x^\ast) + \displaystyle \int_{\partial D} \partial_{n_y} G(x^\ast,y) \mu(y)\, d\sigma_y = f(x^\ast), \quad x^\ast \in \partial D.
\end{equation}
We now make use of \eqref{eq:repres_form} to modify \eqref{eq:DLP_laplace}. One can show the following (see Appendix \ref{ssec:appdx_dlp} for details):
\begin{Proposition}\label{pro:DLP}
Given $x = x^\ast - \ell n_{x^\ast} \in D$ with $x^\ast \in \partial D$, let $\vv$ be a solution of Laplace's equation in $D \subset \mathbb{R}^d$, $d = 2,3$, such that 
\begin{equation}\label{eq:hyp_dlp}
\vv (x^\ast) = 1, \quad \dvvx (x^\ast) = 0.
\end{equation}
The solution of the exterior Dirichlet Laplace problem \eqref{eq:DLP_laplace} admits the modified representation:
\begin{equation}\label{eq:general_DLP}
\begin{aligned}
& u(x)   =  \dsp \int_{\partial D} \partial_{n_y} G(x,y) \mu(y)\left[ 1 - \vv(y)\right]  \, d\sigma_y +  \dsp \int_{\partial D} \partial_{n_y} G(x,y) \left[ \mu(y)  -\right .  \left. \mu(x^\ast)\right]\vv(y)  \, d\sigma_y  \\
& - \mu(x^\ast) \vv(x^\ast) + \mu(x^\ast) \dsp \int_{\partial D} G(x,y) \left[  \dvv(y)- \dvvx(x^\ast)\right]  \, d\sigma_y - \mu(x^\ast) \dvvx(x^\ast),  \quad x \in D.
\end{aligned}
\end{equation}
The modified representation \eqref{eq:general_DLP} has smoother integrands than \eqref{eq:DLP_laplace}, and it addresses the close evaluation problem, in the sense that nearly singular terms vanish as $y \to x^\ast$.
\end{Proposition}
From Proposition \ref{pro:DLP} we can now build auxiliary functions $\vv$ independent of $\mu$, and there exist plenty of candidates: constant, linear, based on the Green's function ($\vv(y) = G(y,x_0)$ with $x_0 \in E$), quadratic ($\vv(y_1, y_2) = 1 + (y_1 - x_1^\ast)(y_2 - x_2^\ast)$, $\vv(y_1, y_2)  = 1 + (y_1 - x_1^\ast)^2 - (y_2 - x_2^\ast)^2$), $\vv(y_1, y_2,y_3) = e^{y_3}(\sin y_1 + \sin y_2)$, etc. 
The solution $\vv \equiv 1$ naturally satisfies the conditions \eqref{eq:hyp_dlp}, and the modified representation \eqref{eq:general_DLP} boils down to 
\begin{equation}\label{eq:modified_dlp_laplace}
u(x)   =  \displaystyle \int_{\partial D} \partial_{n_y} G(x,y) [\mu(y) - \mu (x^\ast)] \, d\sigma_y - \mu(x^\ast) , \quad x \in D.
\end{equation}
The modified representation \eqref{eq:modified_dlp_laplace} {is well-known and widely used \edit{(e.g.~\cite{Hwang13,barnett2014evaluation,ckk2020asymp})}, it} is the simplest representation that naturally addresses the close evaluation problem. {Thus, we do not provide numerical results for this case. Rather, we concentrate on other layer potentials.}
 \subsection{\edit{Modified representation for the Laplace single-layer potential}}\label{ssec:SLP_ds}
{Going back to Problem \eqref{eq:exterior_bvp}, o}ne can show the following (see Appendix \ref{ssec:appdx_slp} for details):
\begin{Proposition}\label{pro:SLP}
Given $x = x^\ast + \ell n_{x^\ast} \in E$ with $x^\ast \in \partial D$, let $\vv$ be a solution of Laplace's equation in $D \subset \mathbb{R}^d$, $d = 2,3$, such that
\begin{equation}\label{eq:hyp_slp}
\dvvx(x^\ast) = 1.
\end{equation}
The solution of the exterior Neumann Laplace problem \eqref{eq:exterior_bvp} admits the modified representation:
\begin{equation}\label{eq:general_SLP}
\begin{aligned}
 u(x)  = \dsp \int_{\partial D}  G(x,y) \rho(y) & \left[ 1 - \dvv (y)\right]  \, d\sigma_y  +  \dsp \int_{\partial D} G(x,y) \left[ \rho(y)  - \rho(x^\ast)\right] \dvv (y)  \, d\sigma_y  \\
&  + \rho(x^\ast) \dsp \int_{\partial D}  \partial_{n_y} G(x,y) \rho(y)\left[\vv(y)-\vv(x^\ast)\right]  \, d\sigma_y, \quad \forall x  \in E.
\end{aligned}
\end{equation}
The modified representation \eqref{eq:general_SLP} {has} smoother integrands {than \eqref{eq:SLP_laplaceN}}.
\end{Proposition}
{Contrary to auxiliary functions provided in Taylor-like methods and density interpolation methods (discussed in Section \eqref{S:motiv}), auxiliary functions $\vv$ do not depend on $\rho$ and rely on only one constrain \eqref{eq:hyp_slp}. Therefore, there is a lot of freedom in choosing $\vv$: given $\mathtt{u}$ a solution of Laplace's equation, then one chooses $\vv := \frac{\mathtt{u}}{\partial_{n_x^\ast} \mathtt{u}(x^\ast)}$ (as long as $\partial_{n_x^\ast} \mathtt{u}(x^\ast) \neq 0$). Candidates may then include:}
 \begin{itemize}
 \item the linear function $\vv(y) = n_{x^\ast} \cdot y$ ;
   \item the function $\vv(y) = 2^{d-1} \pi G(y, x^\ast + n_{x^\ast})$ based on the Green's function ; 
  \item { the quadratic product function $\vv(y) =  \dsp \frac{(y_1 - x_{0,1})(y_2 - x_{0,2})}{n_{x^\ast,1}(x_2^\ast - x_{0,2}) +n _{x^\ast,2}(x_1^\ast - x_{0,1})} $, $x_0 \in D $ };
\item { the quadratic difference function $\vv(y) =  \dsp\frac{1}{2} \frac{(y_1 - x_{0,1})^2 - (y_2 - x_{0,2})^2}{n_{x^\ast,1}(x_1^\ast - x_{0,1}) - n _{x^\ast,2}(x_2^\ast - x_{0,2})} $, $x_0 \in D $ }.
\end{itemize}  
{Note that the above candidates are valid in $\mathbb{R}^d$, one can also consider any of the quadratic functions above in $\mathbb{R}^3$ as a function of $(y_i,y_j)$, $i,j = 1,2,3$, $j \neq i$.} In Section \ref{S:application} we will test \eqref{eq:general_SLP} using {several candidates $\vv$} and make comparisons. 
The modified representation \eqref{eq:general_SLP} adds two terms to compute compared to \eqref{eq:SLP_laplaceN}, it is the price to pay to gain accuracy at close evaluation points. We will make comparative tests to quantify this aspect.
 \subsection{\edit{Modified representation for the Helmholtz double- and single-layer potentials}}\label{ssec:helm_ds}
\noindent  {We consider in this case the sound-soft scattering problem:
\begin{equation}\label{eq:scattering_pb}
\left|
\begin{aligned}
& \text{Find }u \in \mathcal{C}^2 ( E) \cup \mathcal{C}^1(\bar{E}) \text{ such that:}\\
&\Delta u + k^2 u  = 0 \quad\text{in } E,\quad u=f \quad\text{on } \partial D,\quad  \dsp \lim \limits_{R \rightarrow \infty}\int_{|y|=R} | \partial_{n} u- i k u |^2 \, d\sigma_y =0, 
\end{aligned}
\right.
\end{equation}
with some smooth data $f$ associated to the wavenumber $k$. Above, the last condition represents the Sommerfeld radiation condition. The solution of Problem \eqref{eq:scattering_pb} can be represented as a combination of double- and single-layer potentials~\cite{Kress91}:
\begin{equation}\label{eq:DLP-SLP_helm}
u(x) = \dsp \int_{\partial D} \left[ \partial_{n_y} G^H(x,y) - ik G^H(x,y) \right] \mu (y) \, d \sigma_y, \quad x \in E,
\end{equation}
with $G^H$ defined by 
\begin{equation}\label{eq:kernel_helm}
G^H(x,y) = \begin{cases} & \dsp \frac{i}{4} H^{(1)}_0 (k |x -y|), \quad \mbox{for } d =2, \\
& \dsp \frac{1}{4 \pi } \frac{e^{i k |x-y|}}{|x-y|}, \quad \mbox{for } d =3,
\end{cases}
\end{equation}
with $H^{(1)}_0 (\cdot)$ the Hankel function of the first kind, and $\mu$ a continuous density satisfying: 
\begin{equation}\label{eq:BIE_helm}
\begin{aligned}
\frac{1}{2} \mu(x^\ast)  + \dsp \int_{\partial D} \left[ \partial_{n_y} \right. & G^H(x^\ast,y)  \left. - ik G^H(x^\ast,y) \right] \mu (y) \, d \sigma_y = f(x^\ast), \quad x^\ast \in \partial D.
\end{aligned}
\end{equation}
One obtain the following:
 }
\begin{Proposition}\label{pro:helm}
Given $x = x^\ast + \ell n_{x^\ast} \in E$ with $x^\ast \in \partial D$, let $\vv$ be a solution of Helmholtz equation in $D \subset \mathbb{R}^d$, $d = 2,3$, such that
\begin{equation}\label{eq:hyp_helm}
\vv (x^\ast) = 1, \quad \dvvx (x^\ast) = i k.
\end{equation}
Then the solution of the sound-soft scattering problem \eqref{eq:scattering_pb} admits the modified representation:
\begin{equation}\label{eq:general_helm}
\begin{aligned}
u(x)   = \dsp & \int_{\partial D}  \left[ \partial_{n_y}G^H(x,y) - \dvv (y) G^H(x,y) \right]  \left[ \mu(y) - \mu (x^\ast) \right]  \, d\sigma_y \\
&+  \dsp \int_{\partial D} G^H(x,y)\left[\dvv (y) - ik \right]  \mu(y)    \, d\sigma_y  + \mu(x^\ast) \dsp \int_{\partial D}  \partial_{n_y} G^H(x,y) \left[1 - \vv(y)\right]  \, d\sigma_y,  \quad \forall x  \in E.
\end{aligned}
\end{equation}
The modified representation \eqref{eq:general_helm} {has} smoother integrands {than \eqref{eq:DLP-SLP_helm}}.
\end{Proposition}
The proof can be found in Appendix \ref{ssec:appdx_helm}. One can check in particular that plane waves $\vv(y) = e^{i k n_{x^\ast} \cdot (y-x^\ast)}$ do satisfy \eqref{eq:hyp_helm}, whereas Green-based functions like $\vv(y) = G^H(y, x_\ast + n_{x^\ast})$ (up to some constant) cannot. We will use \eqref{eq:general_helm} with plane waves for the numerical examples. 

\section{Numerical examples}\label{S:application}
The accuracy in approximating \eqref{eq:DLP_laplace}--\eqref{eq:modified_dlp_laplace}, \eqref{eq:SLP_laplaceN}--\eqref{eq:general_SLP}, \eqref{eq:DLP-SLP_helm}--\eqref{eq:general_helm} respectively, relies on the resolution of the boundary integral equation \eqref{eq:BIE_laplace}, \eqref{eq:BIE_laplace_NS}, \eqref{eq:BIE_helm} respectively. In what follows we assume that the boundary integral equations are sufficiently resolved\edit{. G}iven the density's resolution, we compare the representations and their modified ones through several examples. All the codes can be found in~\cite{Carvalho-codes}.
\subsection{Exterior Neumann Laplace problem}
\subsubsection{Example 1: exterior Laplace in two dimensions}\label{ssec:Laplace2D}
Since $\partial D$ is a closed smooth boundary, we use the Periodic Trapezoid Rule (PTR) to approximate \eqref{eq:SLP_laplaceN} {and \eqref{eq:general_SLP}, \edit{ where we will use} several $\vv$ according to Proposition \ref{pro:SLP}.}
We consider an exact solution of Problem \eqref{eq:exterior_bvp}:
\[\begin{aligned}
u_{\text{exact}}(x) =  u_{\text{exact}}(x_1, x_2) &=  \frac{x_1 - x_{0,1}}{ |x - x_0|^2}, \quad x_0 = (x_{0,1}, {x_{0,2}}) \in {D},
\end{aligned}\] 
which consists in choosing $g(x^\ast) = \partial_{n_{x^\ast}}u_{\text{exact}}(x^\ast)$, for any $x^\ast \in \partial D$. 
{
All simulations are done outside of a kite-shaped domain using the Periodic Trapezoid Rule with $N = 128$ quadrature points for the following representations:
\begin{itemize}\renewcommand{\labelitemi}{$\bullet$}
\item \textbf{V0:} standard representation \eqref{eq:SLP_laplaceN};
\item \textbf{V1:} modified representation \eqref{eq:general_SLP} with the linear function $\vv_1 (y)= n_{x^\ast} \cdot y $;
\item \textbf{V2:} modified representation \eqref{eq:general_SLP} with the Green-based function\\ $\vv_2(y) =  2 \pi G(y, x^\ast + n^\ast)$; 
\item \textbf{V3:} modified representation \eqref{eq:general_SLP} with the quadratic function\\$\vv_3(y) =  \dsp\frac{1}{2} \frac{y_1 ^2 - y_2^2}{n_{x^\ast,1}x_1^\ast  - n _{x^\ast,2}x_2^\ast }$ ; 
\item \textbf{V4:} modified representation \eqref{eq:general_SLP} with the quadratic function\\ $\vv_4(y) =  \dsp \frac{(y_1 - 5)(y_2 - 5)}{n_{x^\ast,1}(x_2^\ast - 5) +n _{x^\ast,2}(x_1^\ast - 5)} $. 
\end{itemize} 
\begin{figure}[H]
  \centering
 \includegraphics[width=0.3\textwidth]{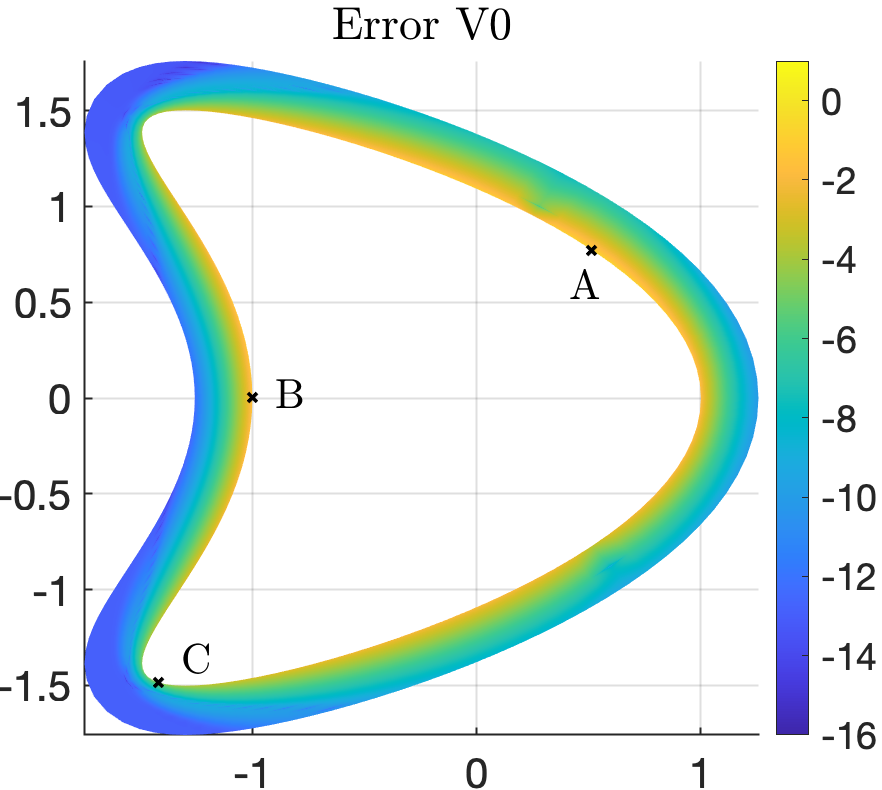}
\includegraphics[width=0.3\textwidth]{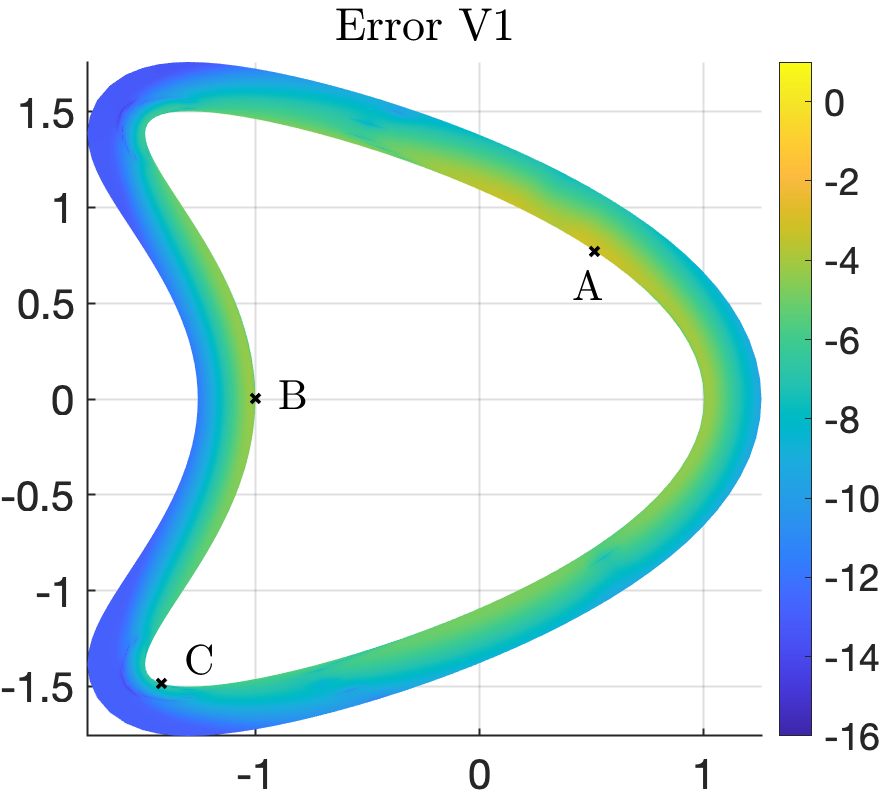}
 \includegraphics[width=0.3\textwidth]{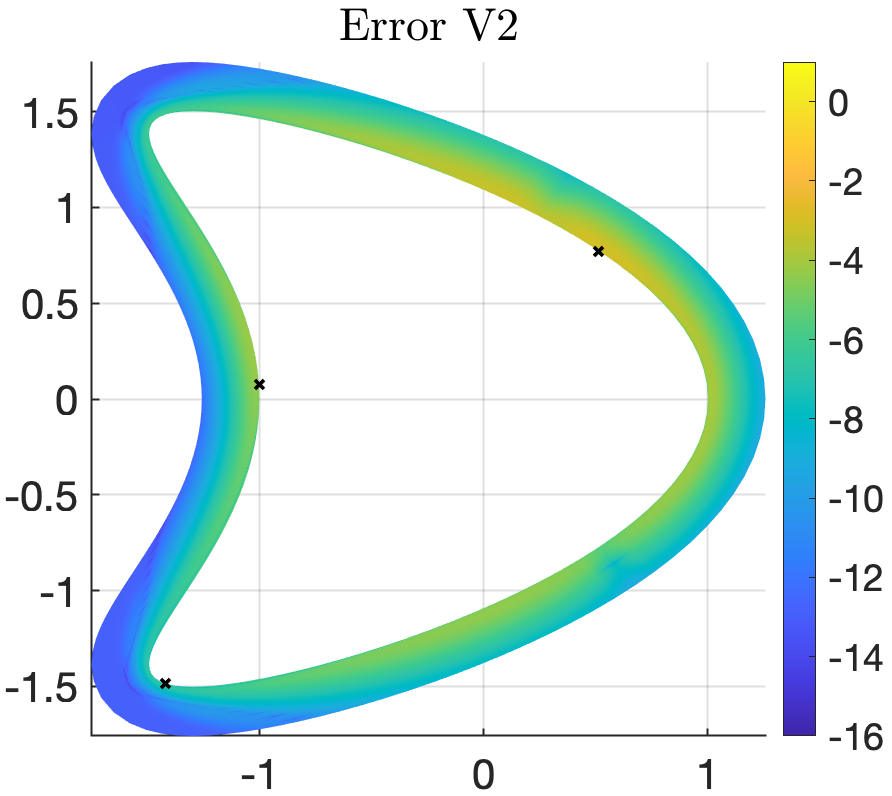}\\
\includegraphics[width=0.3\textwidth]{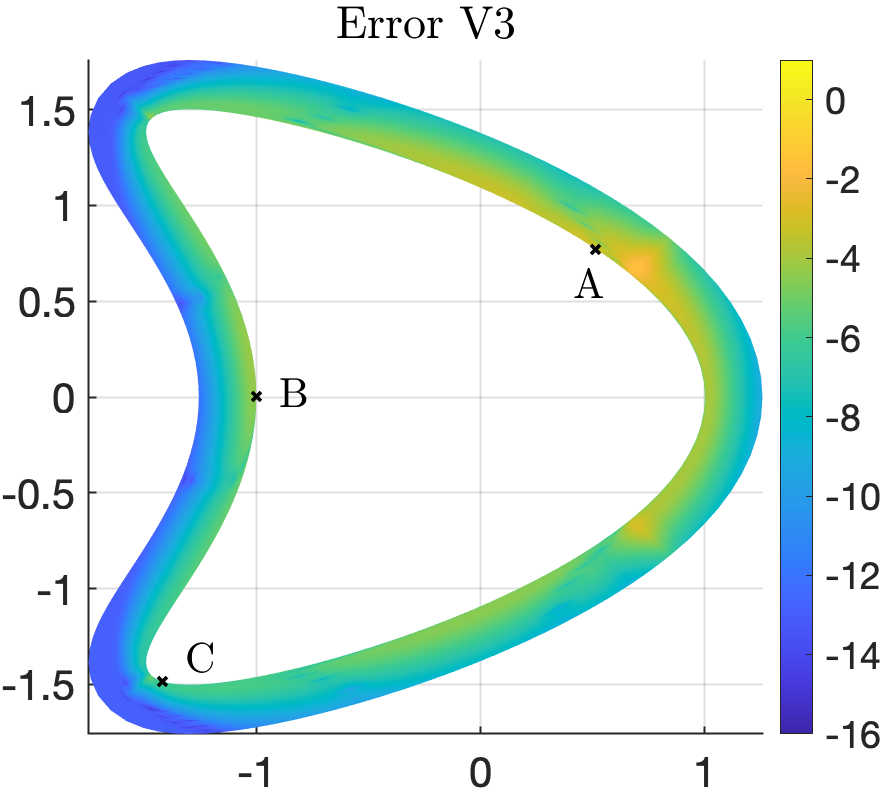}
\includegraphics[width=0.3\textwidth]{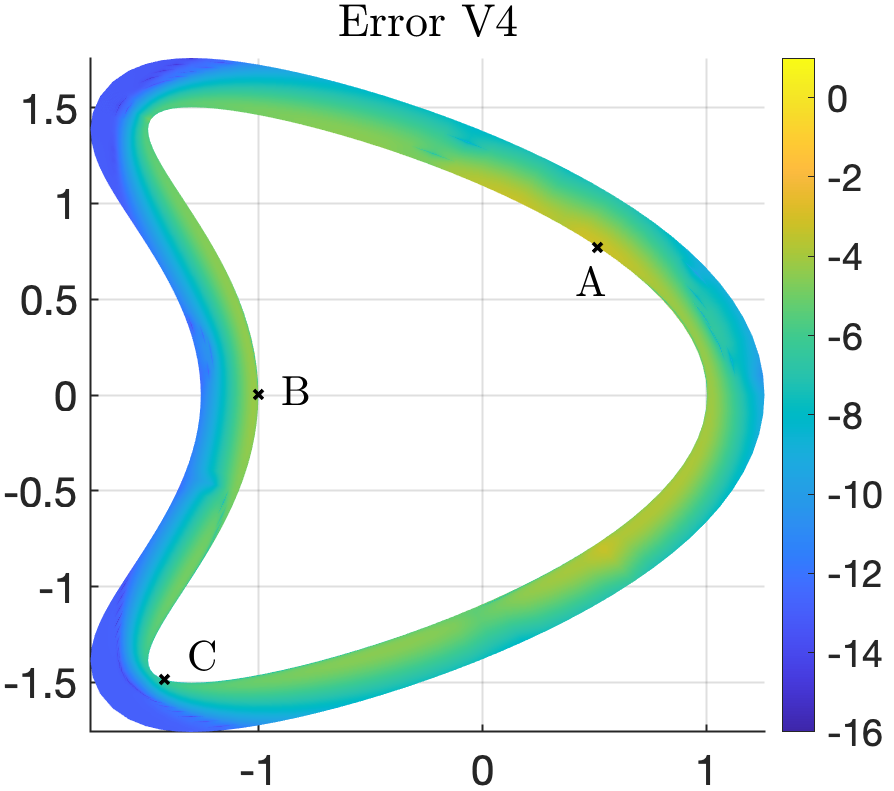}
  \caption{\small {\textbf{Laplace 2D single-layer.} Plots of $\log_{10}$ of the error for the evaluation of the solution of \eqref{eq:exterior_bvp} out of the kite domain defined by the boundary $y(t) =(\cos t + 0.65 \cos(2t) - 0.65, 1.5 \sin t)$, $t \in [0 ,2\pi]$, for the Neumann data, $g = \partial_{n} u_{\text{exact}}$ with $x_0=(0.1,0.4)$, for representations V0, V1, V2, V3, V4 computed using PTR with $N = 128$. Computations are made on a boddy-fitted grid with $N \times 200$ grid points.} }
  \label{img:SLP2D_1_new}
\end{figure}
\begin{figure}[H]
  \centering
 \includegraphics[width=0.32\textwidth]{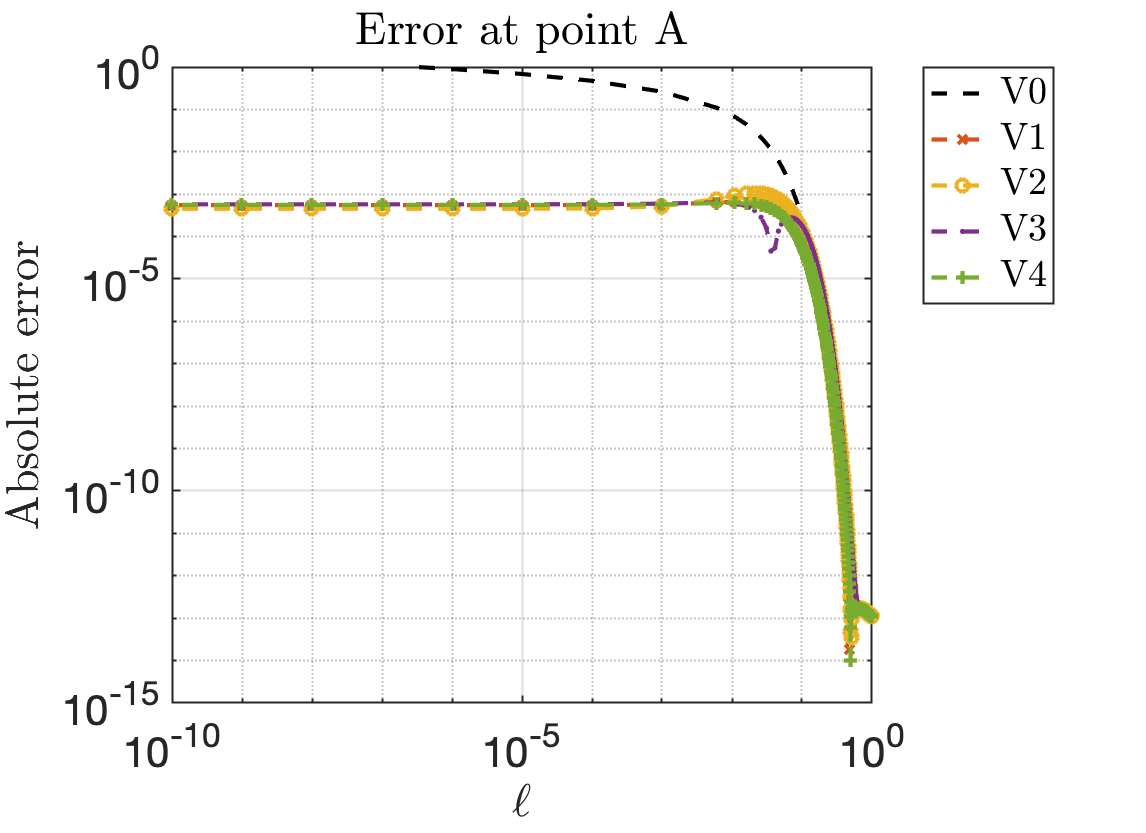}
\includegraphics[width=0.32\textwidth]{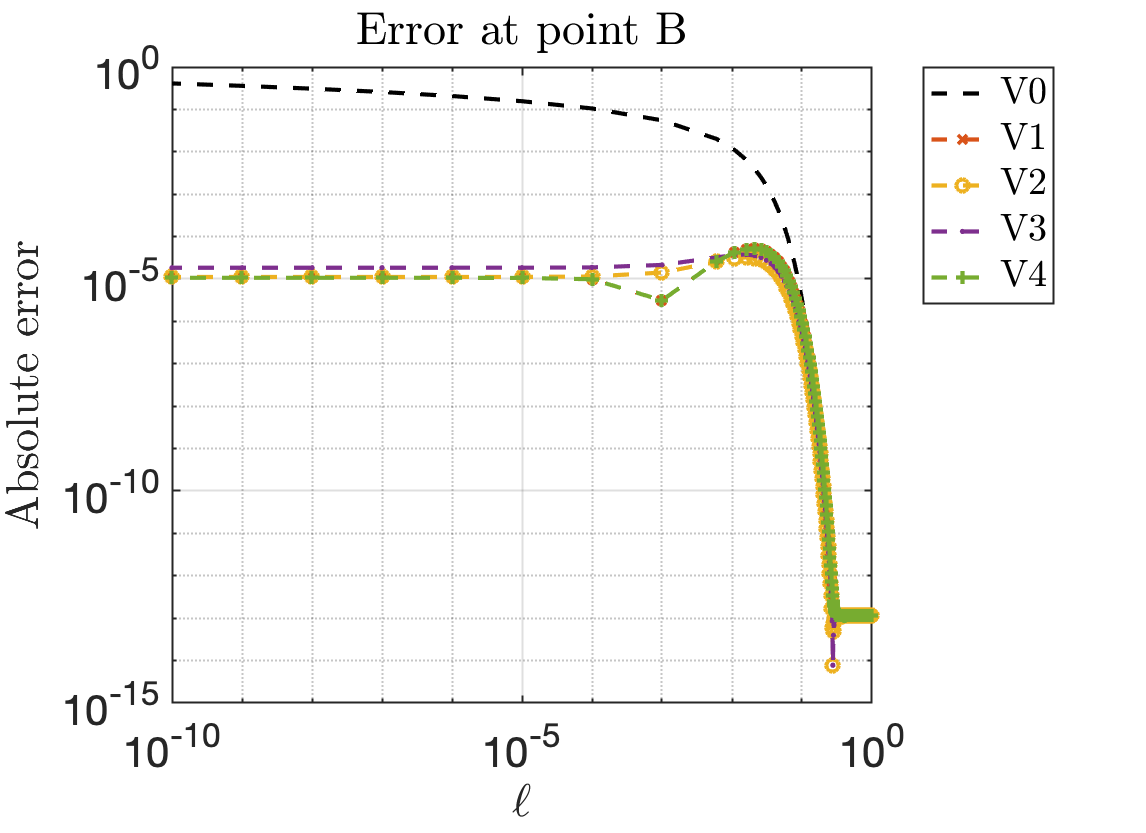}
 \includegraphics[width=0.32\textwidth]{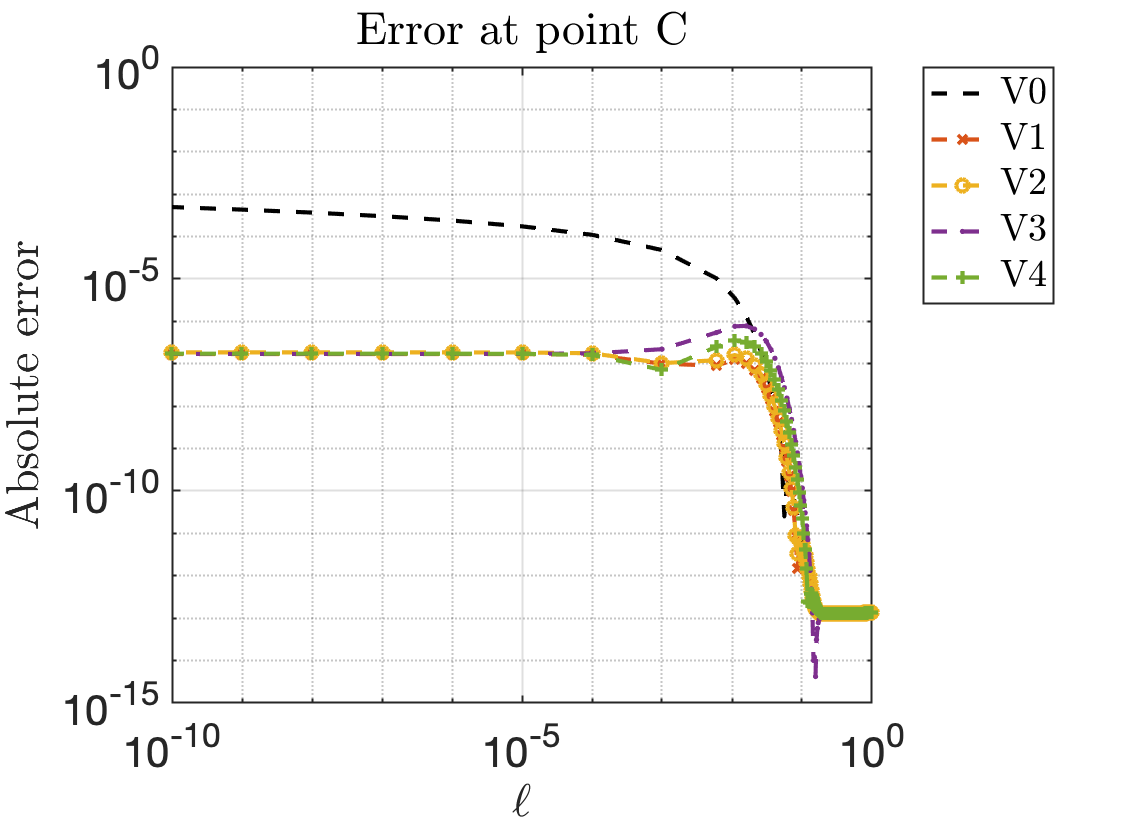}
  \caption{\small {\textbf{Laplace 2D single-layer.} Log-log plots of the errors with respect to $\ell$ made in computing the solution (as described in Figure \ref{img:SLP2D_1_new}) along the normal of the three points A, B, C, plotted as black $\times$'s in Figure \ref{img:SLP2D_1_new}.} }
  \label{img:SLP2D_2_new}
\end{figure}
\begin{figure}[H]
  \centering
 \includegraphics[width=0.32\textwidth]{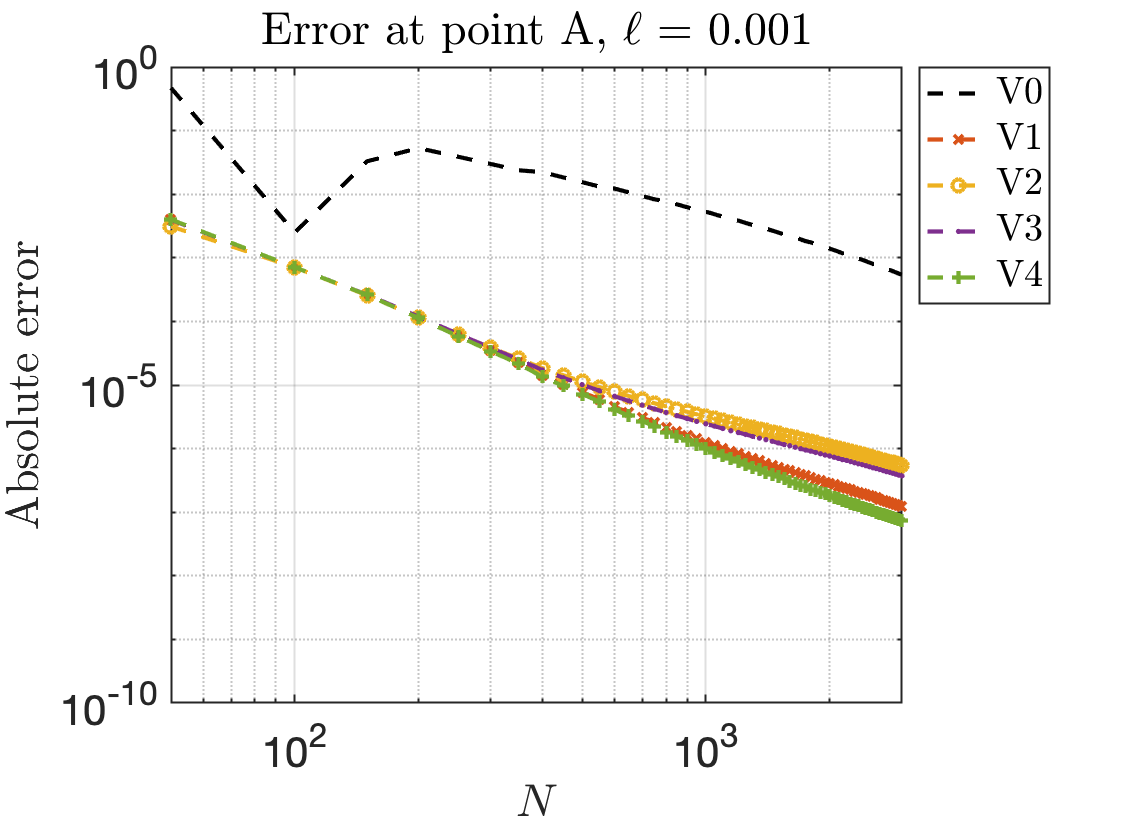}
\includegraphics[width=0.32\textwidth]{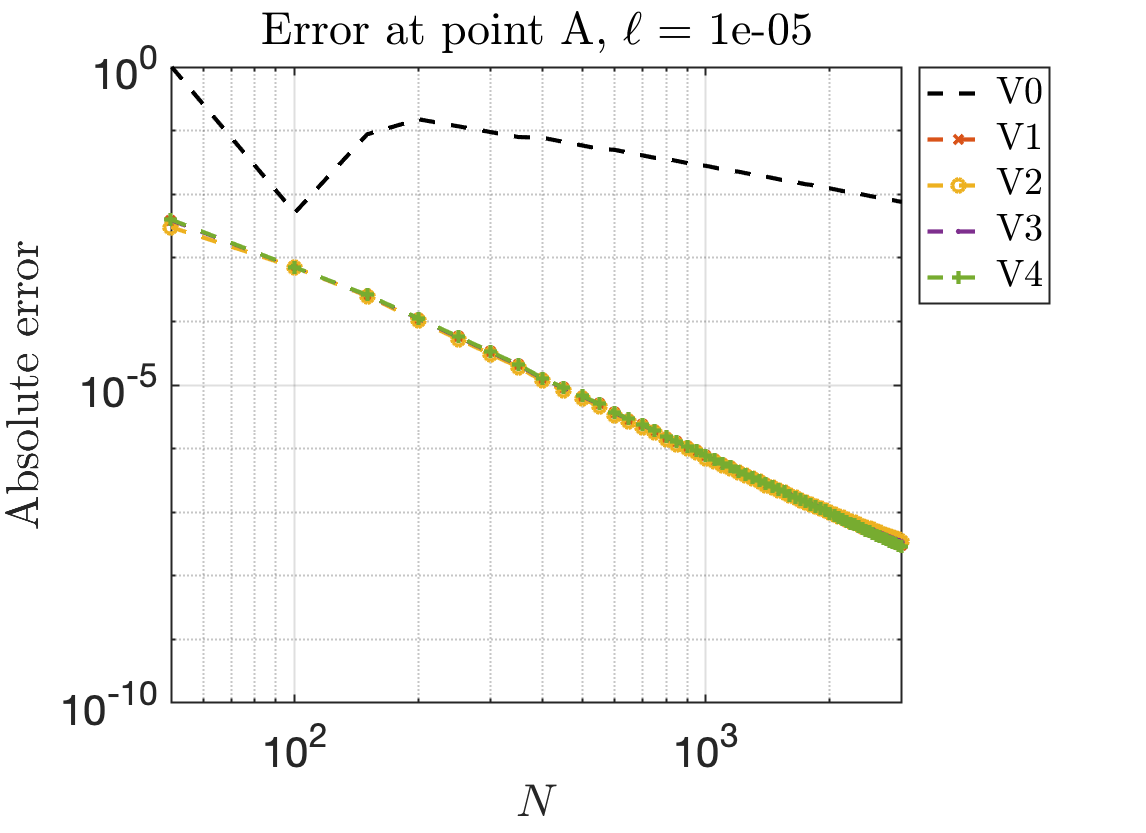}
 \includegraphics[width=0.32\textwidth]{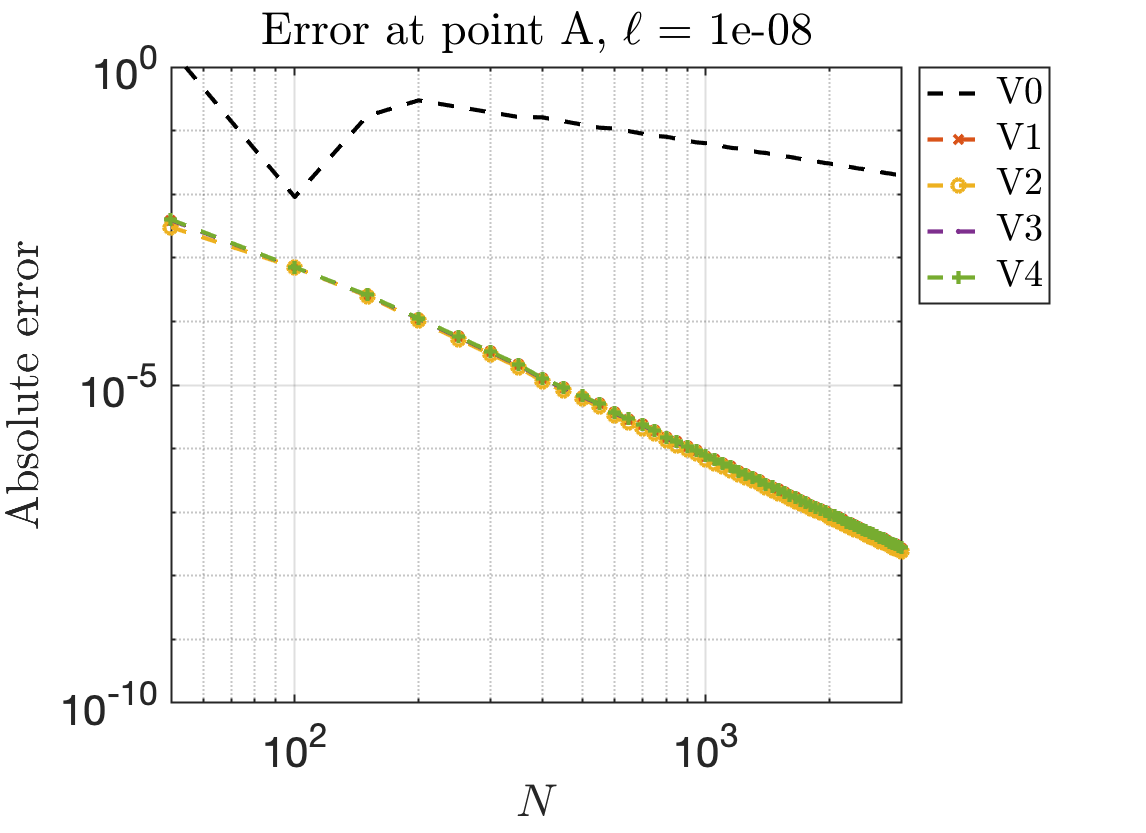}
  \caption{\small {\textbf{Laplace 2D single-layer.} Log-log plots of the errors with respect to $N$ made in computing the solution at some distance $\ell$ along the normal from the point A plotted as black $\times$'s in Figure \ref{img:SLP2D_1_new}.} }
  \label{img:SLP2D_3_new}
\end{figure}
}
We solved \eqref{eq:BIE_laplace_NS} using the {Nyström method based on the} Periodic Trapezoid Rule {(using Matlab classic \textit{backslash})}. The accuracy of all methods is limited by the accuracy of the resolution for $\rho$ {(in particular when considering moderate $N$). This can be assessed by looking the density's Fourier coefficients decay: in this case \edit{the} coefficients decay is bounded by $10^{-5}$ for $N =128$}. Results in Figures \ref{img:SLP2D_1_new} and \ref{img:SLP2D_2_new} show that given $\rho$ resolved, the approximation of the modified representations provide better results overall. Far from the boundary, all methods approximate well the solution. As the evaluation point gets closer to the boundary ($\ell \to 0$), {V0 approximated by PTR} suffers from the close evaluation problem and the error increases {(see \cite{barnett2014evaluation})}. Note that the single-layer potential commonly suffers less from this phenomenon than the double-layer potential~(e.g. \cite{carvalho2018asymptotic}). Using the modified representations {(V1--V4)} allows to reduce the error by a couple of orders of magnitude for the close evaluation problem. {All} modified representations provide a satisfactory correction overall. \edit{We use} a naïve (straightforward) implementation of \eqref{eq:SLP_laplaceN} {and \eqref{eq:general_SLP}} in Matlab, computed on a Mac mini SSD 512Go\edit{. We provide run times in Table \ref{tab:CPU_2D_slp} for various number of quadrature points. Run times do not count the} time to compute the boundary integral equation for $\rho$ (being the same for all methods). 
\begin{table}[h!] 
  \centering
 {\small 
 \begin{tabular}{|c|c|c|c|c|c|}
\hline
Method &  V0 &  V1 & V2 & V3 & V4 \\
\hline 
N = 128 & 0.014 &0.044 & 0.055 &0.045 & 0.05 \\
 \hline
 N = 256 & 0.056& 0.07 & 0.112 & 0.08 & 0.081 \\
  \hline
  N = 512 & 0.12 & 0.192  & 0.263 & 0.2 &  0.19\\
  \hline  
\end{tabular}
}
 \caption{{\textbf{Laplace 2D single-layer.} CPU times (in seconds) for various number of quadrature points and representations. Times account for computing the solution at $N \times 12$ grid points ($\ell  = 10^{-k}$, $k = \llbracket 0, 11 \rrbracket$) on a body-fitted grid.}}\label{tab:CPU_2D_slp}
\end{table}
{Representation V0} is obviously cheaper (less terms to compute) than {V1--V4}, and {V1} is cheaper than {V2--V4} due to simpler terms: there are less operations to conduct to compute {$\vv_1(y)$} than {the other provided auxiliary functions}. \\
To better compare the methods, Figure~\ref{img:SLP2D_3_new} represents log plots of the maximum error with respect to the number of quadrature points $N$ and for various distances $\ell$ from point $A$ (indicated in Figure \ref{img:SLP2D_1_new}). Results show that {modified representations allow to gain a couple of order of magnitude even for moderate $N$ ($N < 100$). Additionally, the error using V0 decreases linearly with the number of quadrature points whereas it is cubic using modified representations.}
While there is no significant difference between the {considered modified representations V1--V4}, one may consider run times {(and simplicity of auxiliary function $\vv$)} to discuss competitiveness. 
Based on above results, overall {representation V1} seems to be the best choice for the best computational cost-accuracy trade-off. 
{Let us emphasize that the focus of this paper is to highlight} the efficacy and simplicity of the proposed modified representations, given a quadrature rule. Our results show that modified representations allow to naturally gain a couple of orders of magnitude in the error, addressing the close evaluation problem even for moderate computational resources. {Additionally, the proposed auxiliary functions are independent of the density $\rho$.} In the next section we investigate the efficacy of {\eqref{eq:general_SLP}} in three dimensions. 

\subsubsection{Example 2: exterior Laplace in three dimensions}\label{ssec:3D_laplace}
\begin{figure}[H]
  \centering
 \includegraphics[width=0.32\textwidth]{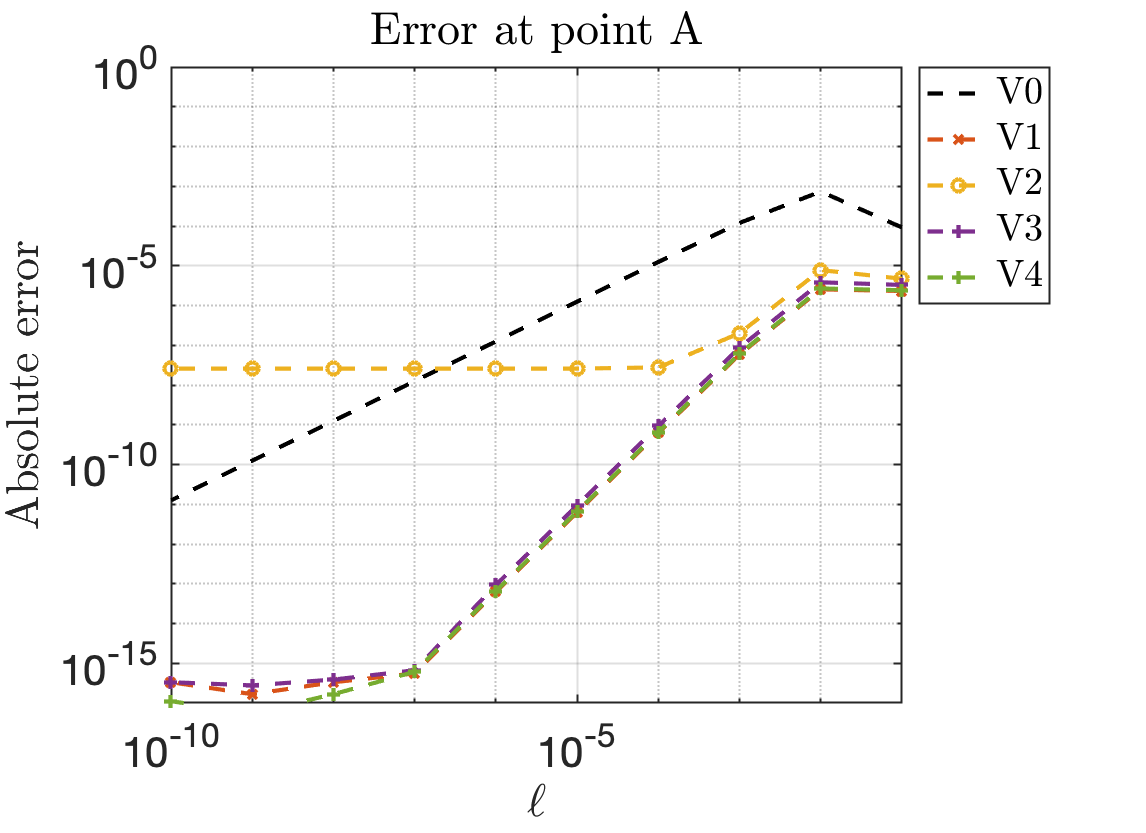}
\includegraphics[width=0.32\textwidth]{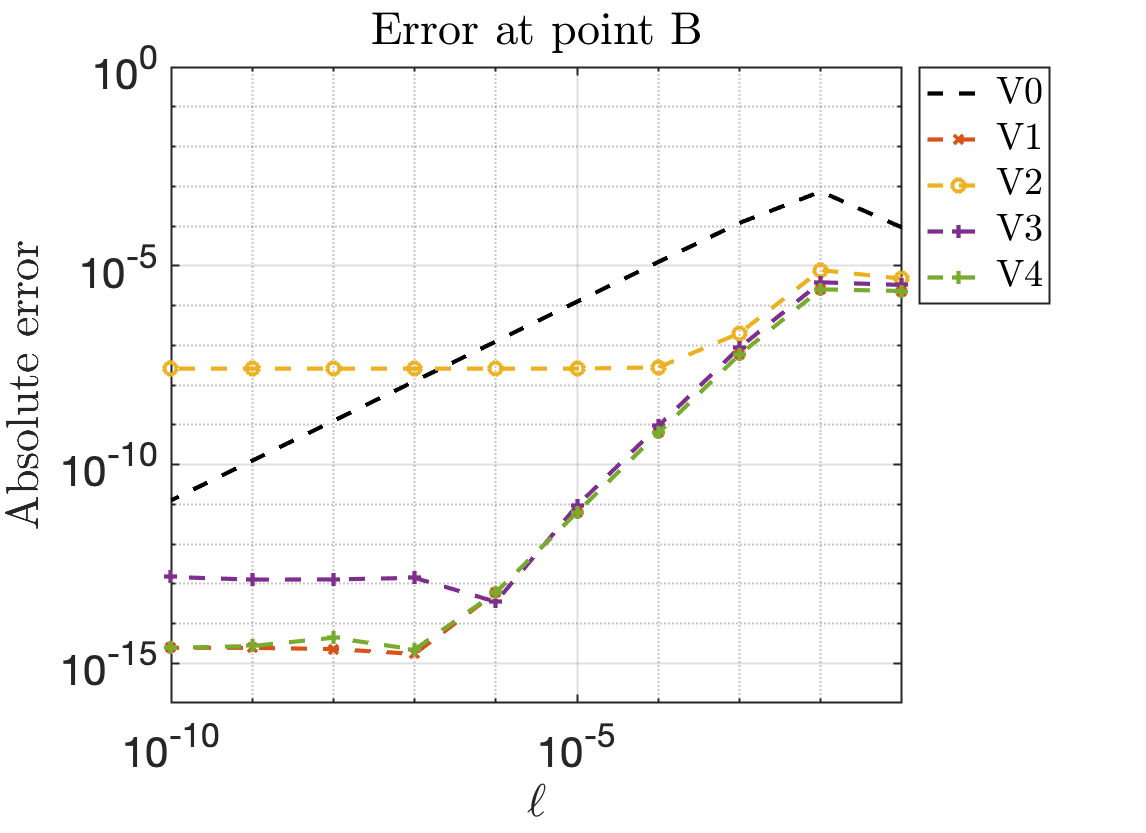}
  \caption{\small {\textbf{Laplace 3D single-layer.} Log-log plots of the errors with respect to $\ell$ made in computing the solution of \eqref{eq:exterior_bvp} for the Neumann data, $g(x^\ast) =-\frac{n_{x^\ast} \cdot (x^\ast - x_0)}{|x^\ast -x_0|^3}$ with $x_0=(0,0,0)$, outside of a sphere a radius 2, along the normal of the point A = $(-0.0065,-0.0327,1.9997)$ (left), of the point B = $(-0.3526,-1.7728,0.8561)$ (right).} }
  \label{img:SLP3D_1_new}
\end{figure}

\begin{figure}[H]
  \centering
 \includegraphics[width=0.32\textwidth]{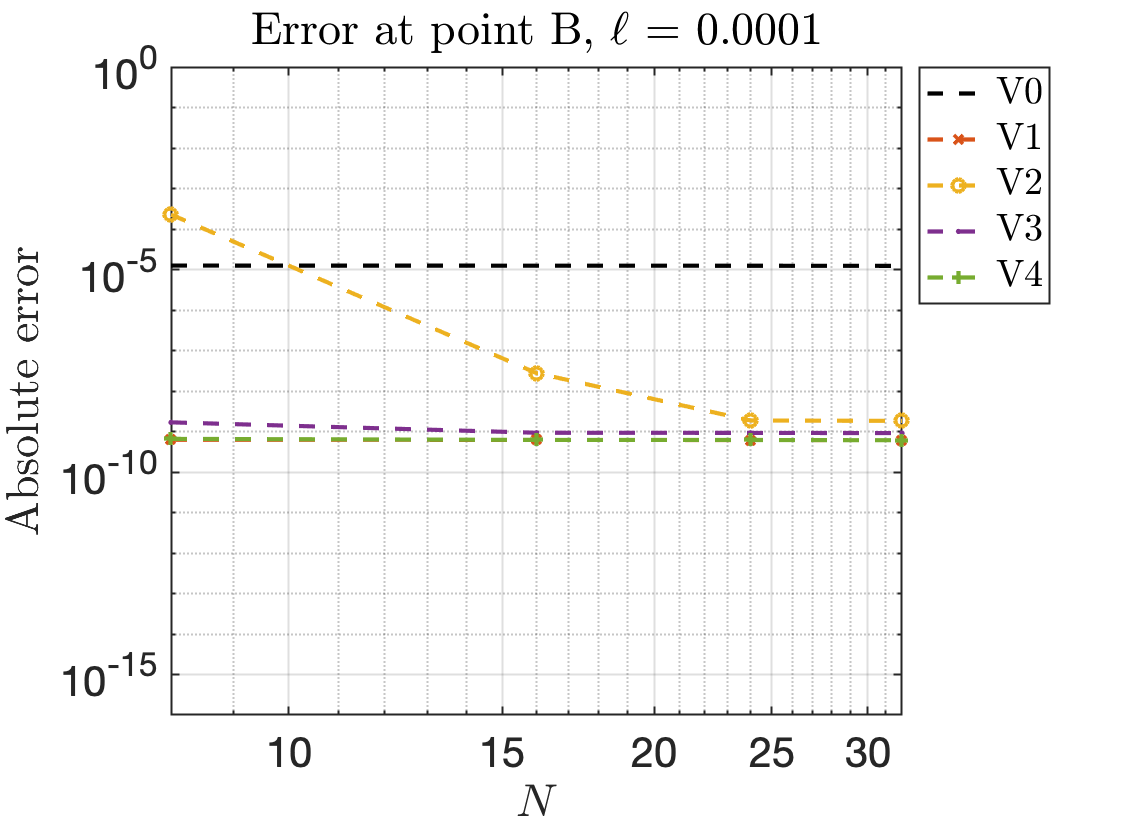}
\includegraphics[width=0.32\textwidth]{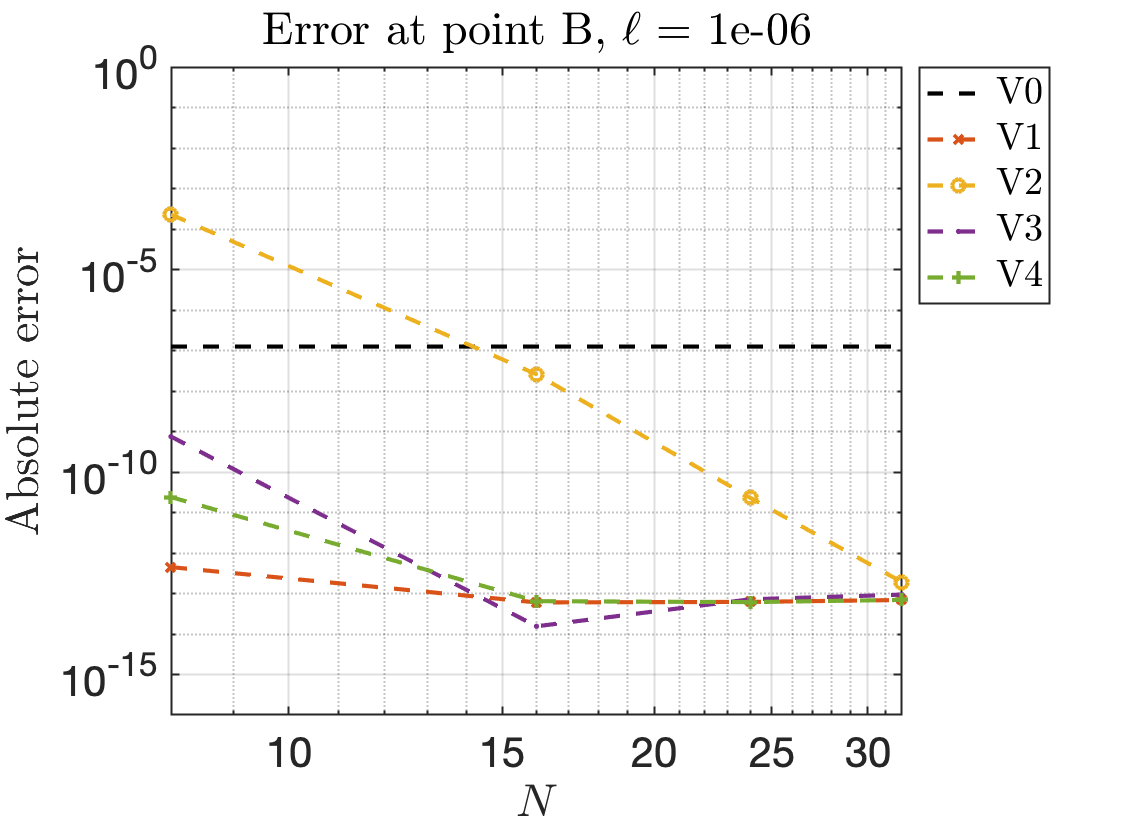}
 \includegraphics[width=0.32\textwidth]{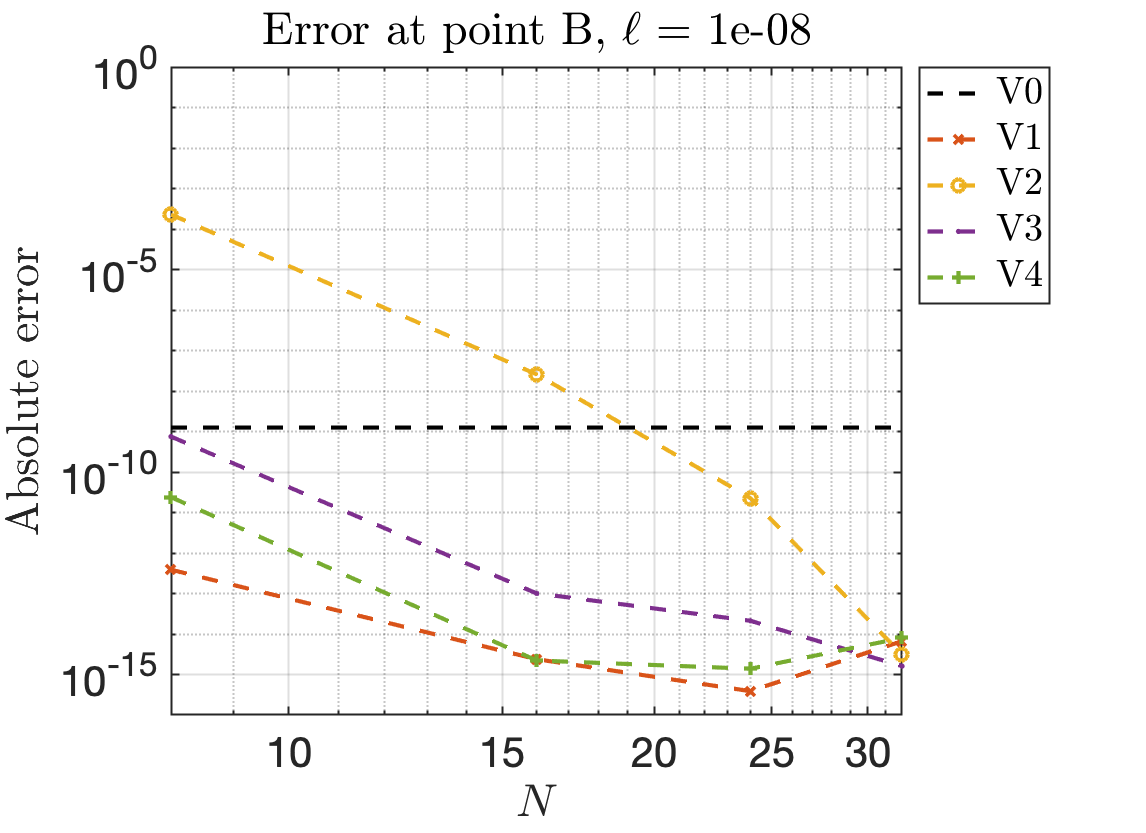}
  \caption{\small {\textbf{Laplace 3D single-layer.} Log-log plots of the errors with respect to $N$ made in computing the solution (as described in Figure \ref{img:SLP3D_1_new}) at some distance $\ell$ along the normal from the point B= $(-0.3526,-1.7728,0.8561)$.} }
  \label{img:SLP3D_2_new}
\end{figure}

{Given a \textit{domain} $D \subset \R^3$ with smooth boundary, w}e assume $\partial D$ to be an analytic, closed, and oriented surface that can be parameterized by $y = y(s,t)$ for $s \in [0,\pi]$ and $t \in [-\pi,\pi]$. Then one can write \eqref{eq:SLP_laplaceN} as
\begin{equation}\label{eq:u3D}
  u(x)=  \int_{-\pi}^{\pi}
  \int_{0}^{\pi} G^(x,y(s,t)) J(s,t) {\rho}(y(s,t)) \sin(s) \mathrm{d}s
  \mathrm{d}t,
\end{equation}
with $J(s,t) = | y_{s}(s,t) \times y_{t}(s,t) |/\sin(s)$ the Jacobian. We now work with a surface integral defined on a sphere, and we use a \textit{three-step method} (see \cite{KKCC19} for details) to approximate \eqref{eq:SLP_laplaceN} and \eqref{eq:general_SLP}. This method corresponds to a modification of the product Gaussian quadrature rule (PGQ)~\cite{atkinson1982numerical}, and it has been shown to be very effective for computing layer potentials in three dimensions at close evaluation points compared to other quadrature methods for nearly singular integrals \cite{KKCC19}. It relies on (i) rotating the local coordinate system so that $x^\ast$ corresponds to the north pole, (ii) use Periodic Trapezoid Rule with $2N$ quadrature points to approximate the integral with respect to $t$, (iii) use Gauss-Legendre with $N$ quadrature points mapped to $(0, \pi)$ (and not $(-1,1)$) to approximate the integral with respect to  $s$. This leads to the approximation:
\[u(x) \approx \frac{\pi^2}{2 N} \sum_{i = 1}^{N} \sum_{j = 1}^{2N} w_{i} \sin( s_{i} ) {F}(s_{i},t_{j}),\]
 with ${F}(s_{i},t_{j}) = G(x,y(s_i,t_j)) J(s_i,t_j) {\rho}(y(s_i,t_j)) $, $t_{j} = -\pi + \pi (j - 1)/N$, $j = 1, \cdots, 2N$, $s_{i} = \pi ( z_{i} + 1 ) / 2$, $i = 1, \cdots, N$ 
with $z_{i} \in (-1,1)$ the $N$-point Gauss-Legendre quadrature rule abscissas with corresponding weights $w_{i}$ for $i = 1, \cdots, N$. One proceeds similarly for \eqref{eq:general_SLP}. We consider an exact solution of Problem \eqref{eq:exterior_bvp}:
\[u_{\text{exact}}(x) =\frac{1}{|x - x_0|}, \quad  x_0 \in {D},\] 
which consists in choosing $g(x^\ast) = \partial_{n_{x^\ast}} u_{\text{exact}}(x^\ast)$, for any $x^\ast \in \partial D$. The efficacy of the three-step method for various geometries (including effects of curvature) is presented in \cite{KKCC19}. {Naively implementing this method has the same computational cost as the PGQ method. We do not focus in this paper on fast implementations but do believe that it is possible to speed up this method using ideas that have been previously developed including the fast multipole method~\cite{GORV21}.} Then for simplicity, results will be computed on a sphere where the resolution of $\rho$ does not require a lot of quadrature points. One can apply the technique for arbitrary closed smooth surfaces, but might be limited by the resolution of \eqref{eq:BIE_laplace_NS}. 
{
All simulations are done outside of a sphere of radius 2 using the three-step method with $N = 16$ for the following representations:
\begin{itemize}\renewcommand{\labelitemi}{$\bullet$}
\item \textbf{V0:} standard representation \eqref{eq:SLP_laplaceN};
\item \textbf{V1:} modified representation \eqref{eq:general_SLP} with the linear function $\vv_1 (y) = n_{x^\ast} \cdot y $;
\item \textbf{V2:} modified representation \eqref{eq:general_SLP} with the Green-based function\\ $\vv_2(y) =  4 \pi G(y, x^\ast + n^\ast)$; 
\item \textbf{V3:} modified representation \eqref{eq:general_SLP} with the quadratic function \\$\vv_3(y) =  \dsp\frac{1}{2} \frac{y_1 ^2 - y_2^2}{n_{x^\ast,1}x_1^\ast  - n _{x^\ast,2}x_2^\ast } $; 
\item \textbf{V4:} modified representation \eqref{eq:general_SLP} with the quadratic product function\\ $\vv_4(y) =  \dsp \frac{(y_1 - 5)(y_2 - 5)}{n_{x^\ast,1}(x_2^\ast - 5) +n _{x^\ast,2}(x_1^\ast - 5)}$ .
\end{itemize} 
Note that there are other quadratic polynomials $\vv$ (as a function of 2 variables instead of 3, see \cite{PeFaTu19} where those polynomials serve as basis for interpolation method). We make here the choice to test using similar functions as in Section \ref{ssec:Laplace2D}.}
We solve \eqref{eq:BIE_laplace_NS} using a Galerkin method and the product Gaussian quadrature rule~\cite{atkinson1982numerical,atkinson1982laplace,atkinson1985algorithm,atkinson1990survey,atkinson1997numerical} (see Appendix \ref{sec:appdx3D} for details)\edit{. The accuracy in approximating V0-V4} is limited by the accuracy of the resolution for $\rho$. {This can be assessed by looking at the coefficients' decay of the density spherical harmonic expansion. In this case the coefficients' decay has reached $10^{-15}$.} Results in Figure \ref{img:SLP3D_1_new} show that given $\rho$ resolved, the approximation of the modified representations provide better results overall{, except for V2 where the error plateaus around $10^{-7}$ for small $\ell$ (providing less accurate results compared to standard representation V0)}. Note that the single-layer potential commonly suffers less from the close evaluation than the double-layer potential, and the chosen method provides already a good approximation. {This is the reason why the error when considering V0 decays as $\ell$ decreases~\cite{KKCC19}.} The modified representations allow to make it even better.
To better assess the efficacy of the modified representations in three dimensions, Figure~\ref{img:SLP3D_2_new} represents log plots of the maximum error with respect to $N \in \{8, 16, 24, 32\}$ (the method uses $2N \times N$ quadrature points) and \edit{for various distances $\ell$} from the point B. Results show that as $\ell \to 0$, {V1--V4} allow to gain {a couple of} orders of magnitude in the error, even for a small $N$. {Note that the error produced by three-step method doesn't seem to depend on $N$, and in this case there are more variations with respect to the choice of auxiliary function $\vv$ than in two dimensions. Here, V1 (the linear function) is the best representation, producing the smallest errors (and the fastest to compute as indicated in Table \ref{tab:CPU_3D_slp}).} Again, the three-step method has been designed to treat nearly-singular integrals\edit{. I}t is the reason why the {method} provides already satisfactory results (given the resolution of $\rho$). The modified representations allow to significantly gain even more accuracy in this case, even with limited computational resources. 
\begin{table}[h!] 
  \centering
 {\small 
\begin{tabular}{|c|c|c|c|c|c|}
\hline
Method &   V0 & V1 & V2 & V3 & V4 \\
 \hline
N = 8 &0.028 &0.029 &0.032 &  0.031 & 0.046 \\
 \hline 
N = 16 & 0.143 & 0.146& 0.148 & 0.150 & 0.142  \\
 \hline 
N = 24 &  0.352 &  0.344 &  0.346 &0.35 &  0.356 \\
 \hline
\end{tabular}
}
 \caption{{\textbf{Laplace 3D single-layer.} CPU times (in seconds) for various number of quadrature points and representations for computing the solution (as described in Figure \ref{img:SLP3D_1_new}) from the points A and B, for $\ell  = 10^{-k}$, $k = \llbracket 0, 11 \rrbracket$.}}\label{tab:CPU_3D_slp}
\end{table}
\subsection{Scattering problem}
Using Proposition \ref{pro:helm}, we compare \eqref{eq:DLP-SLP_helm} with the modified representation \eqref{eq:general_helm} obtained with $\vv (y)= e^{ik n_{x^\ast} \cdot (y-x^\ast)}$: 
\begin{equation}\label{eq:modified_DLP-SLP_helm}
\begin{aligned}
u(x)    = & \dsp \int_{\partial D}   \left[ \partial_{n_y} G^H(x,y) - i k (n_y \cdot n_{x^\ast})e^{i k (n_{x^\ast} \cdot (y- x^\ast))} G^H(x,y) \right] \left [\mu (y)  - \mu(x^\ast)  \right] \, d \sigma_y  \\
& +  ik    \dsp \int_{\partial D}   [(n_y \cdot n_{x^\ast})e^{i k (n_{x^\ast} \cdot (y- x^\ast))}   -1 ]  G^H(x,y)  \mu (y) \, d \sigma_y \\
 &  +\mu(x^\ast)  \dsp \int_{\partial D}  \partial_{n_y} G^H(x,y) [1-e^{i k (n_{x^\ast} \cdot (y- x^\ast))}   ] \, d \sigma_y, \quad  x \in \mathbb{R}^d \setminus \bar{D}.
 \end{aligned}
\end{equation}
\subsubsection{Example 3: scattering in two dimensions}\label{ssec:helm2d}
We consider an exact solution of Problem \eqref{eq:scattering_pb}:
\[u_{\text{exact}}(x) =  \frac{i}{4} H^{(1)}_0(k |x - x_0|), \quad  x_0 \in {D},\] 
which consists in choosing $f(x^\ast) = u_{\text{exact}}(x^\ast)$, for any $x^\ast \in \partial D$. 
{
All simulations are done outside of a star-shaped domain using the Periodic Trapezoid Rule with $N = 256$ quadrature points and $k=15$ for the following representations:
\begin{itemize}\renewcommand{\labelitemi}{$\bullet$}
\item \textbf{V0:} standard representation \eqref{eq:DLP-SLP_helm};
\item \textbf{V1:} modified representation \eqref{eq:modified_DLP-SLP_helm} (i.e. \eqref{eq:general_helm} with the plane wave function $\vv_1 (y)= e^{ik n_{x^\ast} \cdot (y-x^\ast)}$).
\end{itemize} 
}
We solved \eqref{eq:BIE_helm} using Kress product quadrature rule~\cite{Kress91} (see Appendix \ref{sec:appdx2D}). The quadrature rule is well adapted to approximate kernels with a logarithmic singularity. The accuracy of both methods is limited by the resolution for $\mu$ {(the Fourier coefficients' decay of the density is bounded by $10^{-6}$ for $N = 256$ and $k=15$)}. Results in Figures \ref{img:helmholtz2D_1_1} and \ref{img:helmholtz2D_1_2} show that given $\mu$ resolved, the approximation of the modified representation provides better results overall. Similarly to Laplace's examples, both methods approximate well the solution far from the boundary. As the evaluation point gets closer to the boundary ($\ell \to 0$), {V0 approximated with PTR} suffers from the close evaluation problem leading to large errors {(see \cite{barnett2014evaluation})}. Using the modified representation {V1} allows to reduce the error by a couple of order of magnitude for the close evaluation problem.
{
\begin{figure}[h!]
  \centering
  \includegraphics[width=0.42\textwidth]{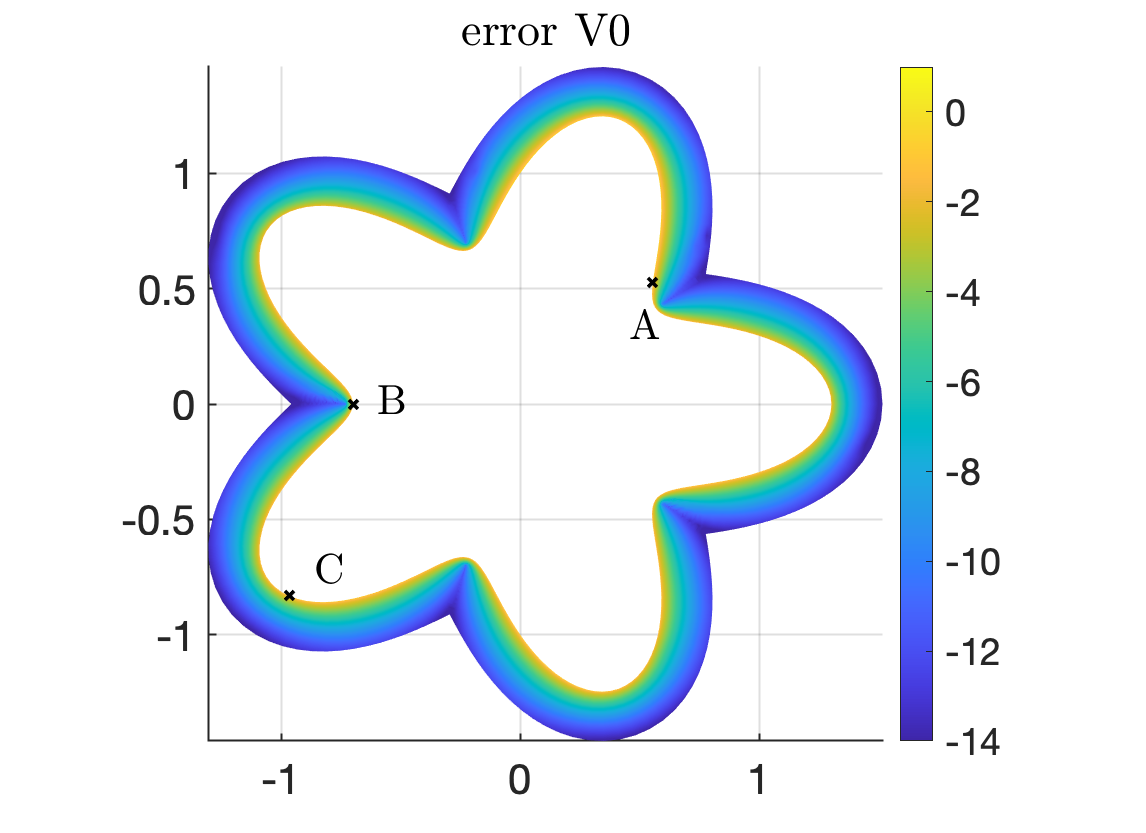}
    \includegraphics[width=0.42\textwidth]{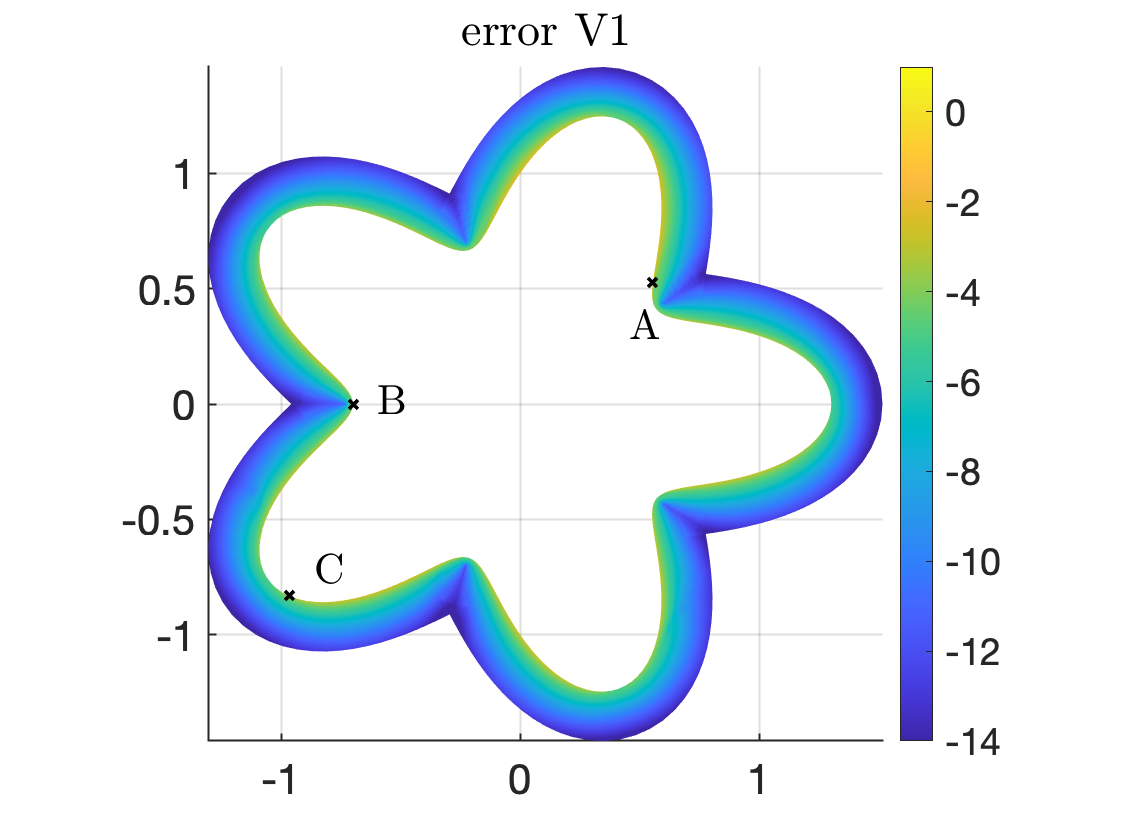}
  \caption{\small {\textbf{Helmholtz 2D.} Plots of $\log_{10}$ of the error for the evaluation of the solution of \eqref{eq:scattering_pb} out of the star domain defined by the boundary $y(t) =(1+ 0.3 \cos 5t) * (\cos t, \sin t)$, $t \in [0 ,2\pi]$, for the Dirichlet data, $f(x^\ast) = \frac{i}{4} H^{(1)}_0(15 |x^\ast - x_0|)$ with $x_0=(0.2,0.8)$, for representations V0, V1, computed using 
    PTR with $N = 256$.}}
  \label{img:helmholtz2D_1_1}
\end{figure}}
{
\begin{figure}[!hbt]
  \centering
  \includegraphics[width=0.32\textwidth]{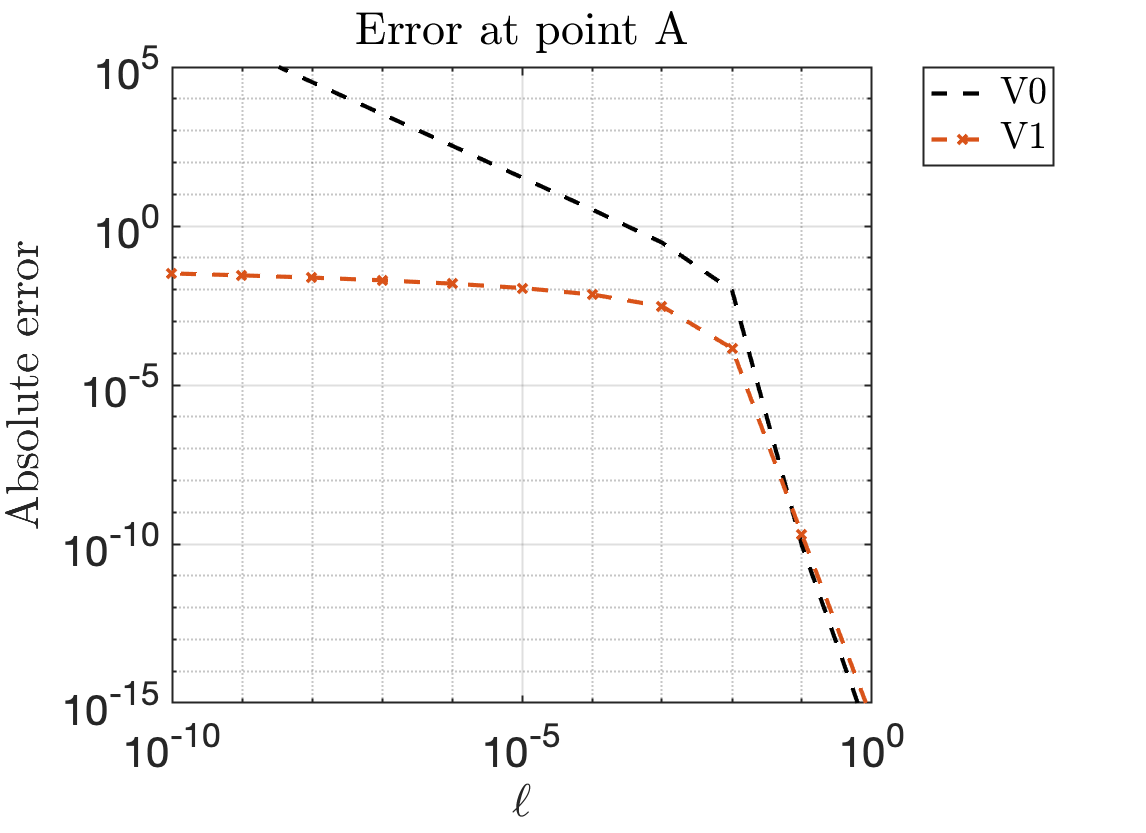}
  \includegraphics[width=0.32\textwidth]{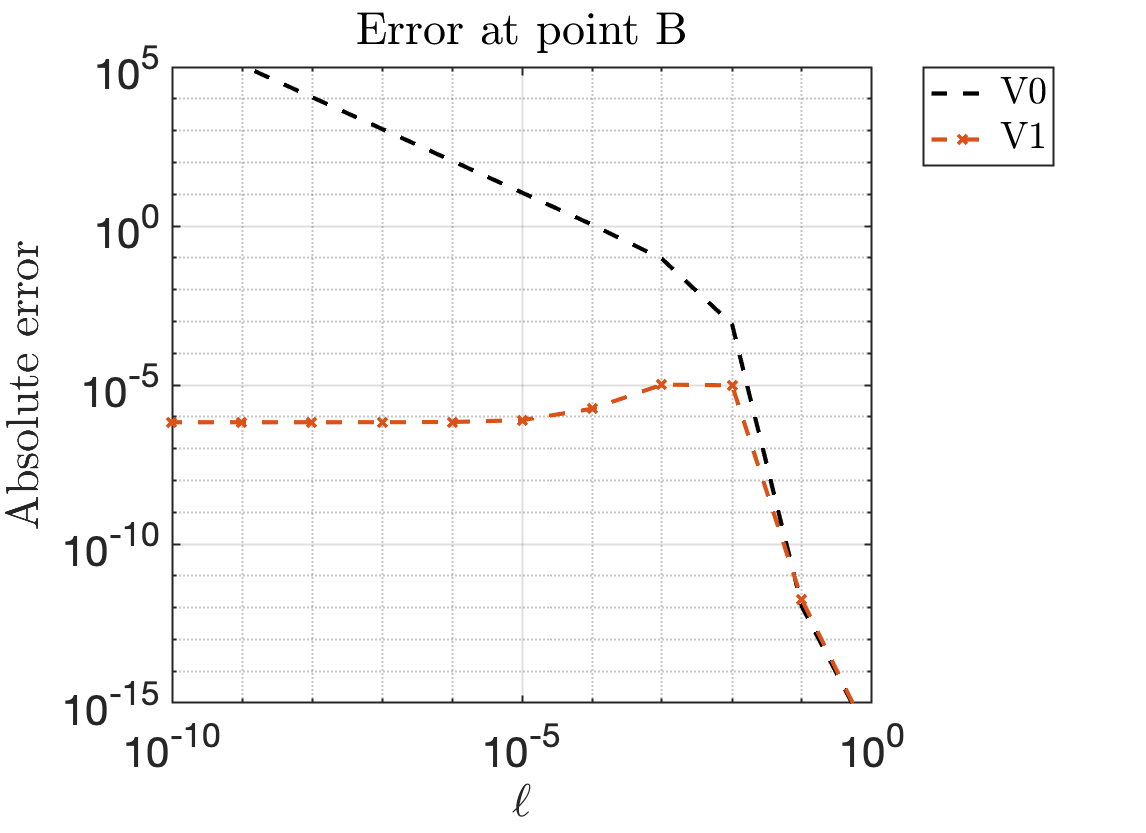}
    \includegraphics[width=0.32\textwidth]{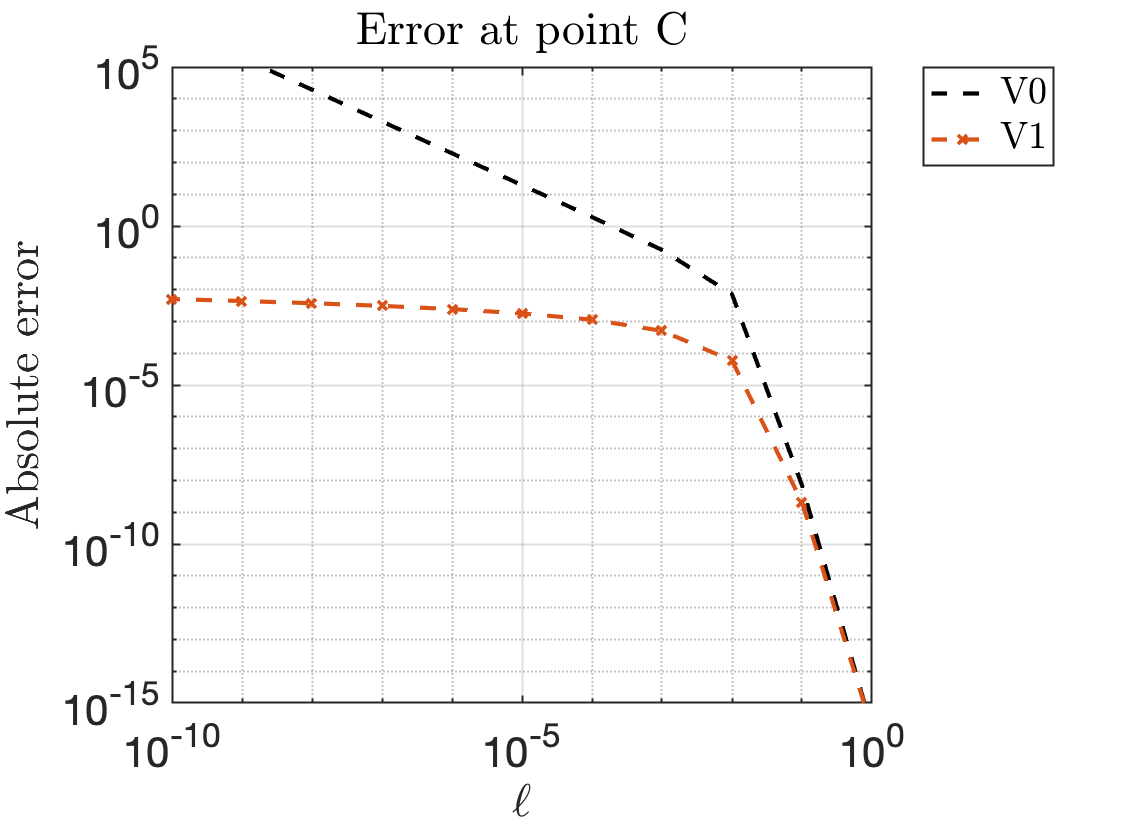}
  \caption{\small {\textbf{Helmholtz 2D.} Log-log plots of the errors made in computing the solution along the normal of the three points A, B, C, plotted as black $\times$'s in Figure \ref{img:helmholtz2D_1_1}. }}
  \label{img:helmholtz2D_1_2}
\end{figure}
\begin{figure}[h!]
  \centering
 \includegraphics[width=0.32\textwidth]{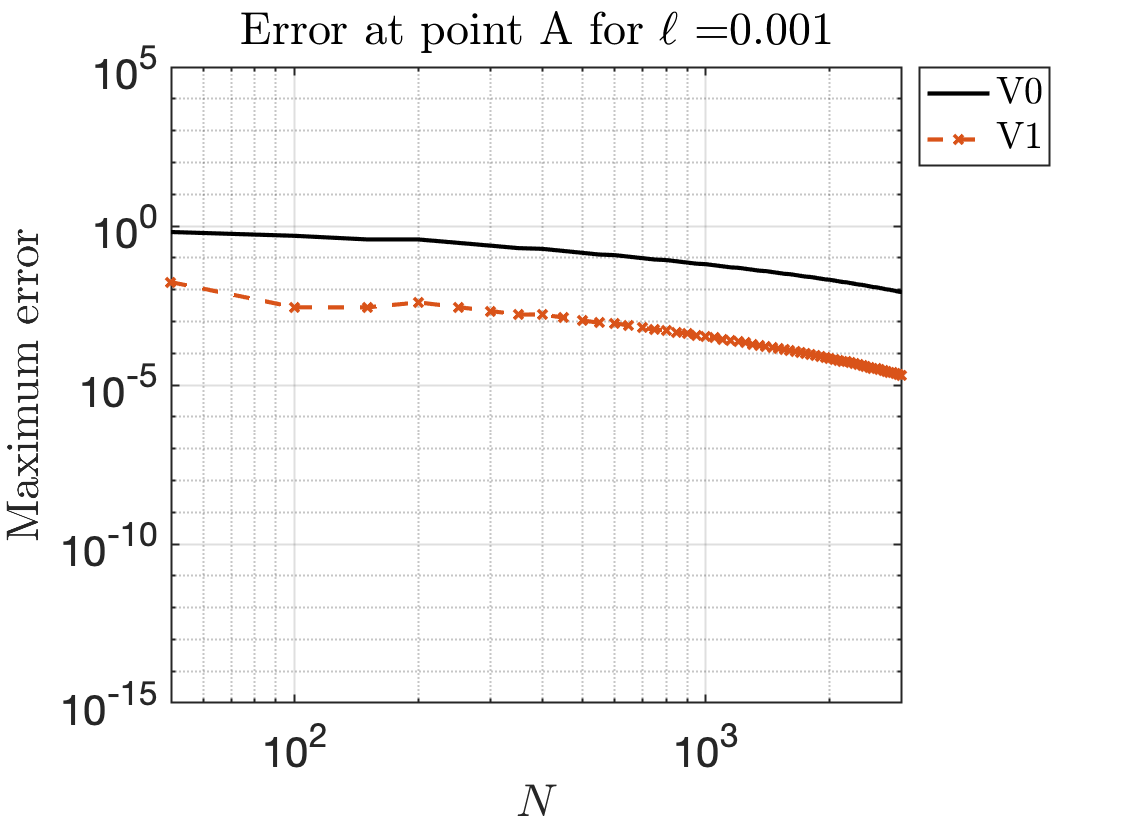}
\includegraphics[width=0.32\textwidth]{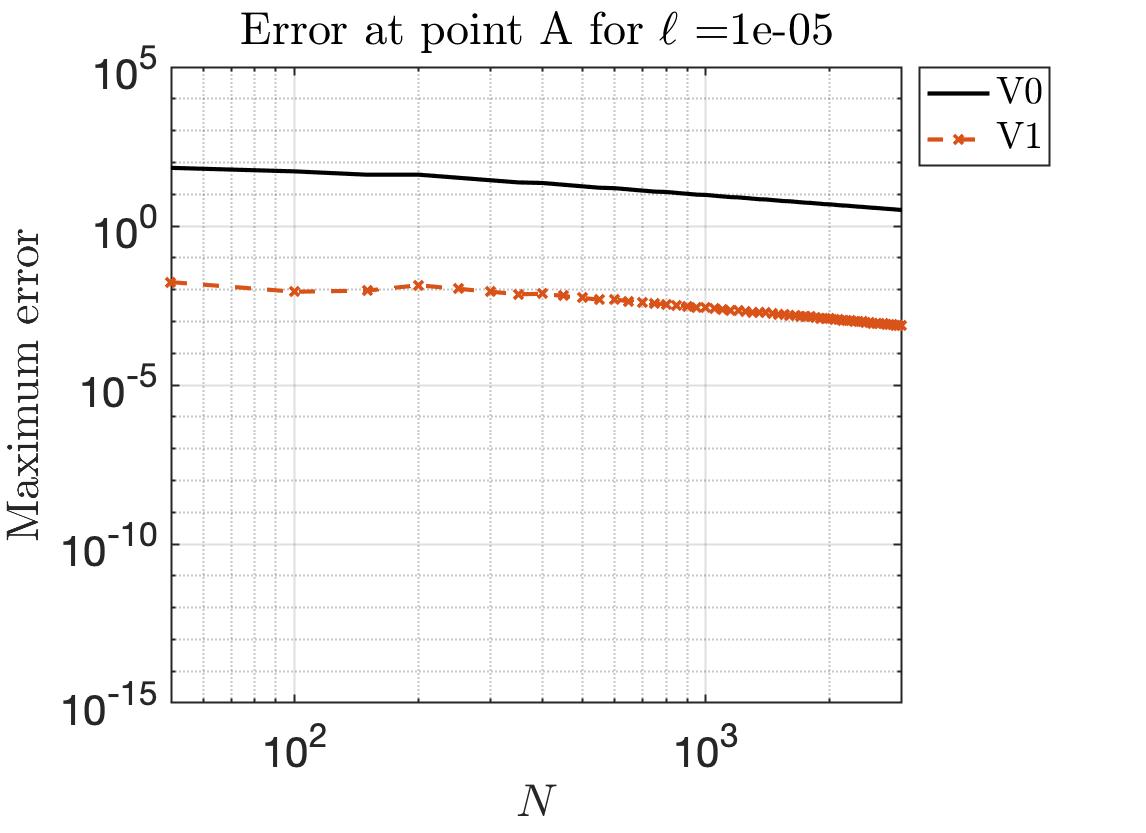}
 \includegraphics[width=0.32\textwidth]{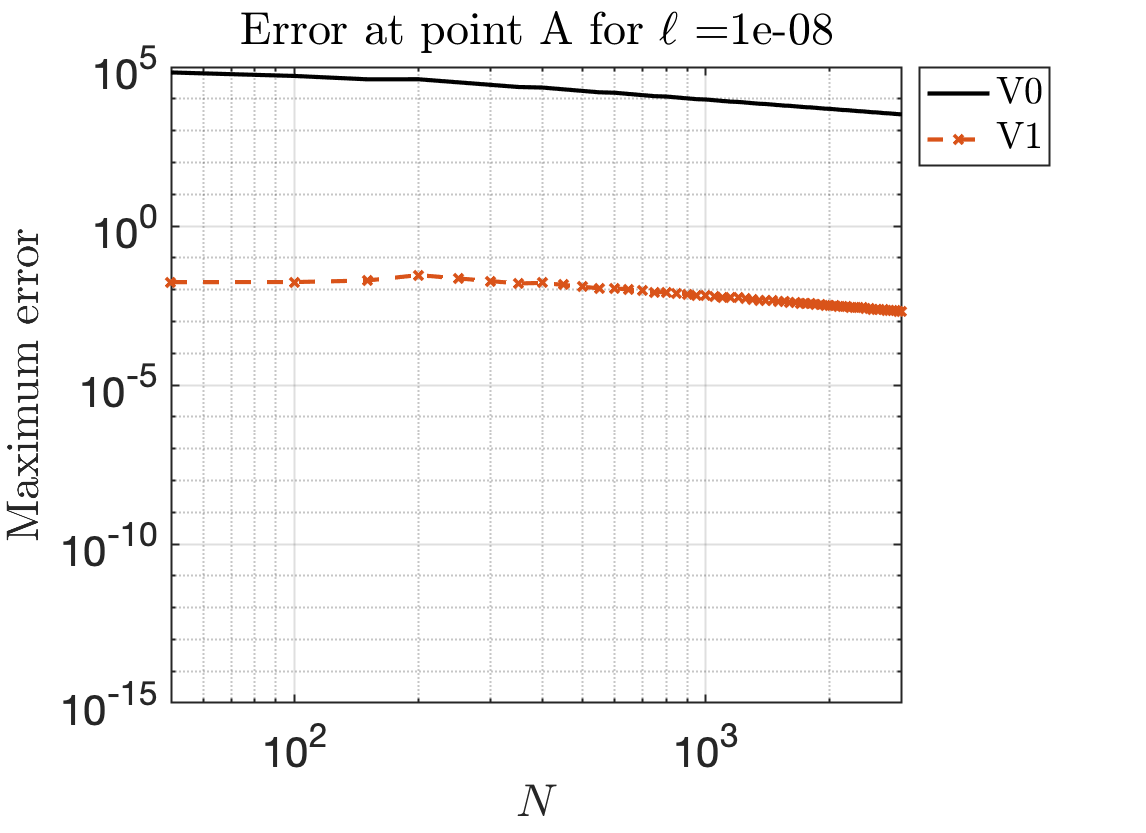}
  \caption{\small {\textbf{Helmholtz 2D.} Log-log plots of the errors with respect to $N$ made in computing the solution at some distance $\ell$ along the normal from the point A plotted as black $\times$'s in Figure \ref{img:helmholtz2D_1_1}.} }
  \label{img:helmholtz2D_1_3}
\end{figure}
}
Figure~\ref{img:helmholtz2D_1_3} represents log plots of the maximum error with respect to the number of quadrature points {$N \in \llbracket 50, 3000\rrbracket$ \edit{and} for various distances $\ell$ from point A (indicated in Figure \ref{img:helmholtz2D_1_1})}. Results show that for any number of quadrature points, the error {when considering V0} explodes as we approach the boundary (error {larger than $10^5$}) while the error with {V1} remains bounded (of the order of $10^{-2}$ {in the case presented above}). In this case standard {rerpresentation V0} strongly suffers from the close evaluation problem, however {the modified representation V1} significantly reduces the error. Even when standard quadrature rules fail to compute the standard representation, the proposed modified one regularizes the solution and provides satisfactory results {without significant additional computational time (as shown in Table \ref{tab:CPU_2D_helm}).
\begin{table}[h!] 
  \centering
 {\small 
\begin{tabular}{|c|c|c|c|}
\hline
Method &  N = 128 & N = 256 & N = 512\\ 
 \hline
	V0 & 0.18 & 0.27 & 0.71\\ 
 \hline
V1 & 0.21 & 0.33 & 0.89 \\
 \hline
\end{tabular}
}
 \caption{{\textbf{Helmholtz 2D.} CPU times (in seconds) for various number of quadrature points and representations. Times account for computing the solution for $N \times 12$ grid points (for $\ell  = 10^{-k}$, $k = \llbracket 0, 11 \rrbracket$) on a body-fitted grid.}}\label{tab:CPU_2D_helm}
 \vspace{-0.5cm}
\end{table}
}
\subsubsection{Example 4: scattering in three dimensions}\label{ssec:helm3d}
We consider an exact solution of \eqref{eq:interior_pb}:
\[u_{\text{exact}}(x) =  \frac{1}{4 \pi} \frac{e^{i k |x - x_0|}}{|x- x_0|}, \quad  x_0 \in {D},\] 
which consists in choosing $f(x^\ast) = u_{\text{exact}}(x^\ast)$, for any $x^\ast \in \partial D$. 
{
All simulations are done outside of an ellipsoid parameterized by $y(s,t) = (2 \cos(t) \sin (s), \sin(t) \sin (s), 2\cos (s))$, $(s,t) \in [0, \pi] \times [-\pi, \pi]$, and using the three-step method with various $N$. This is in order to investigate the technique in the context of limited resolution, namely the coefficients' decay of the density spherical harmonic expansion doesn't reach machine precision. We consider $k = 5$ and the following representations:
\begin{itemize}\renewcommand{\labelitemi}{$\bullet$}
\item \textbf{V0:} standard representation \eqref{eq:DLP-SLP_helm};
\item \textbf{V1:} modified representation \eqref{eq:modified_DLP-SLP_helm}.
\end{itemize} 
\begin{figure}[h!]
    \centering
    \begin{subfigure}[b]{0.32\columnwidth}
        \centering
          \includegraphics[width=\textwidth]{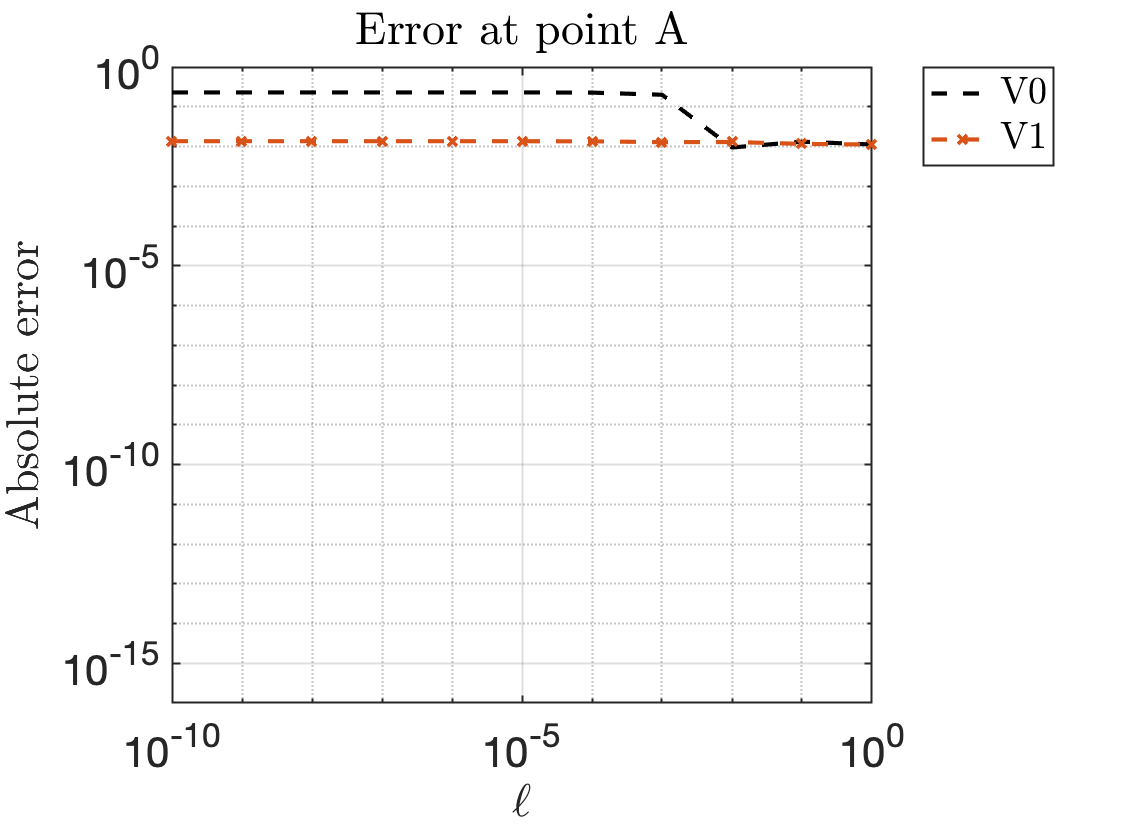}
                  \caption{N = 16}
    \end{subfigure}
        \begin{subfigure}[b]{0.32\columnwidth}
        \centering
          \includegraphics[width=\textwidth]{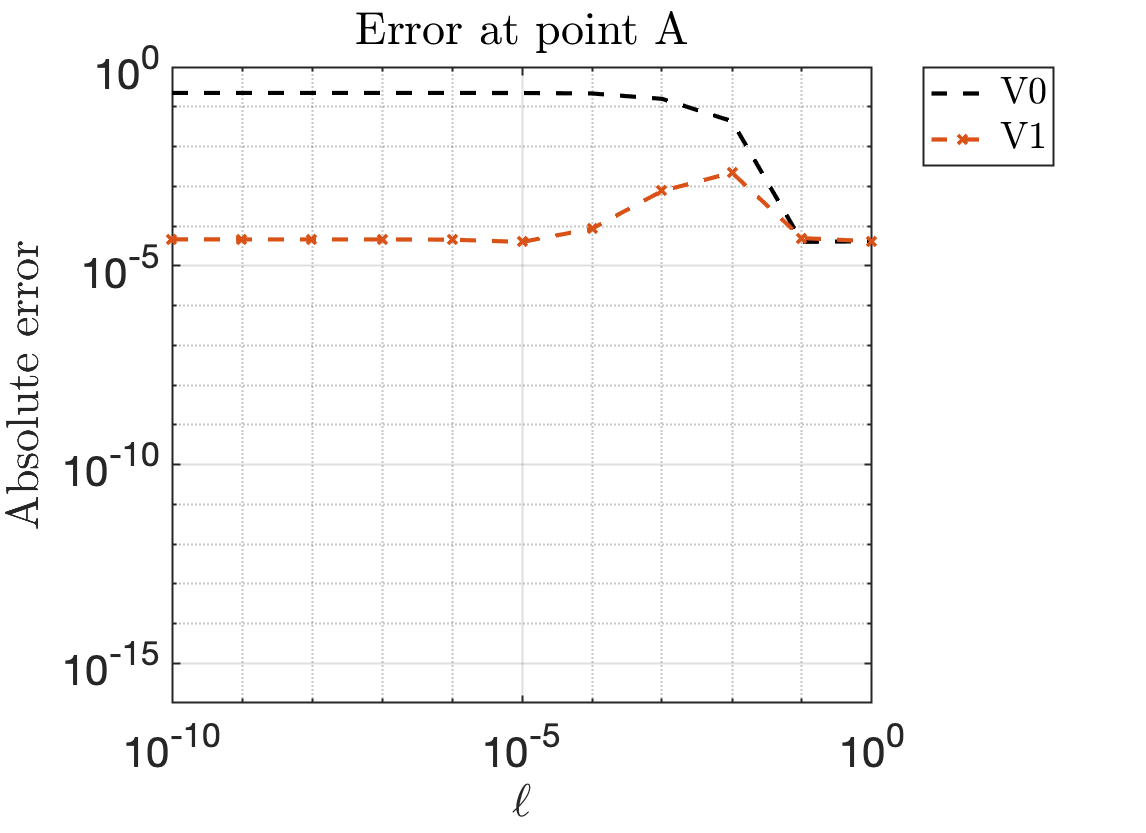}
                  \caption{N = 24}
    \end{subfigure}
        \begin{subfigure}[b]{0.32\columnwidth}
        \centering
          \includegraphics[width=\textwidth]{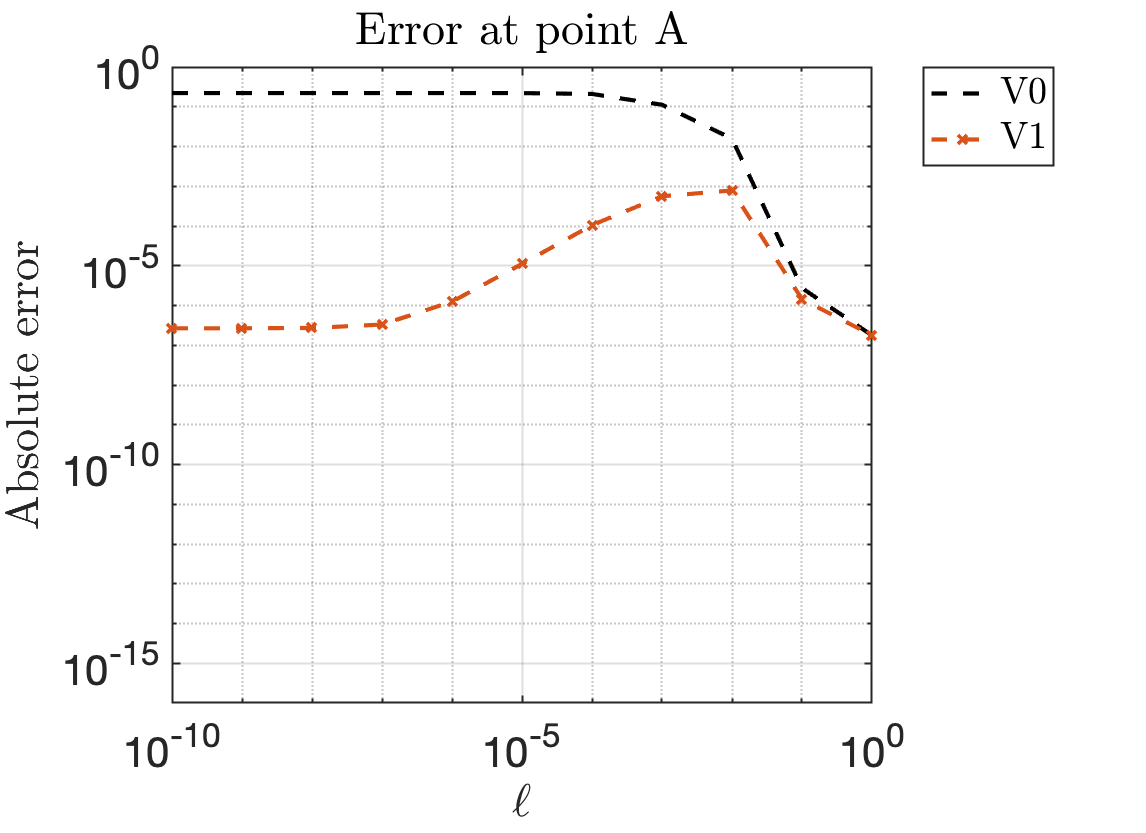}
                  \caption{N = 32}
    \end{subfigure}\vspace{0.3cm}
    \begin{subfigure}[b]{0.32\columnwidth}
        \centering
          \includegraphics[width=\textwidth]{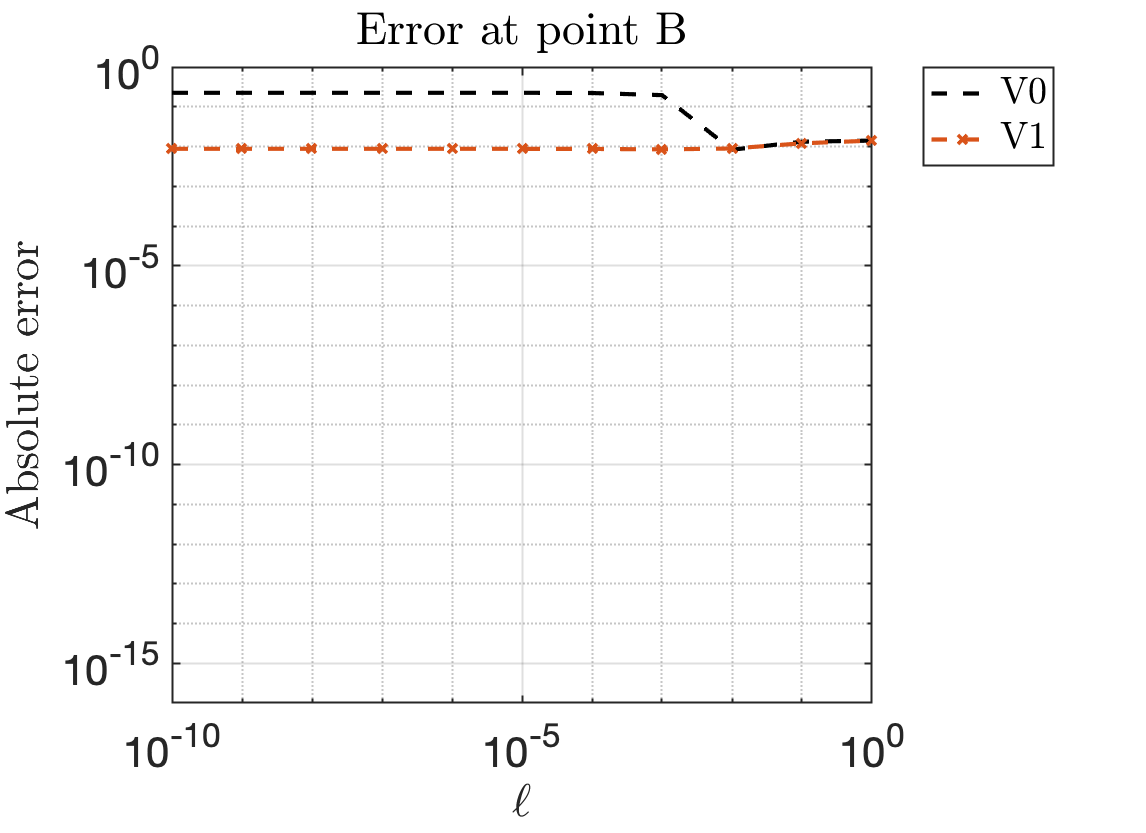}
                  \caption{N = 16}
    \end{subfigure}
        \begin{subfigure}[b]{0.32\columnwidth}
        \centering
          \includegraphics[width=\textwidth]{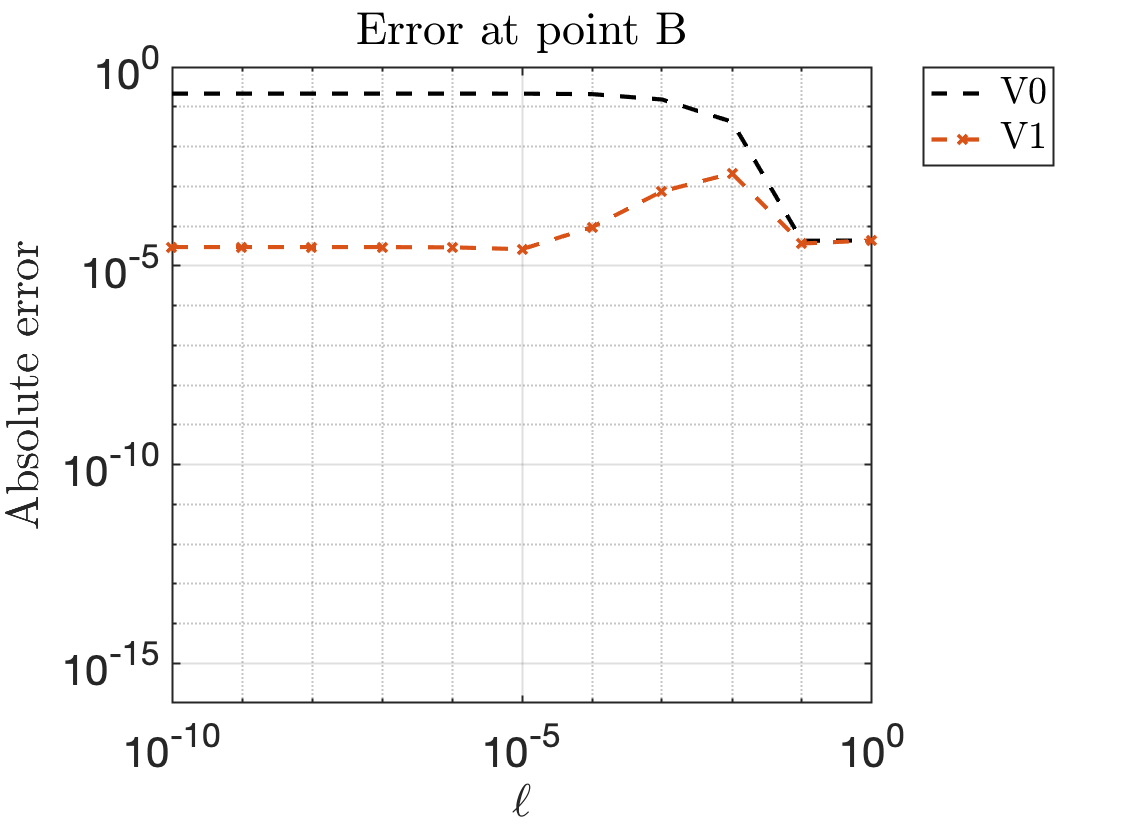}
                  \caption{N = 24}
    \end{subfigure}
        \begin{subfigure}[b]{0.32\columnwidth}
        \centering
          \includegraphics[width=\textwidth]{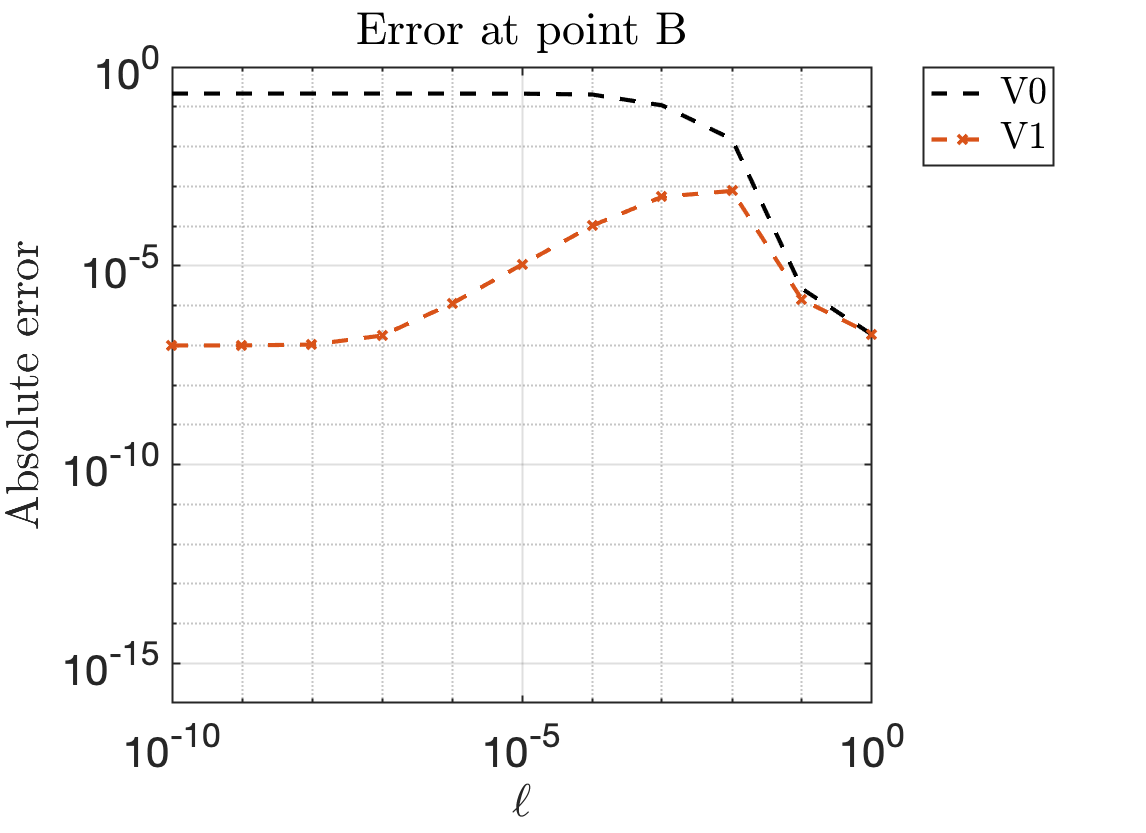}
                  \caption{N = 32}
    \end{subfigure}
  \caption{\small {\textbf{Helmholtz 3D.}} Log-Log of the error along the normal for the evaluation of the solution of \eqref{eq:scattering_pb} out of the ellipsoid parameterized by $y(s,t) = (2 \cos(t) \sin (s), \sin(t) \sin (s), 2\cos (s))$, $(s,t) \in [0, \pi] \times [-\pi, \pi]$, for the Dirichlet data $f(x^\ast) = \frac{1}{4} \frac{e^{i 5 |z -x_0|}}{|x- x_0|}$ with $x_0=(0.1,0.2,0.3)$: at the point A = $(-0.7664,0.0607,1.8433)$ {(top row)}, at the point B = $(-0.0098,-0.0096,1.9999)$ {(bottom row), for various $N$.}}
  \label{img:helmholtz3D_1}
\end{figure}
}
We solved \eqref{eq:BIE_helm} using Galerkin method and the product Gaussian quadrature rule (see Appendix \ref{sec:appdx3D} for details). The accuracy of both methods is limited by the accuracy of the resolution for $\mu$. This limitation can be checked for instance by looking at the density spherical harmonics coefficients' decay: for {$k=5$, the resolution will be capped around $10^{-2}$ for $N = 16$, $10^{-4}$ for $N = 24$, and $10^{-7}$ for $N = 32$}. Results in Figure \ref{img:helmholtz3D_1} show that given $\mu$ resolved, standard representation incurs bigger errors at close evaluation points while the modified representation provides better results overall. Here, the resolution of the boundary integral equation was fairly limited. Figure~\ref{img:helm3D_2_ellipsoid} represents log plots of the maximum error with respect to $N \in \llbracket 8, 32\rrbracket$ (the method uses $2N \times N$ quadrature points) \edit{and} {for various distances $\ell$} from the boundary {from point A}. While the three-step method has been designed to treat nearly-singular integrals and provided satisfactory results for Laplace's problems, the method here requires more quadrature points to achieve accuracy due to the wavenumber (see Section \ref{ssec:freq} for more details). {\edit{T}he standard representation V0 suffers from both the close evaluation problem and the poor density resolution. The modified representation {V1} allows to gain accuracy even with limited resolution (without significant additional computational time as indicated in Table \ref{tab:CPU_3D_helm}).}
{
\begin{figure}[h!]
  \centering
 \includegraphics[width=0.32\textwidth]{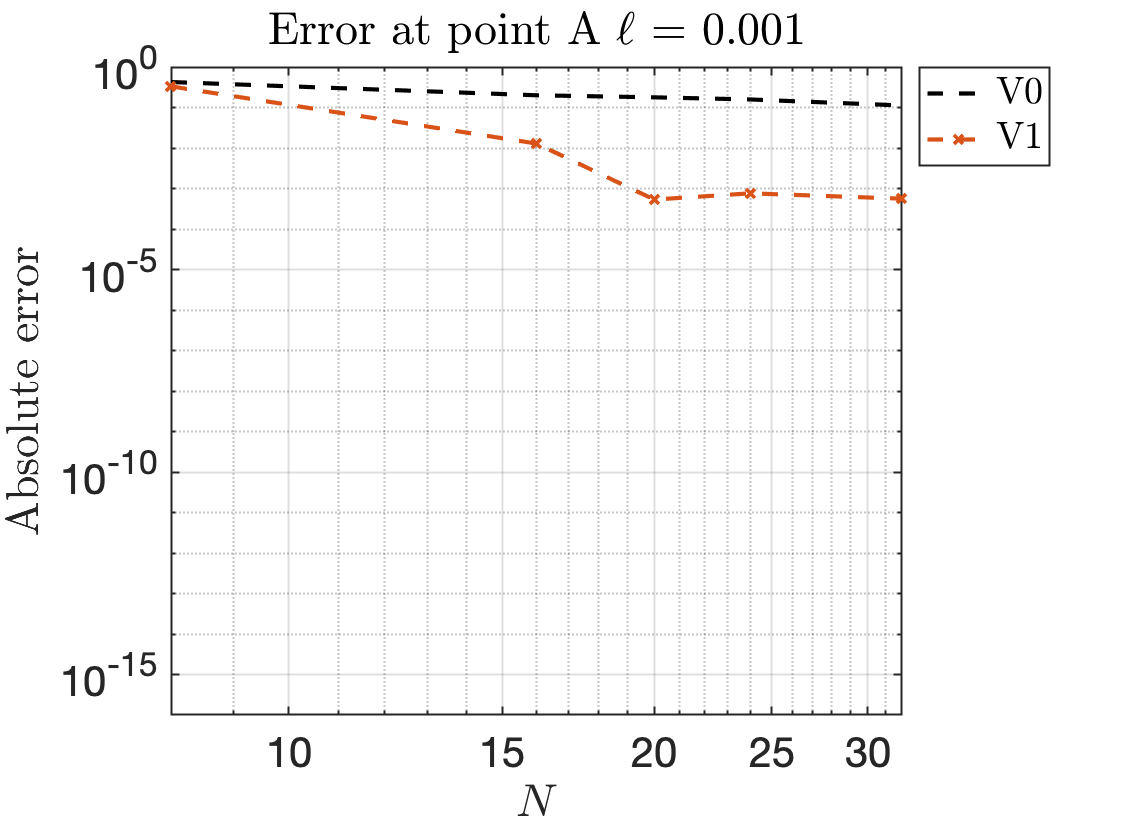}
\includegraphics[width=0.32\textwidth]{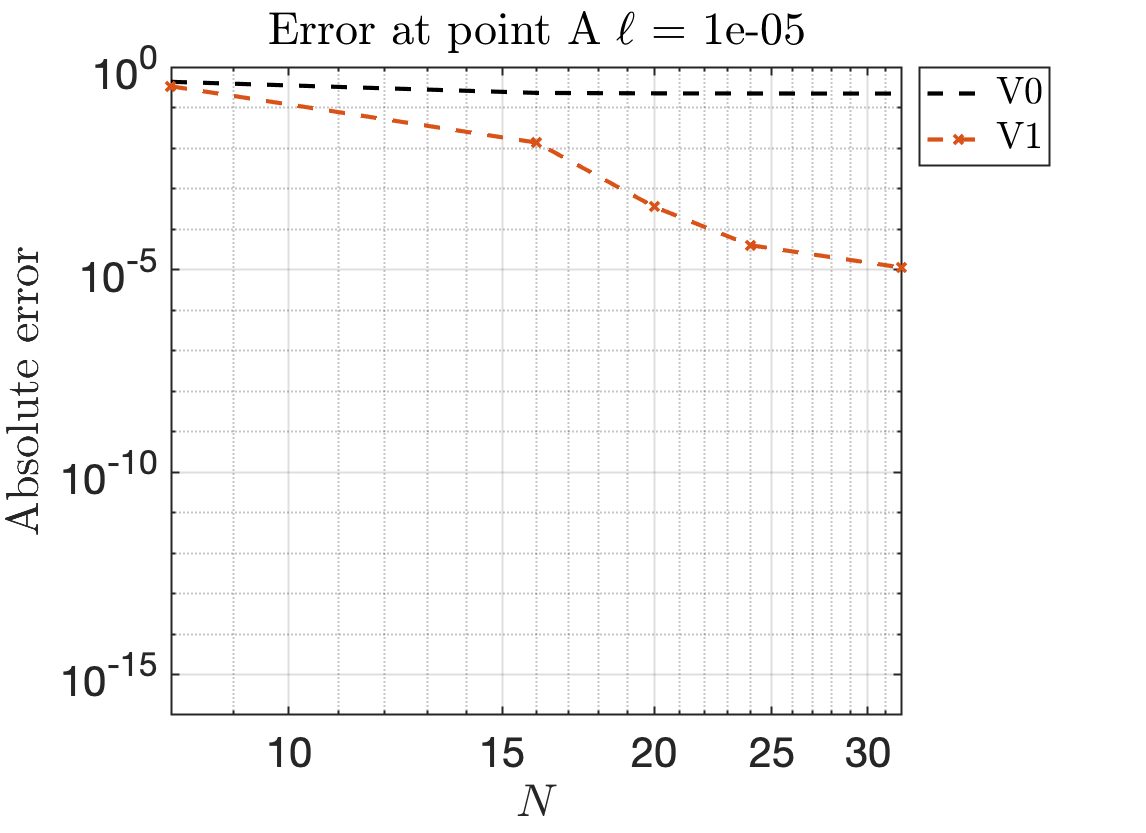}
 \includegraphics[width=0.32\textwidth]{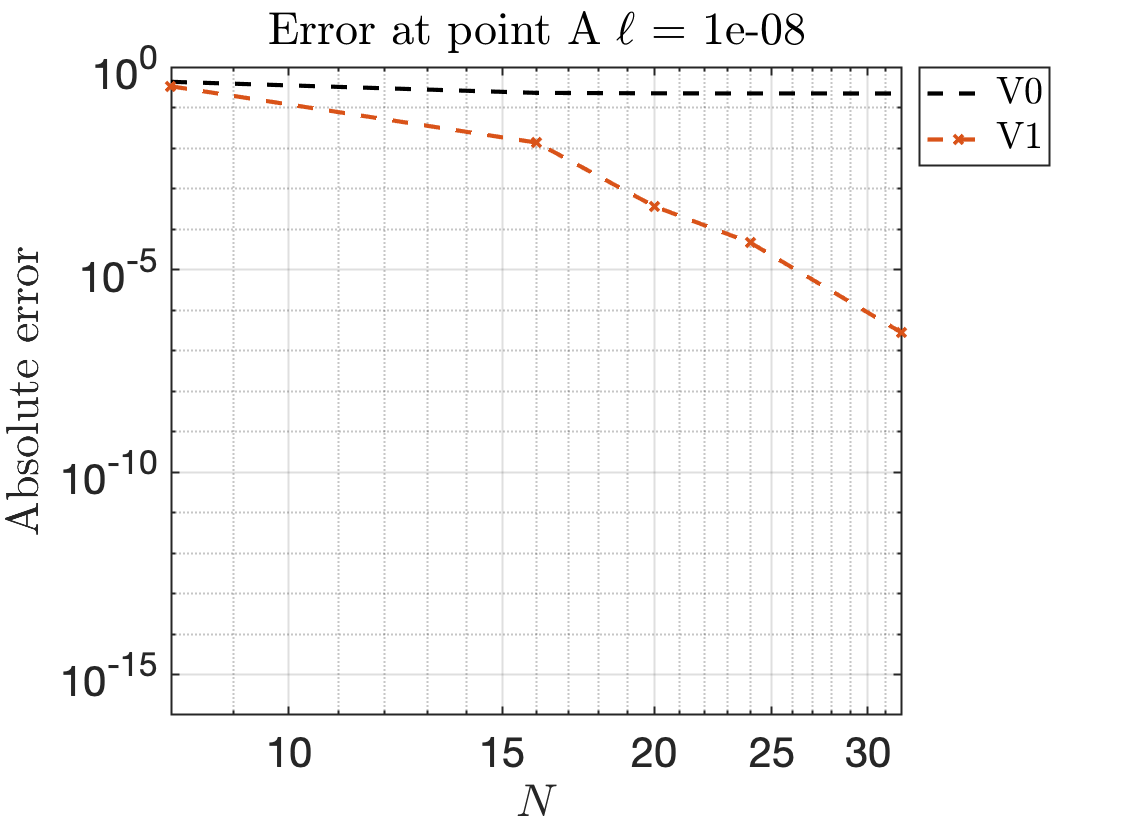}
  \caption{\small {\textbf{Helmholtz 3D.}} {Log-plot of the maximum error for computing the solution as described in Figure~\ref{img:helmholtz3D_1} with $\partial D$ being the ellipsoid parameterized by $y(s,t) = (2 \cos(t) \sin (s), \sin(t) \sin (s), 2\cos (s))$, $(s,t) \in [0, \pi] \times [-\pi, \pi]$, at some distance $\ell$ along the normal from the point A= $(-0.7664,0.0607,1.8433)$.} }
  \label{img:helm3D_2_ellipsoid}
\end{figure}
\begin{table}[h!] 
  \centering
 {\small 
\begin{tabular}{|c|c|c|c|}
\hline
Method &  N = 8 & N =16 & N = 20\\ 
 \hline
	V0 & 0.027 & 0.15 & 0.313\\ 
 \hline
V1 & 0.03 & 0.15 & 0.314 \\
 \hline
\end{tabular}
}
 \caption{{\textbf{Helmholtz 3D.} CPU times (in seconds) for various number of quadrature points and representations. Times account for computing the solution from points A and B, for $\ell  = 10^{-k}$, $k = \llbracket 0, 11 \rrbracket$.}}\label{tab:CPU_3D_helm}
 \vspace{-0.5cm}
\end{table}
}
\subsubsection{High frequency behavior}\label{ssec:freq}
\begin{figure}[H]
    \centering
    \begin{subfigure}[b]{0.32\columnwidth}
        \centering
       \includegraphics[width=\textwidth]{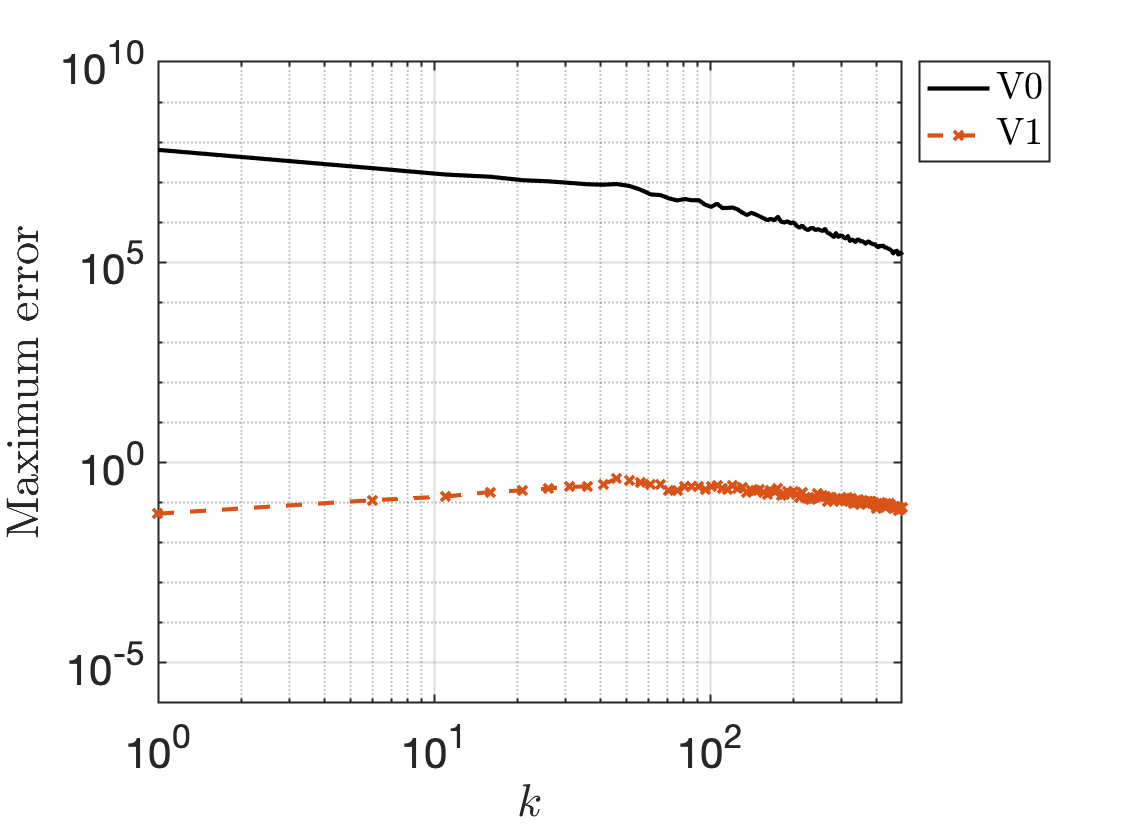}
        \caption{N = 128}
    \end{subfigure}
    \begin{subfigure}[b]{0.32\columnwidth}
        \centering
        \includegraphics[width=\textwidth]{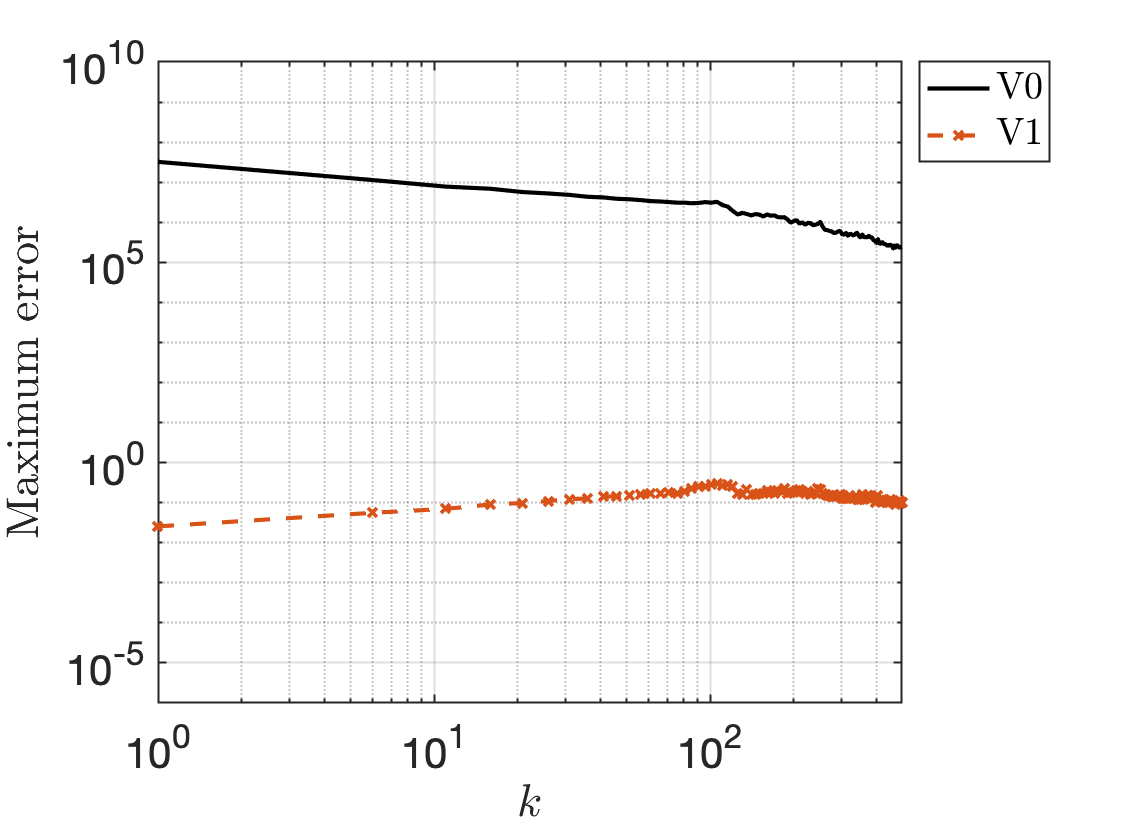}
        \caption{N = 256}
    \end{subfigure}
        \begin{subfigure}[b]{0.32\columnwidth}
        \centering
        \includegraphics[width=\textwidth]{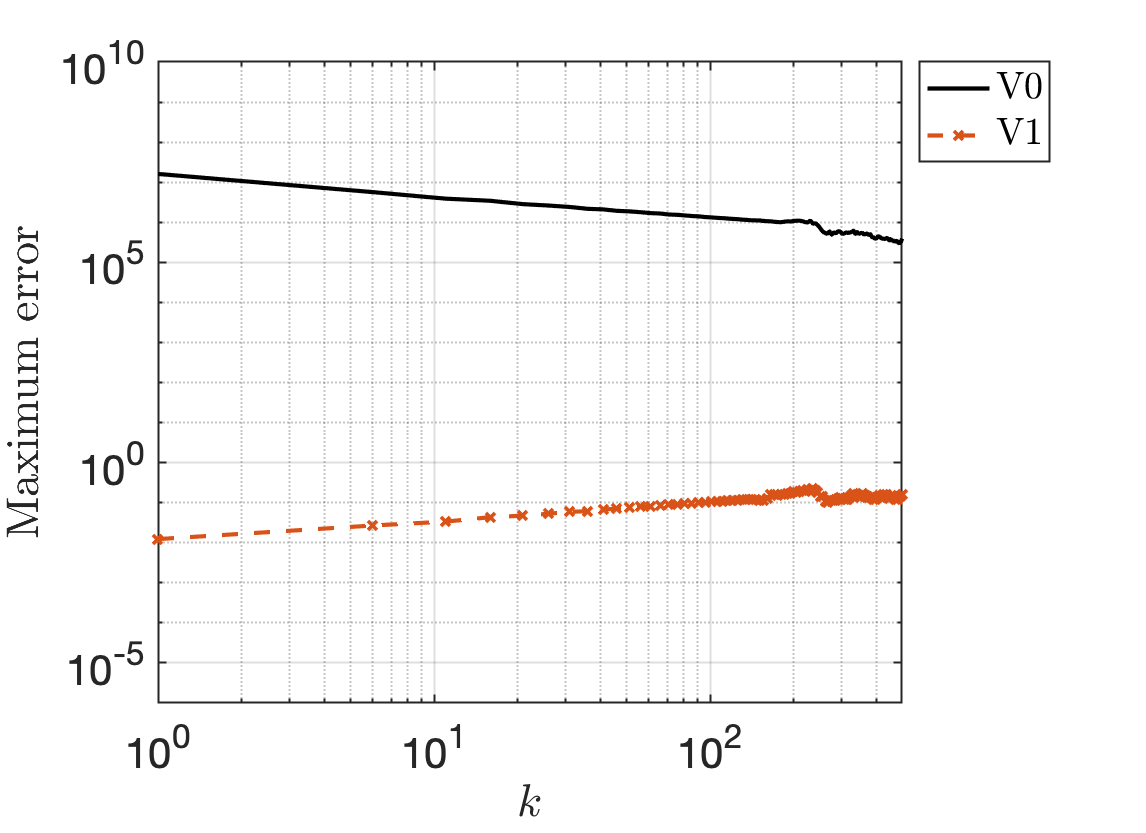}
        \caption{N = 512}
    \end{subfigure}
  \caption{\small {\textbf{Helmholtz 2D.}} {Log-Log of the maximum error in computing the solution of Problem \eqref{eq:scattering_pb} as described in Section \ref{ssec:helm2d}, with respect to the wavenumber $k$, for various number of quadrature points $N$.}}
  \label{img:helmholtz_kcompa_2D}
\end{figure}

\begin{figure}[!hbt]
    \centering
    \begin{subfigure}[b]{0.32\columnwidth}
        \centering
       \includegraphics[width=\textwidth]{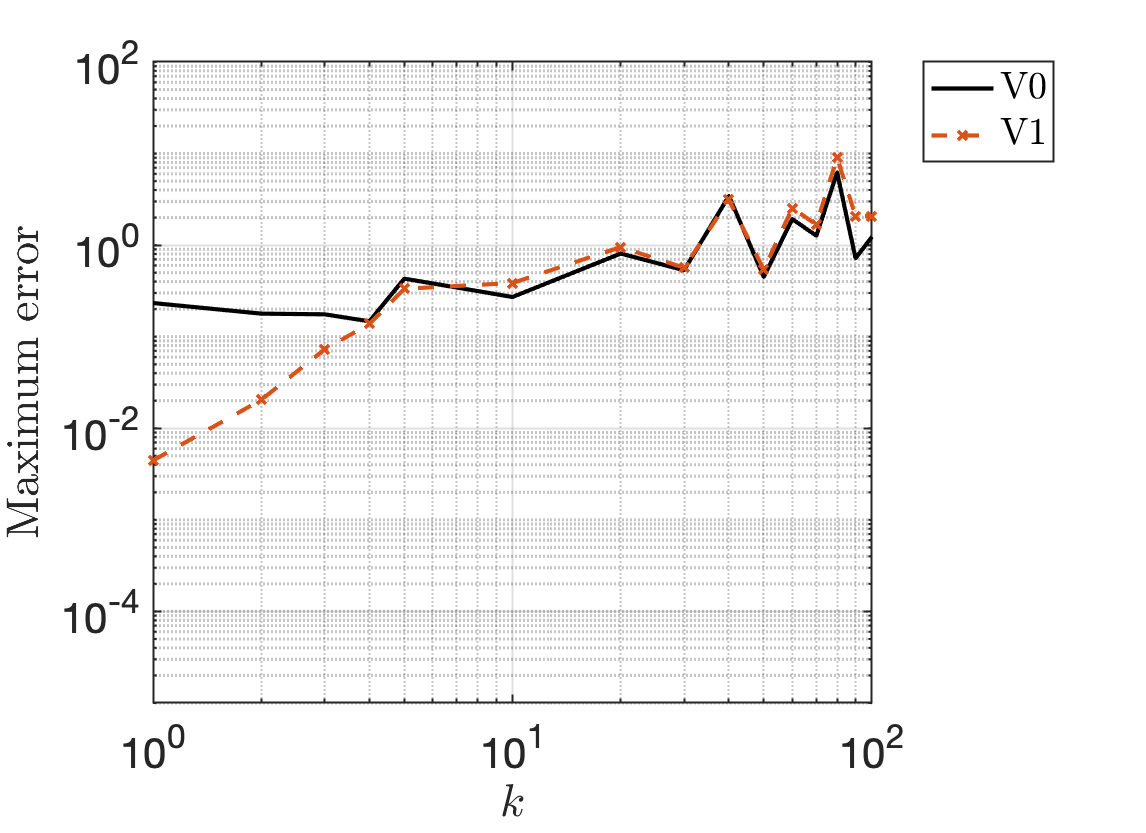}
        \caption{{N = 8}}
    \end{subfigure}
    \begin{subfigure}[b]{0.32\columnwidth}
         \centering
        \includegraphics[width=\textwidth]{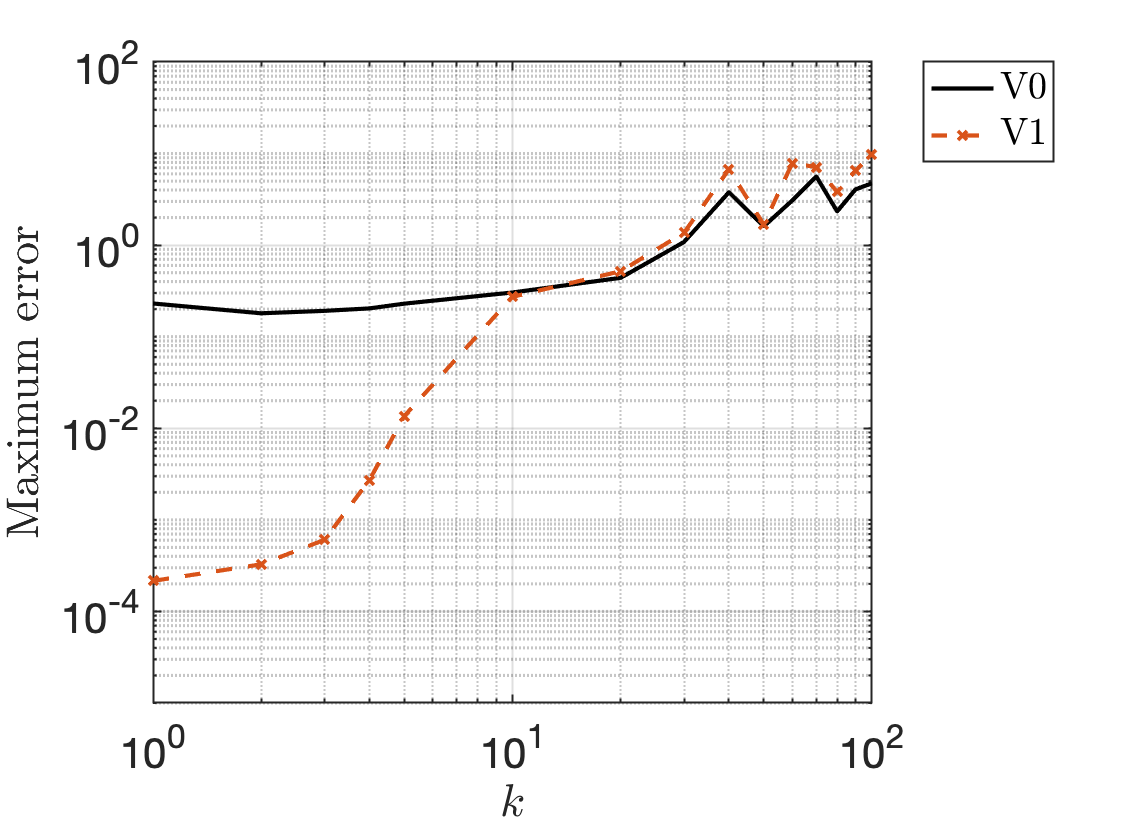}
        \caption{{N=16}}
    \end{subfigure}
        \begin{subfigure}[b]{0.32\columnwidth}
        \centering
        \includegraphics[width=\textwidth]{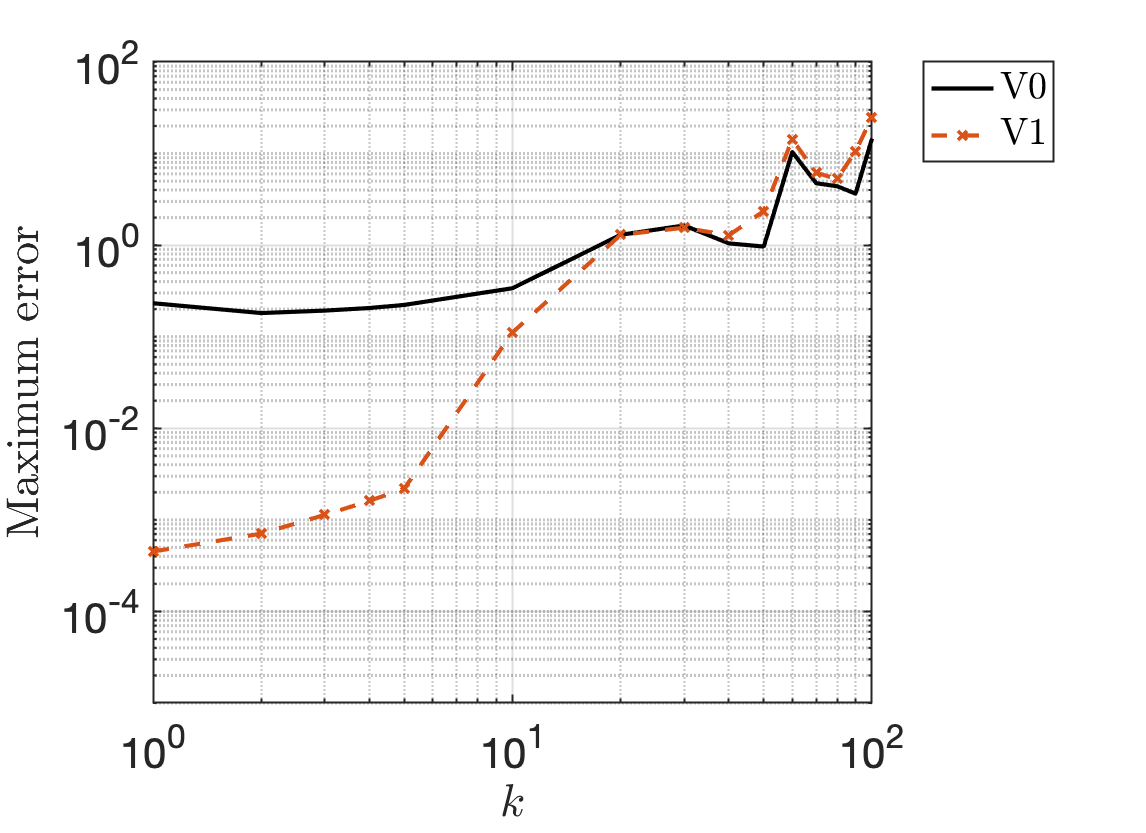}
        \caption{{N=24}}
    \end{subfigure}
  \caption{\small {\textbf{Helmholtz 3D.}} {Log-Log of the maximum error in computing the solution of Problem \eqref{eq:scattering_pb} as described in Section \ref{ssec:helm3d}, with respect to the wavenumber $k$, for various number of quadrature points $N$.}}
  \label{img:helmholtz_kcompa_3D}
\end{figure}

It is well-known that for a fixed number of quadrature points $N$, accuracy is lost for larger wavenumbers $k$. Figures \ref{img:helmholtz_kcompa_2D} and \ref{img:helmholtz_kcompa_3D} represent the high frequency behavior for the Examples 3 and 4, for various $k$ and $N$. We consider the same quadrature rules, exact solution $u_{\text{exact}}$, boundary shapes, as in Sections \ref{ssec:helm2d}, \ref{ssec:helm3d}, but we vary $k$ and/or $N$. The modified representation annihilates some oscillatory behavior by subtracting plane waves along the normal of the evaluation points. It allows then a better approximation for a wider range of wavenumbers (until the number of quadrature points isn't enough), and results in a greater wavenumber stability. Results in Figure \ref{img:helmholtz_kcompa_2D} and \ref{img:helmholtz_kcompa_3D} confirm this phenomenon.

\section{Modified boundary integral equations}\label{S:BIE}

We have used {\eqref{eq:repres_form}} to modify the representation of solution of boundary value problems close to (but not on) the boundary. {One could also use \eqref{eq:repres_form} to avoid weakly }singular integrals in the boundary integral equation as done in BRIEF~\cite{SKKC15}. In the section we present a modified representation of \eqref{eq:BIE_helm}.

\begin{Proposition}\label{pro:BIE_helm}
Given $x^\ast \in \partial D$, let $\vv$ be a solution of Helmholtz equation in $D \subset \mathbb{R}^d$, $d = 2,3$, satisfying conditions \eqref{eq:hyp_helm}.
Then the boundary integral equation \eqref{eq:BIE_helm} admits the modified representation:
\begin{equation}\label{eq:general_BIE_helm}
\begin{aligned}
 & \dsp \int_{\partial D}   \left[ \partial_{n_y}G^H(x^\ast,y)  - \dvv(y) G^H(x^\ast,y) \right] \left[ \mu(y) - \mu (x^\ast) \right]  \, d\sigma_y +  \dsp \int_{\partial D} G^H(x^\ast,y)\left[\dvv(y) - ik \right]  \mu(y)    \, d\sigma_y  \\
 &+ \mu(x^\ast) \dsp \int_{\partial D}  \partial_{n_y} G^H(x^\ast,y) \left[1 - \vv(y)\right]  \, d\sigma_y = f(x^\ast), \quad {\forall x^\ast  \in \partial D}.
\end{aligned}
\end{equation}
The modified representation \eqref{eq:general_BIE_helm} {has smoother integrands than \eqref{eq:BIE_helm}}.
\end{Proposition}
The proof can be found in Appendix \ref{ssec:appdx_helm}. Using again $\vv (y) = e^{ ik n_{x^\ast} \cdot (y-x^\ast)}$, Proposition \ref{pro:BIE_helm} gives us the modified boundary integral equation:
 \begin{equation}\label{eq:modified_BIE_helm}
\begin{aligned}
&  \dsp \int_{\partial D}  \left[ \partial_{n_y}G^H(x^\ast,y) - ik (n_y \cdot n_{x^\ast}) \right. \left. e^{ ik n_{x^\ast} \cdot (y-x^\ast)}G^H(x^\ast,y) \right] \left[ \mu(y) - \mu (x^\ast) \right]  \, d\sigma_y \\
&+ ik  \dsp \int_{\partial D} G^H(x^\ast,y)\left[(n_y \cdot n_{x^\ast}) e^{ ik n_{x^\ast}\cdot (y - x^\ast)}  - 1 \right]  \mu(y)    \, d\sigma_y  \\
&+ \mu(x^\ast) \dsp \int_{\partial D}  \partial_{n_y} G^H(x^\ast,y) \left[1 - e^{ ik n_{x^\ast} \cdot (y-x^\ast)}\right]  \, d\sigma_y  = f(x^\ast), \quad x^\ast \in \partial D.
\end{aligned}
\end{equation}
Equation \eqref{eq:modified_BIE_helm} has no singular integrals {(in the sense its integrands have vanishing singularities)}, in particular it could be approximated using standard quadrature rules such as PTR in two dimensions.
Going back to Examples 3 and 4 presented in Sections \ref{ssec:helm2d} and \ref{ssec:helm3d}, we now compare the approximation of {the representations  \eqref{eq:DLP-SLP_helm}-\eqref{eq:modified_DLP-SLP_helm} where the density $\mu$ has been computed via \eqref{eq:BIE_helm}-\eqref{eq:modified_BIE_helm}. We then have four representations:
\begin{itemize}\renewcommand{\labelitemi}{$\bullet$}
\item \textbf{V0:} standard representation \eqref{eq:DLP-SLP_helm} with previous approximation of \eqref{eq:BIE_helm};
\item \textbf{V1:} modified representation \eqref{eq:modified_DLP-SLP_helm} with previous approximation of \eqref{eq:BIE_helm};
\item \textbf{V2}: standard representation \eqref{eq:DLP-SLP_helm}, approximation of \eqref{eq:modified_BIE_helm} using PTR as Nyström method (2D), using product Gaussian quadrature rule (3D).
\item \textbf{V3}: modified representation \eqref{eq:modified_DLP-SLP_helm}, approximation of \eqref{eq:modified_BIE_helm} using PTR as Nyström method (2D), using product Gaussian quadrature rule (3D).
\end{itemize}
\vspace{-0.5cm}
\begin{figure}[H]
  \centering
        \includegraphics[width=0.32\textwidth]{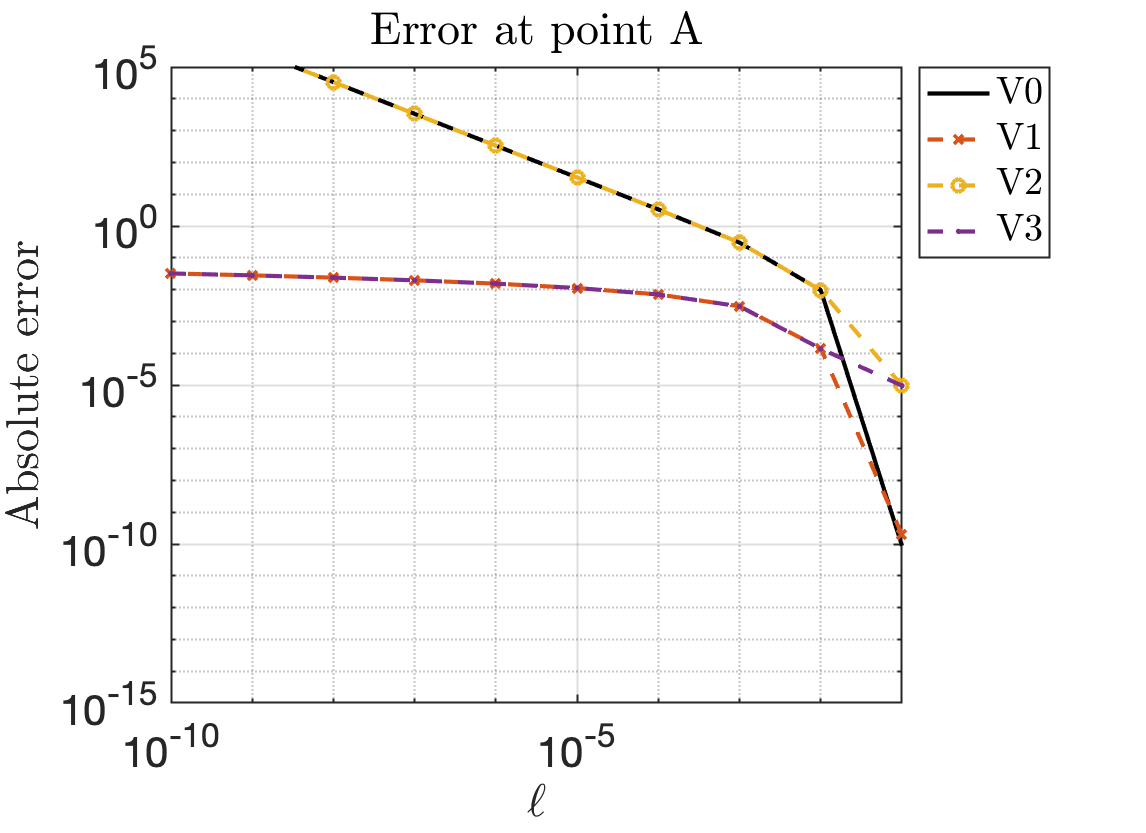}
          \includegraphics[width=0.32\textwidth]{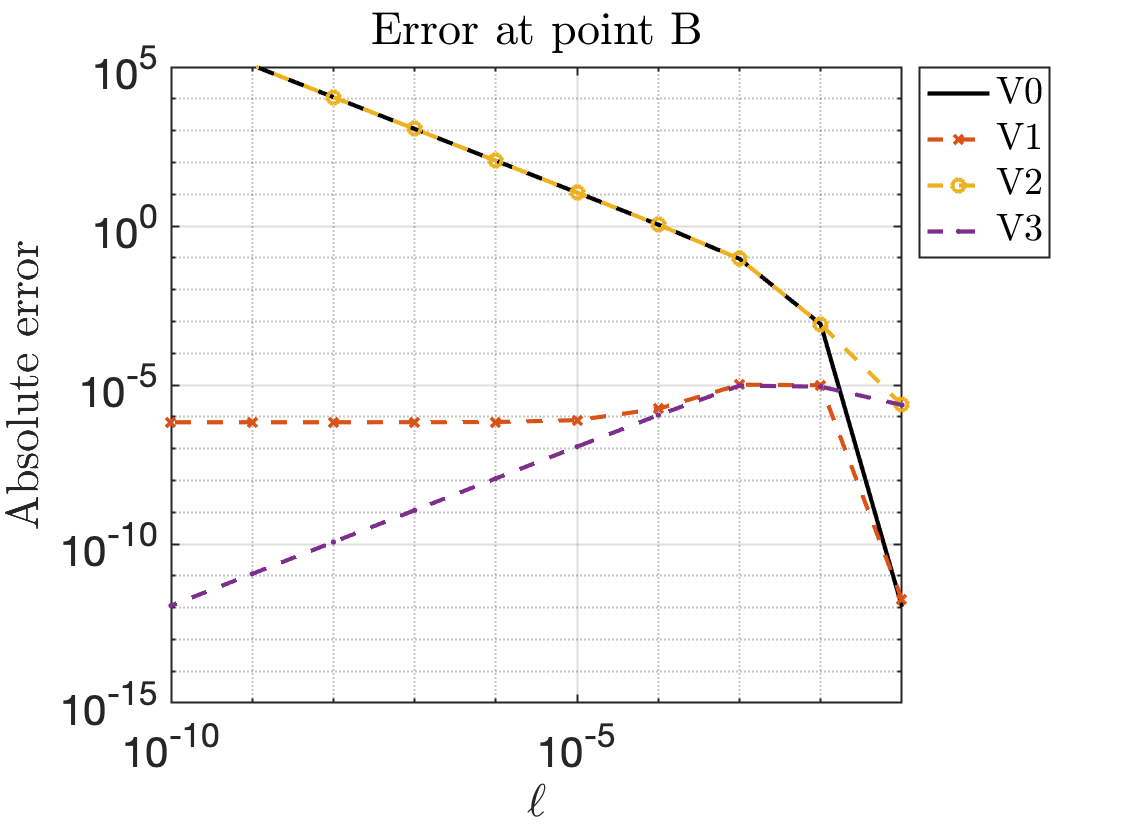}
                 \includegraphics[width=0.32\textwidth]{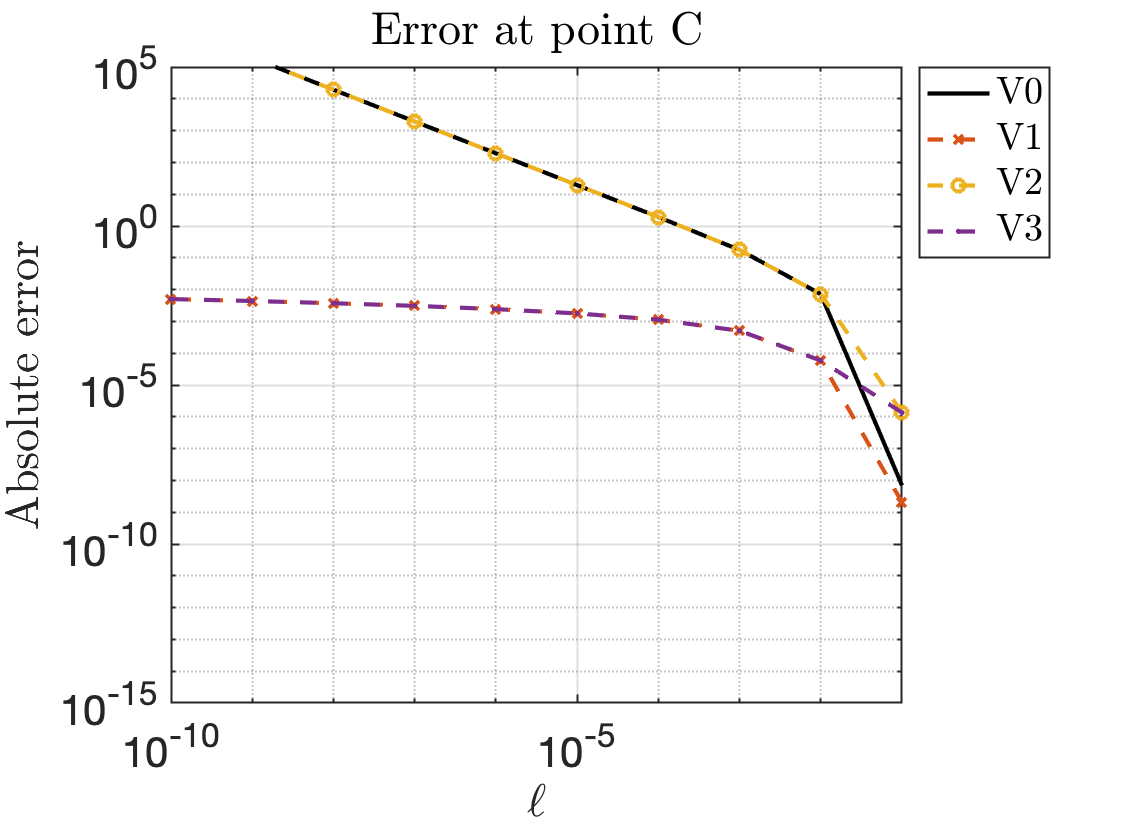}
  \caption{\small {\textbf{Helmholtz 2D.}} Log-Log plot of the error along the normal for the solution of \eqref{eq:scattering_pb} out of the star domain defined by the boundary $y(t) =(1.55+ 0.4 \cos 5t) * (\cos t, \sin t)$, $t \in [0 ,2\pi]$, for the Dirichlet data, $f(x^\ast) = \frac{i}{4} H^{(1)}_0(15 |x^\ast - x_0|)$ with $x_0=(0.2,0.8)$, at the three points A,B, C plotted as black $\times$'s in Figure \ref{img:helmholtz2D_1_1}.}
  \label{img:helmholtz2D_BIE}
\end{figure}
\begin{table}[h!] 
  \centering
 {\small 
\begin{tabular}{|c|c|c|c|}
\hline
Method &  N = 128 & N  = 256 & N = 512\\ 
 \hline
 \eqref{eq:BIE_helm} with Kress product rule & 0.12 & 0.45 & 1.70\\ 
 \hline
 \eqref{eq:modified_BIE_helm} with PTR & 0.09 & 0.302 & 1.16 \\
 \hline
\end{tabular}
}
 \caption{{\textbf{Helmholtz 2D.} CPU times (in seconds) for various number of quadrature points to compute the solution of the boundary integral equation.}}\label{tab:CPU_2D_helm-bie}
 \vspace{-0.5cm}
\end{table}
}
Figure \ref{img:helmholtz2D_BIE} represents the results in two dimensions and illustrates how the resolution of $\mu$ limits the approximation of the solution of \eqref{eq:scattering_pb}. 
Far from the boundary the error made using {V2-V3} cannot be better than order $10^{-6}$. This limitation is due to the poor resolution of $\mu$ using {Nyström method based on PTR to approximate \eqref{eq:modified_DLP-SLP_helm}}. {This can be assessed by looking at the density Fourier coefficients' decay, which caps at $10^{-6}$ for $N = 256$}. However, as the evaluation point gets closer to the boundary ($\ell \to 0$), {V3} yields competitive \edit{(sometimes better)} results. {Additionally, the use of Nyström PTR allows to reduce CPU times as indicated in Table \ref{tab:CPU_2D_helm-bie}.} The modified boundary integral equation \eqref{eq:modified_BIE_helm} can be approximated using standard quadrature rules such as Periodic Trapezoid Rule (note that {Nyström PTR} was not possible to use to solve for \eqref{eq:BIE_helm} due to singular integrals). Its resolution may be limited but it offers interesting corrections for the close evaluation problem using simple quadrature rules {as well as faster solvers.}

\begin{figure}[H]
  \centering
           \includegraphics[width=0.4\textwidth]{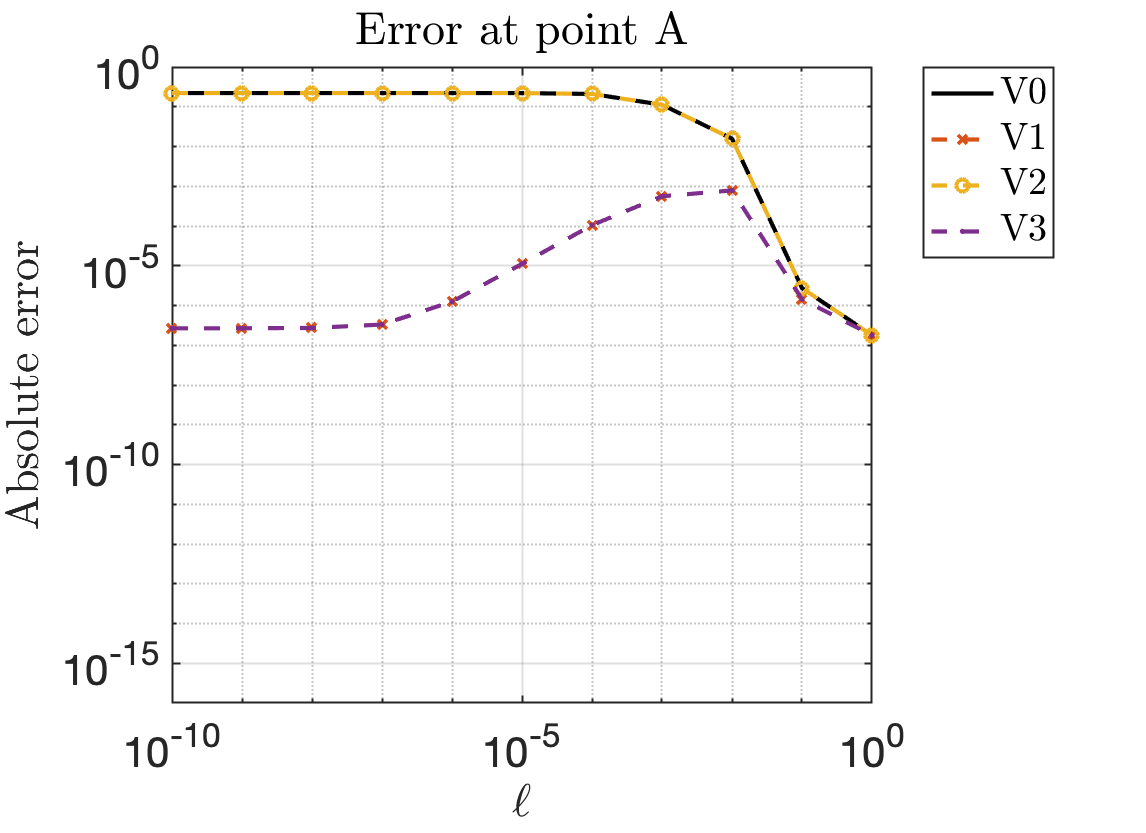}
                    \includegraphics[width=0.4\textwidth]{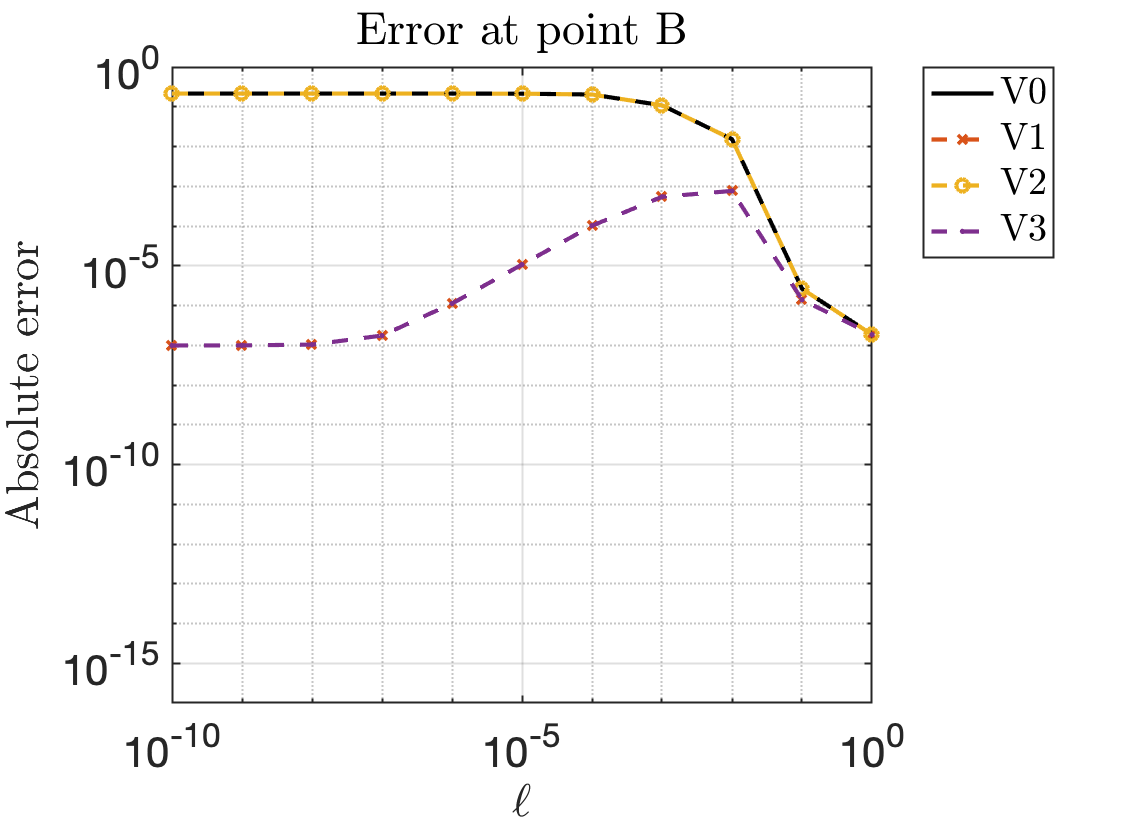}
  \caption{\small {\textbf{Helmholtz 3D.}} Log-Log plot of the error {for the problem described in Figure \ref{img:helmholtz3D_1} using $N = 32$, and for the four representations (standard or modified, off and on boundary).}}
  \label{img:helmholtz3D_BIE_1}
\end{figure}

Results in Figure \ref{img:helmholtz3D_BIE_1} \edit{show that the resolution of the solution using both methods yields the same accuracy in three dimensions}. The product Gaussian quadrature rule is an open quadrature at the singular point $y = x^\ast$ (see Appendix \ref{sec:appdx3D})\edit{. Thus,} the modification introduced in \eqref{eq:modified_DLP-SLP_helm} doesn't affect the approximation. The product Gaussian quadrature rule is a well-used, efficient, easy to implement method, but one could consider a closed quadrature rule to study the effect of \eqref{eq:modified_BIE_helm} more closely.

\section{Conclusion}\label{S:conclu}

In this paper we have provided modified representations for Laplace and Helmholtz layer potentials to address the close evaluation problem in several boundary value problems.  Similar to Gauss' law, we take advantage of one auxiliary function, solution of the partial differential equation at stake. Similar technique has been used in the context of BRIEF and density interpolation\edit{. O}ur approach provides guidelines on how to develop them {independently of the density}, and valuable insights into the layer potentials inherent nearly singular behavior. Several examples in two and three dimensions have been presented and demonstrated the efficiency of the modified representations. Given a quadrature rule, the modified representation of the solution provides a better approximation by several orders of magnitude even with limited computational resources. This assumes that the density, solution of the boundary integral equation, is sufficiently well-resolved. The modified boundary integral equation has no singular behaviors anymore, and allows us to use standard quadrature rules that do not treat singularities. \\
We have provided general modified representations, one can use them with any solution of their choice as long as they follow the provided guidelines to address the close evaluation. One can use this technique to modify any other wave problems, including sound-hard, penetrable obstacles. Future work includes applying those techniques to plasmonic scattering problems \cite{Helsing,ammari}, deriving an asymptotic analysis to quantify the limit behavior of the error as the evaluation point approaches the boundary, as well as extensions to other partial differential equations such as Stokes problems and others.


\appendix
\section{Kress product quadrature}\label{sec:appdx2D}
In this section we provide a brief summary about the Kress product quadrature rule~\cite{Kress91} used to compute the density $\mu$, solution of \eqref{eq:BIE_helm}, in two dimensions. Denoting the parameterization of $\partial D$ as $y(t)$, $t \in (0, 2 \pi)$, and denoting $x^\ast  = y(t^\ast)$, we compactly rewrite \eqref{eq:BIE_helm}
\begin{equation}
  \frac{1}{2} \mu(t^\ast) + \dsp \int_0^{2 \pi} {K}(t,t^\ast) \mu(t) \, dt = f(t^\ast),
  \label{eq:BIE-mu0}
\end{equation}
with the abuse of notation ${K}(t,t^\ast) = \left( \partial_{n_y} G^H(x^\ast, y(t)) - ik G^H(x^\ast,y(t))\right) |y'(t)|$, $\mu(t) = \mu(y(t))$, and $f(t) = f(y(t))$. The Kress product quadrature rule is well adapted for weakly singular integrals involving kernel with a logarithmic singularity. To that aim one rewrites:
\[{K}(t,t^\ast) = {K}_1(t,t^\ast) \log \left(4 \sin ^2\left( \frac{t^\ast - t}{2}\right) \right) + {K}_2(t,t^\ast), \]
with smooth functions $K_1$, $K_2$ (the expression of $K_1$, $K_2$ can be found in \cite{Kress91}). Then one discretizes the integral using $N = 2n$ quadrature points as follows:
\[
\begin{aligned}
& \dsp \int_0^{2 \pi} {K}(t,t^\ast) \mu(t) \, dt  \approx \sum \limits_{k=0}^{2n-1}   \left( R_k^{(n)}(t^\ast) K_1(t^\ast, t_k) + \frac{\pi}{n} K_2(t^\ast, t_k) \right) \mu(t_k),
\end{aligned}\]
with $t_k = \frac{\pi k}{n}$, $k = 0, \dots, 2n-1$, and $R_k^{(n)}(t^\ast)$ the weights 
\[ R_k^{(n)}(t^\ast) = \dsp -\frac{2\pi}{n} \sum \limits_{j=1}^{n-1} \frac{1}{j} \cos(j(t^\ast - t_k)) - \frac{\pi}{n^2} \cos \left(n(t^\ast -t_k) \right), \quad k = 0, \dots, 2n-1.\] 
\section{Galerkin approximation}\label{sec:appdx3D}
In this section we provide a brief summary about the Galerkin approximation used to compute the solutions of \eqref{eq:BIE_laplace}, \eqref{eq:BIE_laplace_NS}, \eqref{eq:BIE_helm} and \eqref{eq:modified_BIE_helm} in three dimensions. First, we compactly write \eqref{eq:BIE_laplace}, \eqref{eq:BIE_laplace_NS}, \eqref{eq:BIE_helm} and \eqref{eq:modified_BIE_helm} as
\begin{equation}
  \mathscr{K} [\psi] = F ,
  \label{eq:BIE-mu}
\end{equation}
with $\psi$ denoting the density (i.e. $\mu$, $\rho$), and $F$ denoting the Dirichlet or Neumann data. We introduce the approximation for $\psi$
\begin{equation}
  \psi(y(\theta,\varphi))  \approx \sum_{n
  = 0}^{N - 1} \sum_{m = -n}^{n} Y_{nm}(\theta,\varphi) \hat{\psi}_{nm},
  \label{eq:mu-expansion}
\end{equation}
with $y(\theta, \varphi)$, $\theta \in (0, \pi)$, $\varphi \in (-\pi,\pi)$ a parameterization of the boundary $\partial D$, $\{ Y_{nm}(\theta,\varphi) \}_{n,m}$ the orthonormal set of
spherical harmonics. For $x^\ast \in \partial D$, we write $x^\ast = y(\theta^\star, \varphi^\star)$. Note that $N$ in \eqref{eq:mu-expansion}
corresponds also to the same order of the quadrature rule used to approximate \eqref{eq:BIE_laplace}, \eqref{eq:BIE_laplace_NS}, \eqref{eq:BIE_helm} and \eqref{eq:modified_BIE_helm}. Substituting \eqref{eq:mu-expansion} into \eqref{eq:BIE-mu} and taking
the inner product with $Y_{n'm'}(\theta^{\star},\varphi^{\star})$, we
obtain the Galerkin equations
\begin{equation}
\sum_{n = 0}^{N-1} \sum_{m =
    -n}^{n} \langle Y_{n'm'}, \mathscr{K}[Y_{nm}] \rangle
  \hat{\psi}_{nm} = \langle Y_{n'm'}, F \rangle.
  \label{eq:DLP-Galerkin}
\end{equation}

We construct the $N^{2} \times N^{2}$ linear system for the unknown
coefficients, $\hat{\psi}_{n'm'}$ resulting from \eqref{eq:DLP-Galerkin}
evaluated for $n' = 0, \cdots, N - 1$ with corresponding values of
$m'$. To compute the inner products,
$\langle Y_{n'm'}, \mathscr{K}[Y_{nm}] \rangle$ and
$\langle Y_{n'm'}, F \rangle$, we use the product Gaussian quadrature
rule for spherical integrals~\cite{atkinson1982numerical}. This corresponds to approximate the integral with respect to $\varphi$ using $N$
Gauss-Legendre quadrature points, and the integral with respect to $\theta$ using a $2N$ Periodic Trapezoid Rule
points. One can proceed as in the three-step method (see Section \ref{ssec:3D_laplace}, and \cite{KKCC19} for more details), by adding a rotation of the local coordinate system so that $x^\ast$ corresponds to the north pole, and by using the $N$
Gauss-Legendre quadrature points mapped to $(0,\pi)$ and not $(-1,1)$.\\
For \eqref{eq:BIE_laplace} we have 
\[ \begin{aligned}
&  \mathscr{K} [Y_{nm}](\theta^{\star},\varphi^{\star}) =  - \frac{1}{2} Y_{nm}(\theta^{\star},\varphi^{\star})  + \int_{-\pi}^{\pi} \int_{0}^{\pi}  \partial_{n_y} G^L(\theta^{\star},\varphi^{\star},\theta, \varphi) 
    J(\theta,\varphi) \sin (\theta) 
  Y_{nm}(\theta,\varphi) \mathrm{d}\theta
  \mathrm{d}\varphi.
  \end{aligned}
  \]
For \eqref{eq:BIE_laplace_NS} we make use of the adjoint $\mathscr{K}^\star$ of $\mathscr{K}$. Using Gauss' law we write \\
$\sum_{n = 0}^{N-1} \sum_{m =
    -n}^{n} \langle  \mathscr{K}^\star [Y_{n'm'}], Y_{nm} \rangle
  \hat{\psi}_{nm} = \langle Y_{n'm'}, F \rangle$
with 
\[ \begin{aligned}
& \mathscr{K}^\star [Y_{n'm'}](\theta,\varphi) =  \int_{-\pi}^{\pi} \int_{0}^{\pi}  \partial_{n_x^\ast} G^L(\theta^{\star} ,\varphi^{\star},\theta, \varphi)
    J(\theta^\star,\varphi^\star) \sin (\theta^\star) 
 [ Y_{n'm'}(\theta^\star,\varphi^\star) -Y_{n'm'}(\theta,\varphi) ] \mathrm{d}\theta^\star
  \mathrm{d}\varphi^\star.
  \end{aligned}
  \]
  For \eqref{eq:BIE_helm} we have
\[ \begin{aligned}  & \mathscr{K} [Y_{nm}](\theta^{\star},\varphi^{\star}) =  \frac{1}{2} Y_{nm}(\theta^{\star},\varphi^{\star}) +\int_{-\pi}^{\pi} \int_{0}^{\pi} \left[ \partial_{n_y} G^H(\theta^{\star},\varphi^{\star},\theta, \varphi) - ik G^H(\theta^{\star} ,\varphi^{\star},\theta, \varphi)\right] \\
&  \hspace{10cm} \quad 
    J(\theta,\varphi) \sin (\theta) 
  Y_{nm}(\theta,\varphi) \mathrm{d}\theta
  \mathrm{d}\varphi,
  \end{aligned}
  \]
and for \eqref{eq:modified_BIE_helm} we have
\[
\begin{aligned}
 \mathscr{K}_m [Y_{nm}] & (\theta^{\star},\varphi^{\star}) =   \int_{-\pi}^{\pi} \int_{0}^{\pi} \left[ \partial_{n_y} G^H(\theta^{\star},\varphi^{\star},\theta, \varphi) - ik (n_y \cdot n_{x^\ast})\right.  \left.  e^{ik (n_{x^\ast}\cdot (y(\theta, \varphi) - y(\theta^\star, \varphi^\star))} G^H(\theta^{\star}, \varphi^{\star},\theta, \varphi)\right] \\
 &   \hspace{6.8cm} \quad  J(\theta,\varphi) \sin (\theta) 
 \left[ Y_{nm}(\theta,\varphi) -Y_{nm}(\theta^\star,\varphi^\star) \right] \mathrm{d}\theta
  \mathrm{d}\varphi \\
  &+ i k \int_{-\pi}^{\pi} \int_{0}^{\pi}  [( n_y \cdot n_{x^\ast})e^{ik (n_{x^\ast} \cdot (y(\theta, \varphi) - y(\theta^\star, \varphi^\star))} -1]  G^H(\theta^{\star} \varphi^{\star},\theta, \varphi)
    J(\theta,\varphi) \sin (\theta) 
  Y_{nm}(\theta,\varphi) \mathrm{d}\theta
  \mathrm{d}\varphi\\
  & +   Y_{nm}(\theta^\star,\varphi^\star) \int_{-\pi}^{\pi} \int_{0}^{\pi}  [1 - e^{ik (n_{x^\ast} \cdot (y(\theta, \varphi) - y(\theta^\star, \varphi^\star))}] \partial_{n_y}G^H(\theta^{\star} \varphi^{\star},\theta, \varphi)
    J(\theta,\varphi) \sin (\theta) 
 \mathrm{d}\theta
  \mathrm{d}\varphi .
  \end{aligned}
  \]

\section{Proof of modified representations}\label{sec:appdx_proofs}

\subsection{Modified double-layer potential \eqref{eq:general_DLP}}\label{ssec:appdx_dlp}
Given $\vv$ solution of Laplace's equation in $D \subset \mathbb{R}^d$, $d =2,3$, and for $x \in D$ we write $x = x^\ast - \ell n_{x^\ast}$, with $x^\ast \in \partial D$. Then we write \eqref{eq:DLP_laplace} as:
\[
\begin{aligned}
 u(x) &  = \dsp \int_{\partial D} \partial_{n_y} G(x,y) \mu(y)\left[ 1 - \vv(y)\right]  \, d\sigma_y +\dsp \int_{\partial D} \partial_{n_y} G(x,y) \mu(y) \vv(y)  \, d\sigma_y  \\
& = \dsp \int_{\partial D} \partial_{n_y} G(x,y) \mu(y)\left[ 1 - \vv(y)\right]  \, d\sigma_y +\dsp \int_{\partial D} \partial_{n_y} G(x,y) [ \mu(y) - \mu(x^\ast) ] \vv(y)  \, d\sigma_y   \\
&+ \mu(x^\ast) \dsp \int_{\partial D} \partial_{n_y}G(x,y)   \vv(y) -  G(x,y)  \dvv(y)   \, d\sigma_y + \mu(x^\ast) \dsp \int_{\partial D}   G(x,y) \dvv(y)   \, d\sigma_y\\
\end{aligned}
\]
Using \eqref{eq:repres_form} the third term becomes $- \mu(x^\ast) \vv(x^\ast)$. Then 
\[
\begin{aligned}
u(x)  & =  \dsp \int_{\partial D} \partial_{n_y} G(x,y) \mu(y)\left[ 1 - \vv(y)\right]  \, d\sigma_y +\dsp \int_{\partial D} \partial_{n_y} G(x,y) [ \mu(y) - \mu(x^\ast) ] \vv(y)  \, d\sigma_y   - \mu(x^\ast) \vv(x^\ast) \\
&+ \mu(x^\ast) \dsp \int_{\partial D}   G(x,y)  [ \dvv(y) - \dvvx(x^\ast)]   \, d\sigma_y  + \mu(x^\ast) \dvvx(x^\ast) \dsp \int_{\partial D} G(x,y)  \, d\sigma_y
\end{aligned}
\]
which is \eqref{eq:general_DLP}, after using \eqref{eq:repres_form} for the last term.
\subsection{Proof of Proposition \ref{pro:SLP}}\label{ssec:appdx_slp}
In this section we derive \eqref{eq:general_SLP}. Given $\vv$ solution of Laplace's equation in $D \subset \mathbb{R}^d$, $d =2,3$, and for $x \in E$ we write $x = x^\ast +\ell n_{x^\ast}$, with $x^\ast \in \partial D$. Then we write \eqref{eq:SLP_laplaceN} as:
\[
\begin{aligned}
u(x)  &  = \dsp \int_{\partial D}  G(x,y) \rho(y)\left[ 1 - \partial_{n_y} \vv(y)\right]  \, d\sigma_y +\dsp \int_{\partial D} G(x,y) \rho(y) \dvv(y)  \, d\sigma_y  \\
& = \dsp \int_{\partial D} G(x,y) \rho(y)\left[ 1 - \dvv(y)\right]  \, d\sigma_y +\dsp \int_{\partial D} G(x,y) [ \rho(y) - \rho(x^\ast) ] \dvv(y)  \, d\sigma_y   \\
& + \rho(x^\ast) \dsp \int_{\partial D}  G(x,y) \dvv(y) - \partial_{n_y} G(x,y)  \vv(y)   \, d\sigma_y + \rho(x^\ast) \dsp \int_{\partial D} \partial_{n_y}G(x,y) \vv(y)   \, d\sigma_y\\
\end{aligned}
\]
Using \eqref{eq:repres_form}, the third term vanishes. Then 
\[
\begin{aligned}
 u(x)  &
 = \dsp \int_{\partial D} G(x,y) \rho(y)\left[ 1 - \dvv(y)\right]  \, d\sigma_y +\dsp \int_{\partial D} G(x,y) [ \rho(y) - \rho(x^\ast) ] \dvv(y)  \, d\sigma_y\\
 & +\dsp \rho(x^\ast) \int_{\partial D}  \partial_{n_y} G(x,y) [\vv(y) - \vv(x^\ast) ] \, d\sigma_y  +\dsp \rho(x^\ast) \vv(x^\ast) \int_{\partial D} \partial_{n_y} G(x,y)   \, d\sigma_y   
\end{aligned}
\]
The last term vanishes using \eqref{eq:repres_form} then one obtains \eqref{eq:general_SLP}.
\subsection{Proof of Propositions \ref{pro:helm}, \ref{pro:BIE_helm}}\label{ssec:appdx_helm}
In this section we derive \eqref{eq:general_helm}, \eqref{eq:general_BIE_helm}. Given $\vv$ solution of the Helmholtz equation in $D \subset \mathbb{R}^d$, $d =2,3$, and for $x \in E$ we write $x = x^\ast +\ell n_{x^\ast}$, with $x^\ast \in \partial D$. Then we write \eqref{eq:DLP-SLP_helm} as:
\begin{equation} \label{eq:deriv_helm}
\begin{aligned}
& u(x)    = \dsp \int_{\partial D}  \left[ \partial_{n_y}G^H(x,y) - \dvv(y) G^H(x,y) \right]  \mu(y)  \, d\sigma_y  + \dsp \int_{\partial D}  \left[ \dvv(y)  - ik  \right] G^H(x,y) \mu(y)  \, d\sigma_y \\
& = \dsp \int_{\partial D}  \left[ \partial_{n_y}G^H(x,y) - \dvv(y) G^H(x,y) \right]  [\mu(y)- \mu (x^\ast)]  \, d\sigma_y \\
& + \mu (x^\ast) \dsp \int_{\partial D}  \left[ \partial_{n_y}G^H(x,y) \vv(y) -  G^H(x,y) \dvv(y) \right]  \, d\sigma_y\\
&+ \dsp \int_{\partial D}  \left[ \dvv(y)  - ik  \right] G^H(x,y) \mu(y)  \, d\sigma_y + \mu (x^\ast) \dsp \int_{\partial D} \partial_{n_y}G^H(x,y) \left[  (1 - \vv(y) \right]  \, d\sigma_y 
\end{aligned}
\end{equation}
Using \eqref{eq:repres_form}, the third term vanishes, then one obtains \eqref{eq:general_helm}. One proceeds similarly starting with \eqref{eq:BIE_helm}: one can show that the layer potentials in \eqref{eq:BIE_helm} correspond to \eqref{eq:deriv_helm} for $x = x^\ast \in \partial D$. Finally, \eqref{eq:repres_form} gives that the third term boils down to $- \frac{1}{2} \mu(x^\ast) \vv(x^\ast)$, which finishes the proof.

\textbf{Acknowledgements. }
This research was supported by National Science Foundation Grant: DMS-1819052. The author would like to thank S. Khatri, A. D. Kim, M. Bonnet, and Z. Moitier for fruitful discussions and feedback.


\begin{thebibliography}{}

  	\bibitem{akselrod2014probing} Akselrod~G.~M.,Argyropoulos~C.,
	  Hoang~T~B., Cirac{\`\i}~C., Fang~C., Huang~J., Smith~D.~R.,
	  Mikkelsen~M.~H., Probing the mechanisms of large Purcell
	  enhancement in plasmonic nanoantennas, Nat. Photonics 8 (2014)
	  835--840.
	  
    \bibitem{barnett2015spectrally} Barnett~A.~H., Wu~B., Veerapaneni~S.,
  Spectrally accurate quadratures for evaluation of layer potentials
  close to the boundary for the 2D Stokes and Laplace equations, SIAM
  J. Sci. Comp. 37 (2015) B519--B542.
  
      \bibitem{keaveny2011applying} Keaveny~E.~E., Shelley~M.~J., Applying a
    second-kind boundary integral equation for surface tractions in
    {S}tokes flow, J. Comput. Phys. 230 (2011) 2141--2159.
  
  	  \bibitem{marple2016fast} Marple~G.~R., Barnett~A., Gillman~A.,
  Veerapaneni~S., A fast algorithm for simulating multiphase flows
  through periodic geometries of arbitrary shape, SIAM J.
  Sci. Comput. 38 (2016) B740--B772.
  

	\bibitem{mayer2008label} Mayer~K.~M., Lee~S., Liao~H., Rostro~B.~C.,
  Fuentes~A., Scully~P.~T., Nehl~C.~L., Hafner~J.~H., A label-free
  immunoassay based upon localized surface plasmon resonance of gold
  nanorods, ACS Nano 2 (2008) 687--692.
  
  
  \bibitem{novotny2011antennas} Novotny~L., Van~Hulst~N., Antennas for
    light, Nat. Photonics 5 (2011) 83--90.

  \bibitem{sannomiya2008situ} Sannomiya~T., Hafner~C., Voros~J., In situ
    sensing of single binding events by localized surface plasmon
    resonance, Nano Lett. 8 (2008) 3450--3455.
    
     \bibitem{smith2009boundary} Smith~D.~J., A boundary element
  regularized Stokeslet method applied to cilia-and flagella-driven
  flow, Proc. R. Soc. Lond. A 465 (2009) 3605--3626.
  
 \bibitem{schwab1999extraction}
Schwab~C., Wendland~W., On the extraction technique in boundary integral
  equations, Math. Comput. 68~(225) (1999) 91--122.
  
    \bibitem{beale2001method}
Beale~J.~T., Lai~M.~C., A method for computing nearly singular integrals, SIAM
  J. Numer. Anal. 38~(6) (2001) 1902--1925.

\bibitem{beale2016simple}
J.~T. Beale, W.~Ying, J.~R. Wilson, A simple method for computing singular or
  nearly singular integrals on closed surfaces, Commun. Comput. Phys. 20~(3)
  (2016) 733--753.
  
  \bibitem{helsing2008evaluation}
Helsing~J., Ojala,~R. On the evaluation of layer potentials close to their
  sources, J. Comput. Phys. 227~(5) (2008) 2899--2921.
  
 \bibitem{af2016fast}
af~Klinteberg~ L., Tornberg~A.-K., A fast integral equation method for solid
  particles in viscous flow using quadrature by expansion, J. Comput. Phys. 326
  (2016) 420--445.

\bibitem{af2017error}
af~Klinteberg~L., Tornberg~A.-K., Error estimation for quadrature by expansion
  in layer potential evaluation, Adv. Comput. Math. 43~(1) (2017) 195--234.
  
\bibitem{barnett2014evaluation}
Barnett~A.~H., Evaluation of layer potentials close to the boundary for
  {L}aplace and {H}elmholtz problems on analytic planar domains, SIAM J. Sci.
  Comput. 36~(2) (2014) A427--A451.
  
    \bibitem{epstein2013convergence}
Epstein~C.~L., Greengard~L., Kl{\"o}ckner~A.~K., On the convergence of local
  expansions of layer potentials, SIAM J. Numer. Anal. 51~(5) (2013)
  2660--2679.
  
  \bibitem{klockner2013quadrature}
Kl{\"o}ckner~A., Barnett~A., Greengard~L., O'Neil~M., Quadrature by expansion:
  A new method for the evaluation of layer potentials, J. Comput. Phys. 252
  (2013) 332--349.
  
\bibitem{rachh2017fast}
Rachh~M., Kl{\"o}ckner~A., O'Neil~M., Fast algorithms for quadrature by
  expansion i: Globally valid expansions, J. Comput. Phys. 345 (2017) 706--731.
  
\bibitem{wala20183DQBX}
Wala~M., Kl{\"o}ckner~A., A fast algorithm for {Q}uadrature by {E}xpansion in
  three dimensions, J. Comput. Phys. 388 (2019) 655-689.
  
    \bibitem{GORV21} 
  Greengard~L., O'Neil~M., Rachh~M., Vico~F., Fast multipole methods for the evaluation of layer potentials with locally-corrected quadratures, Journal of Computational Physics: X 10 (2021) 100092.

  \bibitem{Perez18}
Pérez-Arancibia~C., A plane-wave singularity subtraction technique for the classical Dirichlet and Neumann combined field integral equations, Appl. Numer. Math. 123 (2018) 221-240.

\bibitem{PeFaTu19}
Pérez-Arancibia~C., Faria~L., Turc,~C. Harmonic density interpolation methods for high-order evaluation of Laplace layer potentials in 2D and 3D, J. Comput. Phys. 376 (2019) 411-434.

\bibitem{PeTuFa19}
Pérez-Arancibia~C., Turc~C., Faria~L., Planewave density interpolation methods for 3D Helmholtz boundary integral equations, SIAM J. Sci. Comp. 41~(4) (2019) A2088-A2116.
  
\bibitem{carvalho2018asymptotic}
Carvalho~C., Khatri~S., Kim~A.~D., Asymptotic analysis for close evaluation of
  layer potentials, J. Comput. Phys. 355 (2018) 327--341.

\bibitem{ckk2020asymp}
Carvalho~C., Khatri~S., Kim~A.~D., Asymptotic approximation for the close
  evaluation of double-layer potentials, SIAM J. Sci. Comp. 42~(1) (2020) A504-A533.
  
\bibitem{KKCC19}
Khatri~S., Kim~A.~D., Cortes~R., Carvalho~C.,
Close evaluation of layer potentials in three dimensions, J. Comput. Phys, 423 (2020) 109798.    
  
\bibitem{Hwang13} 
 Hwang~W. S., A regularized boundary integral method in potential theory, Computer Methods in Applied Mechanics and Engineering. 259~(123) (2013) 9. 
 
 \bibitem{Liu99}
Liu~Y.~J., Rudolphi~T.~J., New identities for fundamental solutions and their applications to non-singular boundary element formulations, Comp. Mech. 24(1999) 286-292.
   
\bibitem{KSC12}
Klaseboer~E., Sun~Q., Chan~D. Y., Non-singular boundary integral methods for fluid mechanics applications, Journal of Fluid Mechanics. 696~(468) (2012) 78.

\bibitem{SKKC14}
Sun~Q., Klaseboer~E., Khoo~B.-C., Chan~D. Y., A robust and non-singular formulation of the boundary integral method for the potential problem, Engineering Analysis with Boundary Elements, 1~(43) (2014) 117-23.  
  
  \bibitem{SKKC15}
  Sun~Q., Klaseboer~E., Khoo~B.-C., Chan~D. Y., Boundary regularized integral equation formulation of the Helmholtz equation in acoustics, Royal Society open science. 2~(1) (2015) 140520.

  \bibitem{FaPeBo21}
  Faria~L. M., Pérez-Arancibia~C., Bonnet~M., General-purpose kernel regularization of boundary integral equations via density interpolation, Computer Methods in Applied Mechanics and Engineering 378 (2021) 113703.  
  
\bibitem{CoKr13}
Colton~D., Kress~R.,
\textit{Integral equation methods in scattering theory},
SIAM, 2013.

\bibitem{guenther1996partial}
 Guenther~R.~B., Lee~J.~W., Partial Differential Equations of Mathematical
  Physics and Integral Equations, Dover Publications, 1996.
  
\bibitem{atkinson1997numerical}
Atkinson~K.~E., The Numerical Solution of Integral Equations of the Second
  Kind, Cambridge University Press, 1997.
  
  \bibitem{bremer2010nonlinear}
Bremer~J., Gimbutas~Z., Rokhlin~V., A nonlinear optimization procedure for
  generalized gaussian quadratures, SIAM J. Sci. Comp.  32~(4)
  (2010) 1761--1788.

\bibitem{bruno2001fast}
Bruno~O.~P., Kunyansky~L.~A., A fast, high-order algorithm for the solution of
  surface scattering problems: basic implementation, tests, and applications,
  J. Comput. Phys. 169~(1) (2001) 80--110.
  
\bibitem{ganesh2004high}
Ganesh~M., Graham~I., A high-order algorithm for obstacle scattering in three
  dimensions, Journal of Computational Physics 198~(1) (2004) 211--242. 
  
\bibitem{Kress91}
Kress~R., {Boundary integral equations in time-harmonic acoustic scattering}, Math. Comp. Mod. 15 (1991) 229-243.

\bibitem{Kress89}
Kress~R., \textit{Linear Integral Equations}, Springer, 1989.

\bibitem{Carvalho-codes}
Carvalho~C., Subtraction-techniques codes, 2020 {https://doi.org/10.5281/zenodo.3934284
}.
    
\bibitem{atkinson1982numerical}
Atkinson~K.~E., Numerical integration on the sphere, ANZIAM J. 23~(3) (1982)
  332--347.

\bibitem{atkinson1982laplace}
 Atkinson~K.~E., The numerical solution {L}aplace's equation in three
  dimensions, SIAM J. Numer. Anal. 19~(2) (1982) 263--274.

\bibitem{atkinson1985algorithm}
 Atkinson~K.~E., Algorithm 629: An integral equation program for {L}aplace's
  equation in three dimensions, ACM Trans. Math. Softw. 11~(2) (1985) 85--96.

\bibitem{atkinson1990survey}
Atkinson~K.~E., A survey of boundary integral equation methods for the
  numerical solution of {L}aplace's equation in three dimensions, in: Numerical
  {S}olution of {I}ntegral {E}quations, Springer (1990) 1--34.
   
\bibitem{ammari}	  
Ammari~H., Millien~P., Ruiz~M., Zhang~H., Mathematical analysis of plasmonic nanoparticles: the scalar case., Archive for Rational Mechanics and Analysis, 2 (2017) 597-658.

  \bibitem{Helsing}
  Helsing~J., Karlsson~A., An extended charge-current formulation of the electromagnetic transmission problem, SIAM J. App. Math. 80 (2020) 951-976.

\end{thebibliography}
\end{document}